%% file: paper_variational.tex
\documentclass{elsarticle}
\usepackage{amssymb,amsmath,amsthm}
\usepackage{tikz}
\usepackage{pgfplots,pgfplotstable}
\usepackage[hypcap=true]{subcaption}
\usepackage{comment}

\begin{document}

\title{Bounds-constrained finite element approximation of time-dependent partial differential equations\tnoteref{t1}}
\tnotetext[t1]{This work is supported by the U.S. National Science Foundation grant \#2410736.}

\author[1]{Robert C. Kirby\corref{cor1}}
\ead{robert\_kirby@baylor.edu}

\author[1]{John D. Stephens}
\ead{john\_stephens2@baylor.edu}

\cortext[cor1]{Corresponding author}

\affiliation[1]{organization={Department of Mathematics,
Baylor University},
addressline={1410 S.~4th St.},
city={Waco},
postcode={76706},
state={TX},
country={USA}}

\begin{abstract}
  \input{abstract.tex}
\end{abstract}

\begin{keyword}
Runge-Kutta method \sep finite element method \sep variational inequality \sep bounds constraints
\end{keyword}

\maketitle

\input{body_variational.tex}

\clearpage
\bibliographystyle{plain}
\bibliography{paper}

\end{document}

%% file: abstract.tex
Finite element methods provide accurate and efficient methods for the numerical solution of partial differential equations by means of restricting variational problems to finite-dimensional approximating spaces.
However, they do not guarantee enforcement of bounds constraints inherent in the original problem.
Previous work enforces these bounds constraints by replacing the variational equations with variational inequalities.
We extend this approach to collocation-type Runge--Kutta methods for time-dependent problems, obtaining (formally) high order methods in both space and time.
By using a novel reformulation of the collocation scheme, we can guarantee that the bounds constraints hold uniformly in time.
Numerical examples for a model of phytoplankton growth, the heat equation, and the Cahn--Hilliard system are given.

%% file: body_variational.tex
\section{Introduction}
Although finite elements provide an accurate, efficient, and theoretically robust toolkit for the numerical solution of partial differential equations,
they do not natively enforce bounds constraints associated with those problems.
Many factors contribute to these difficulties.
Representing range-restricted polynomials is known to be quite difficult~\cite{beatson1982restricted, lasserre2007sum, nesterov2000squared}.
Nochetto and Wahlbin~\cite{nochetto2002positivity} show that no linear operator can reproduce polynomials of degree greater than one while maintaining positivity, and maximum principles are quite delicate even for linear polynomials~\cite{draganescu2005failure}.

One may replace the variational problems defining finite element methods with variational inequalities that build bounds constraints into the feasible set.
Such approaches based on piecewise linear polynomials~\cite{chang2017variational, layton1996oscillation} give accurate solutions with uniform enforcement of the bounds constraints.
To obtain higher-order approximations, our work in~\cite{kirby2024high} uses the variational inequalities with constraints applied to the control net of Bernstein polynomials.  This guarantees uniform bounds constraints and gives good numerical results, but theoretical convergence rates are still an open question.
Alternatively, it is possible to apply bounds constraints only at the Lagrange nodes~\cite{barrenechea2024nodally}.  Here, optimal convergence rates are known, but the computed solution may violate constraints away from the Lagrange nodes.

We also note that proximal Galerkin methods~\cite{keith2023proximal} represent another, fundamentally different, approach to bounds enforcement.
Here, the original problem is replaced by solving a sequence of regularized and reformulated problems, each of which can be guaranteed to preserve bounds.
Although these methods are not formulated cleanly as a single variational inequality, the methods natively support high orders of approximation and provide scalable, mesh-independent convergence rates.  These are an important alternative approach.

Prior work on bounds constraints has primarily focused on stationary problems, and here we turn to evolutionary equations.
One may write single-stage methods such as backward Euler as variational problems for each time step and discretize in space by a bounds-constrained method.
However, single-stage methods are limited to low order.  Backward Euler gives only first order accuracy.  The implicit midpoint rule gives second order, but is not L-stable and so less suited for parabolic equations.
While multi-step methods such as BDF can achieve higher order by solving equations quite similar to backward Euler (and hence, posing a variational inequality to enforce bounds constraints), they are subject to the Dahlquist Barrier~\cite{dahlquist1963special}.

Even if spatial discretization led to bounds constrained ordinary differential equations, temporal discretization would still present a problem.
For example Bolley and Crouzeix~\cite{bolley1978conservation} show that a method that preserves positivity uniformly in the time step size cannot exceed first order.
However, techniques to circumvent this barrier (e.g.~\cite{horvath1998positivity, blanes2022positivity}) for ODE do not directly apply to the case of PDE discretized with variational inequalities.

In this work, we focus on implicit Runge--Kutta methods, with special emphasis on a subclass of fully implicit methods arising from continuous collocation.
Methods within this family can provide high orders of approximation without order barriers.
Although such methods require solving a rather complicated stage-coupled algebraic system at each time step, recent work on solvers~\cite{clines2022efficient,mmg,staff2006preconditioning,vanlent2005} can make them competitive with more commonly-practiced alternatives~\cite{abu2022monolithic,pazner2017stage,southworth2022fast1,southworth2022fast2}.
The Irksome project~\cite{farrell2021irksome,kirby2024extending} sits on top of Firedrake~\cite{Rathgeber:2016} to efficiently apply Runge--Kutta methods to finite element spatial discretizations, and it supports many of the cited approaches to fast solvers.

To obtain bounds-constrained discretization of time-dependent problems, we combine implicit Runge--Kutta methods with variational inequalities.
We replace the stage-coupled variational equation for computing the stage variables with variational inequalities.
For collocation-type RK methods, these stages are actually the coefficients of a collocating polynomial in the Lagrange basis.
Hence, the variational inequality guarantees that bounds constraints (at least, in some discrete sense) will hold at each collocation time.
Moreover, we show how to reformulate the RK method to compute the Bernstein form of the collocating polynomial.
This guarantees that bounds constraints will hold uniformly in time.

In Section~\ref{sec:setting}, we give an abstract presentation of a family of evolution equations  and give three particular model applications - an ODE model of phytoplankton growth, the heat equation, and the Cahn--Hilliard system.
Then, we formulate the Runge--Kutta and finite element discretization (with either variational equalities or inequalities) in Section~\ref{sec:method}.
Both Lagrange and Bernstein formulations of the collocation methods are given.
Numerical results evaluating these methods are presented in Section~\ref{sec:numres}, and conclusions and directions for future work are discussed in Section~\ref{sec:conc}.

\section{Problem setting}
  \label{sec:setting}
We let $\mathcal{V}$, $\mathcal{H}$ be Hilbert spaces with $\mathcal{V}$ compactly embedded in $\mathcal{H}$.
We let $(\cdot, \cdot)$ denote the $\mathcal{H}$ inner product
and $\langle \cdot , \cdot \rangle$ the duality pairing on $\mathcal{V}^\prime \times \mathcal{V}$.  
For $t \geq 0$, we let $F_t:\mathcal{V} \rightarrow \mathcal{V}^\prime$ be some bounded mapping (not necessarily linear), and we consider the abstract evolution equation of seeking some differentiable $y: \mathbb{R}^+ \rightarrow \mathcal{V}$ such that

\begin{equation}
  \label{eq:ev}
(y^\prime, v) = \langle F_t(y), v \rangle, \ \ \ v \in \mathcal{V},
\end{equation}
together with an appropriate initial condition
\begin{equation}
  y(0) = y_0 \in \mathcal{V}.
\end{equation}

In many applications, $y$ takes on values in some closed, convex subset $\mathcal{K} \subset \mathcal{V}$ at each time.
For example, physical or mathematical considerations may imply that $y$ takes on values in some restricted range, such as nonnegative functions.
Standard discretizations of~\eqref{eq:ev}, such as Galerkin finite elements in space combined with Runge-Kutta or some other scheme for time marching, need not produce approximations that respect these constraints, as seen in~\cite{bolley1978conservation, nochetto2002positivity}.

Instead, we proceed in the sense of Rothe~\cite{rothe1930zweidimensionale} by first discretizing~\eqref{eq:ev} in time to obtain a sequence of infinite-dimensional problems, then obtain spatial discretization via  variational inequalities~\cite{barrenechea2024nodally,chang2017variational,kirby2024high,layton1996oscillation}.
We consider a broader class of implicit Runge-Kutta methods.
The convergence of certain strongly stable methods to the solution of~\eqref{eq:ev} as the time step goes to zero is established for a very broad class of parabolic problems in~\cite{ostermann2002convergence}.
Discretizing the stage-coupled problem in space then leads to a fully discrete scheme.
A standard Galerkin scheme gives the same result as discretizing first in space and then applying the Runge-Kutta scheme to the discrete system.
Working with Runge-Kutta methods arising from collocation, spatial discretization with variational inequalities leads to bounds constraints at the discrete stages of the Runge-Kutta method, and we also derive a reformulation of the collocation scheme in terms of the Bernstein basis that leads to the collocating polynomial satisfying bounds constraints uniformly in time.

Before proceeding to formulate our temporal discretization along these lines,
we first present some model problems of interest.  Our first example is a nonlinear ODE model of phytoplankton growth on two nutrients.
This is a test problem presented in~\cite{bruggeman2007} for nonnegativity-preserving integrators for ODE. Seek $C$, $N$, $P$, $D:[0,T]\rightarrow \mathbb{R}$ such that

\begin{equation}\label{eq:PhytoSys}
\begin{split}
\frac{dC}{dt} &= -a r_{max} \frac{C}{K_C + C}\frac{N}{K_N + N} P,\\
\frac{dN}{dt} &= -b r_{max} \frac{C}{K_C + C}\frac{N}{K_N + N}P,\\
\frac{dP}{dt} &= r_{max}\frac{C}{K_C + C}\frac{N}{K_N + N}P - eP,\\
\frac{dD}{dt} &= eP,
\end{split}
\end{equation}
with initial conditions $C(0) = C_0$, $N(0) = N_0$, $P(0)=P_0$, and $D(0) = D_0$. In this model, $C$ is a carbon source, $N$ is a nitrogen source, $P$ is the population of phytoplankton, and $D$ is the accumulated detritus. Here, $a$, $b$, $K_C$, $K_N$, $r_{max}$, and $e$ are physical constants.

Then, we consider two PDE applications.  To fix notation,
we take $L^2(\Omega)$ to be the standard space of square-integrable functions over some domain $\Omega \in \mathbb{R}^d$ and $H^1(\Omega)$ its subspace of functions with square-integrable weak derivatives of order 1.  Additionally, $H^1_0(\Omega) \subset H^1(\Omega)$ contains only functions vanishing on $\partial \Omega$.

Our first PDE to consider is the heat equation, posed on $\Omega \subset \mathbb{R}^d$ for $1 \leq d \leq 3$.
Let  $f: [0, T] \rightarrow L^2(\Omega)$. 
We seek $u: [0, T] \rightarrow \mathcal{V}$ such that 
\begin{equation}\label{eq:heat_general}
    \left( u', v \right) + \left( \nabla u, \nabla v \right) = \left(f, v \right)
\end{equation}
for all $v \in \mathcal{V}$ and $0 \leq t \leq T$.  We close the problem with an initial condition
\begin{equation}
  u(0) = u_0 \in \mathcal{V},
\end{equation}
and the Dirichlet boundary condition
\begin{equation}\label{eq:heat_bc}
u(t)|_{\partial \Omega} = g \qquad \text{for } t\in [0, T].
\end{equation}

Next, the Cahn-Hilliard equation models phase separation in a binary fluid posed on a domain $\Omega \subset \mathbb{R}^d$ for $d=2,3$.
We take $\mathcal{V} = H^1(\Omega) \times H^1(\Omega)$,
and we write the system for the order parameter $c:\Omega \rightarrow [-1,1]$ and chemical potential $\mu : \Omega \rightarrow \mathbb{R}$ as
\begin{equation}\label{eq:cahn_hilliard_general}
  \begin{split}
    \left(\frac{\partial c}{\partial t}, v \right) &= - M\left(\nabla \mu, \nabla v\right), \\
    \left(\mu, w \right) &= \left( F'(c), w\right) + \epsilon^2\left( \nabla c, \nabla w \right) \\
  \end{split}
\end{equation}
for all $v, w \in H^1(\Omega)$ and all $0 < t \leq T$.
We utilize the Flory-Huggins type logarithmic potential 
\begin{equation}\label{eq:CH_logarithmic_potential}
F(s) = \frac{\theta_0}{2}[(1 + s)\ln(1 + s) + (1 - s)\ln(1 - s)] - \frac{\theta_c}{2}s^2
\end{equation}
\cite{flory1942, huggins1942}, where $\theta_0$ and $\theta_c$ are positive constants, $\theta_0 < \theta_c$.

Here, the order parameter $c$ indicates the relative amount of fluid in each state, with $\pm 1$ indicating one state or the other and values in between indicating a mixture.  The parameter $M$ indicates the mobility and can, in more general settings, depend on $c$.  The parameter $\epsilon$ is a constant indicating the size of free energy for a given concentration gradient.
The system is closed with boundary conditions, which we take as homogeneous Neumann or periodic conditions on both $c$ and $\mu$.

Physically, only values of $c$ between -1 and 1 make sense, and under mild assumptions it has been shown (e.g.~\cite{debussche1995cahn}) that the analytical solution to~\eqref{eq:cahn_hilliard_general} exists and is contained in the interval $(-1, 1)$ a.e. for all times $t > 0$. This is an interesting and nontrivial target application for our bounds-constrained Runge-Kutta methods.

\section{Method development}
\label{sec:method}
In this section, we more concretely present our approach to obtaining fully-discrete and bounds-preserving PDE approximations.
First, we describe collocation-type Runge--Kutta methods for the evolution equation.
These give stage-coupled variational problems for the stages at each time level, which can then be discretized by a Galerkin method.
Bounds constraints can be enforced on the stages by means of replacing the variational problems with variational inequalities.
While this method ensures bounds constraints hold at each stage, we can also reformulate the scheme in terms of the Bernstein coefficients of the collocating polynomial rather than its values at the time levels to obtain a uniform-in-time bounds constraint.

\subsection{Runge--Kutta schemes}
We partition the time interval $[0, T]$ into intervals
$(t^n, t^{n+1})$ with $t^0 = 0$ and $t^N = T$.  As a notational convenience, we will assume a uniform time step size of $t^{n+1} - t^n = k$, although there is no practical limitation of Runge-Kutta methods to uniform step sizes.

Typically, an $s$-stage Runge-Kutta method is encoded by a Butcher tableau
\begin{equation}
  \begin{array}{c|c}
    \mathbf{c} & A \\ \hline
    & \mathbf{b}
  \end{array},
\end{equation}
where $\mathbf{b}, \mathbf{c} \in \mathbb{R}^s$ and $A \in \mathbb{R}^{s \times s}$.

To define the variational problem associated with a Runge-Kutta scheme for~\eqref{eq:ev}, we define by $\mathcal{V}^s$ the $s$-way product $\mathbf{\mathcal{V}}^s = \prod_{i=1}^s \mathcal{V} = \mathcal{V} \times \mathcal{V} \times \dots \times \mathcal{V}$.

Given an approximation to $y$ at time $t^n$, a Runge-Kutta scheme for~\eqref{eq:ev} seeks $Y \in \mathbf{\mathcal{V}}^s$ such that
\begin{equation}
  \label{eq:stagedef}
  \left(Y_i, v_i\right)
  = \left( y^n, v_i \right) + k \sum_{j=1}^s A_{ij}
  \langle F_{t^n + \mathbf{c}_j k}(Y_j), v_i \rangle
\end{equation}
for all $v_i \in V, 1 \leq i \leq s$.
The stage values $Y_j$ approximate the solution $y$ at time $t^n +\mathbf{c}_j h$.
We write this abstractly as a variational problem seeking $Y \in \mathcal{V}^s$ such that
\begin{equation}
  \label{eq:abstractstage}
  \langle \mathcal{F}(Y), V\rangle = 0, \ \ \ V \in \mathcal{V}
\end{equation}

Then, $y^{n+1} \in V$ arises from the much simpler variational problem
\begin{equation}
  \label{eq:updatestage}
  \left( y^{n+1}, v \right)  = \left( y^n, v \right)
  + k \sum_{j=1}^s b_j \langle F_{t_n + \mathbf{c}_j k}(Y_j), v \rangle.
\end{equation}
This is be posed variationally since $F_t$ maps $V$ into its dual.
When the Runge-Kutta method is \emph{stiffly accurate}, it holds that $y^{n+1}$ coincides with the final stage value $Y_s$ so that~\eqref{eq:updatestage} is not needed in practice.

Of particular interest here are the subset of Runge-Kutta methods that arise from continuous collocation, such as the classical Gauss-Legendre and RadauIIA families.
A continuous collocation method seeks a polynomial $u: [t^n, t^{n+1}] \rightarrow \mathcal{V}$ such that
\begin{equation}
  \label{eq:contcoll}
  \begin{split}
    u(t^n) & = y^n \\
    \left( u^\prime(t^n + \mathbf{c}_i k), v \right) & = \left\langle F_{t^n+\mathbf{c}_ik}(u(t^n+\mathbf{c}_i k)), v \right\rangle, \ \ \ v \in \mathcal{V}
  \end{split}
\end{equation}
Any continuous collocation method can be written as a Runge-Kutta method~\cite{wanner1996solving}.
Using the stage value formulation of Runge-Kutta we have described, it is possible to construct the collocating polynomial $u$ through interpolation.
We first consider the case where in the Butcher tableau we have $\mathbf{c}_1 \neq 0$.
Define $\mathbf{c}_0 = 0$, and $\overline{Y} \in \mathbf{V}^{s+1}$ by
\begin{equation}
  \label{eq:ybar}
\overline{Y} = \begin{bmatrix} y^n \\ Y \end{bmatrix},
\end{equation}
and we index this vector from 0 rather than 1.
Then, we  let $\{ \ell_i(\tau) \}_{i=0}^s$ be the Lagrange polynomials associated with points $\mathbf{c}_i$ on $[0, 1]$, and the collocating polynomial is just

\begin{equation}
  \label{eq:interp}
  u(t^n + \tau k) = \sum_{i=0}^s \overline{Y}_i \ell_i(\tau).
\end{equation}

As an alternative to~\eqref{eq:updatestage}, we note that $y^{n+1}$ can be computed by evaluating the collocating polynomial
\begin{equation}
  y^{n+1} = u(t^n + k),
\end{equation}
and when $\mathbf{c}_s = 1$, we recover the stiffly accurate case of $y^{n+1} = Y_s$.

When $\mathbf{c}_1 = 0$, such as occurs in the LobattoIIIA family, the collocating polynomial must satisfy a more general confluent interpolation problem:
\begin{equation}
  \label{eq:confinterp}
  \begin{split}
    u(t^n + \mathbf{c}_i k) & = Y_i, \ \ \ 1 \leq i \leq s \\
    ( u^\prime(t^n), v ) & = \left\langle F_{t^n}(u(t^n)), v \right\rangle, \ \ \ v \in \mathcal{V}
  \end{split}
\end{equation}
One can build the collocating polynomial of degree $s$  by
\begin{equation}
  \label{eq:confinterp_stages}
  u(t^n + \tau k) = u^\prime(t^n) \tilde{\ell}(\tau) + \sum_{i=1}^s Y_i \widehat{\ell}_i(\tau),
\end{equation}
where $\tilde{\ell}$ vanishes at all $c_i$ and has $\tilde{\ell}^\prime(0) = 1$ and
the $\{ \widehat{\ell}_i \}_{i=1}^s$ have vanishing derivative at 0 and satisfy $\widehat{\ell}_i(\mathbf{c}_j) = \delta_{ij}$. 
The primary use case here would be the LobattoIIIA family, which includes Crank-Nicolson and generalizes it to high orders.
However, this family lacks the L-stability of the RadauIIA family for dissipitave problems or the symplecticity of the Gauss-Legendre family.

\subsection{Spatial discretization}
A Galerkin method for ~\eqref{eq:abstractstage} may be posed by taking some finite-dimensional subspace $\mathcal{V}_h \subset \mathcal{V}$
and $\mathcal{V}_h^s$ the $s$-way product of $\mathcal{V}_h$ with itself.
Then, we seek $Y_h \in \mathcal{V}_h^s$ such that
\begin{equation}
  \label{eq:abstractgalerkin}
  \langle \mathcal{F}(Y_h), V_h \rangle = 0, \ \ \ V_h \in \mathcal{V}_h^s.
\end{equation}

A fully discrete method is obtained by posing some initial state $y^0_h \in \mathcal{V}_h$.  For each $n \geq 0$, having obtained $y_h^n \in \mathcal{V}_h$, we solve~\eqref{eq:abstractgalerkin} to obtain stage values, and then compute
$y^{n+1}_h$ by
\begin{equation}
  \label{eq:discreteupdatestage}
  \left( y_h^{n+1}, v_h \right)  = \left( y^n, v_h \right)
  + k \sum_{j=1}^s b_j \langle F_{t_n + \mathbf{c}_j k}(Y_{h,j}), v_h \rangle, \ \ \ v_h \in \mathcal{V}_h,
\end{equation}
or, if appropriate, just set $y_h^{n+1} = Y_{h,s}$.

Our target PDE applications are posed on (subspaces of) $H^1$ or products thereof, so the standard spaces of continuous piecewise polynomials are appropriate here.
We let $\{\mathcal{T}_h\}_h$ be a quasiuniform family of triangulations of the problem domain $\Omega$.
Let $\mathcal{W}_h \subset H^1(\Omega)$ consist of continuous piecewise polynomials of degree $r$ with $r \geq 1$ and $\mathcal{W}_{h, 0} \subset H^1_0(\Omega) \cap \mathcal{W}_h$  its subspace with vanishing trace.
For the heat equation, we take 
$\mathcal{V}_h = \mathcal{W}_{h, 0}$, and
for the Cahn-Hilliard equation, we take $\mathcal{V}_h = \mathcal{W}_h \times \mathcal{W}_h$ since Neumann or periodic conditions are imposed on each variable.

\subsection{Enforcing bounds constraints}
Although a standard Galerkin method delivers an approximation that is well-defined, independent of the choice of basis, techniques for enforcing bounds constraints depend rather intimately on this choice.
Determining the positivity of a generic multivariate polynomial is in fact NP-hard~\cite{lasserre2007sum}, so here we restrict ourselves to some workable, if imperfect approaches.

Let $\mathcal{B} = \{ \psi_i \}_{i=1}^{\dim\mathcal{W}_h}$ be a basis for $\mathcal{W}_h$, which we will typically take either to be the piecewise Lagrange or Bernstein polynomials.
Then, we define the set
\begin{equation}
  \mathcal{J}^{\mathcal{B}, I}_h
  = \left\{ \sum_{i=1}^{\dim \mathcal{W}_h} c_i \psi_i : c_i \in I \right\},
\end{equation}
to comprise only members of $\mathcal{W}_h$ whose basis coefficients relative to $\mathcal{B}$ satisfy the bounds constraints.

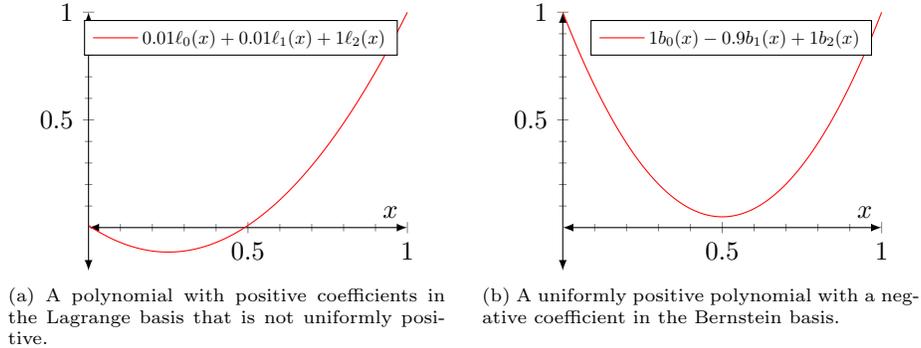
\begin{figure}
  \centering
  \begin{subfigure}[t]{.48\textwidth}
    \centering
    \begin{tikzpicture}
      \centering
        \begin{axis}[width=\textwidth, grid style={line width=.1pt, draw=gray!10}, legend style={nodes={scale=0.7, transform shape}}, minor tick num=4, major grid style={line width=.2pt,draw=gray!50}, axis lines=middle,
            axis line style={latex-latex},
            legend pos=north east,
            ymin=-0.2,
            ymax=1,
            ytick={0.0, 0.5, 1.0},
            yticklabels={$0$, $0.5$, $1$},
            xtick={0,0.5,1}, 
            xticklabels={$0$, $0.5$, $1$},
            xlabel=$x$]
          \addplot[domain=0:1, samples=100, color=red] {0.01*2*(x-0.5)*(x-1) - 0.01*4*x*(x-1) + 2*x*(x-0.5)};
          \addlegendentry{$0.01\ell_0(x) + 0.01\ell_1(x) + 1\ell_2(x)$} ;
        \end{axis}
    \end{tikzpicture}
    \subcaption{A polynomial with positive coefficients in the Lagrange basis that is not uniformly positive.}
    \label{fig:edgecaselagrange}
  \end{subfigure}\hspace{0.04\textwidth}%
  \begin{subfigure}[t]{.48\textwidth}
    \centering
    \begin{tikzpicture}
      \centering
        \begin{axis}[width=\textwidth, grid style={line width=.1pt, draw=gray!10}, legend style={nodes={scale=0.7, transform shape}}, minor tick num=4, major grid style={line width=.2pt,draw=gray!50}, axis lines=middle,
            axis line style={latex-latex},
            legend pos=north east,
            ymin=-0.2,
            ymax=1,
            ytick={0.0, 0.5, 1.0},
            yticklabels={$0$, $0.5$, $1$},
            xtick={0,0.5,1}, 
            xticklabels={$0$, $0.5$, $1$},
            xlabel=$x$]
          \addplot[domain=0:1, samples=100, color=red] {(1-x)^2 -1.8*x*(1-x)+x^2};
          \addlegendentry{$1b_0(x) -0.9b_1(x) + 1b_2(x)$} ;
        \end{axis}
    \end{tikzpicture}
    \subcaption{A uniformly positive polynomial with a negative coefficient in the Bernstein basis.}
    \label{fig:edgecasebernstein}
  \end{subfigure}
  \caption{Edge case polynomials in the Lagrange and Bernstein bases.}
  \label{fig:edgecasepolys}
\end{figure}

When $\mathcal{B}$ consists of piecewise Lagrange polynomials, $\mathcal{J}_h^{\mathcal{B}, I}$ contains polynomials that satisfies the bounds constraints at the interpolating nodes.
In this case, every member of $\mathcal{W}_h$ with range in $I$ lies in $\mathcal{J}_h^{\mathcal{B}, I}$, but so do many that do not -- the bounds are readily violated between nodes.  Figure~\ref{fig:edgecaselagrange} shows a simple example of a univariate polynomial on $[0, 1]$ for which this is the case.

\begin{figure}
  \begin{center}
    \begin{tikzpicture}[scale=3.0]
      \draw[->] (0, 0) -- (1.1, 0) node[right] {$x$};
      \draw[->] (0, 0) -- (0, 1.1) node[above] {$y$};
      \draw [smooth, samples=100, domain=0:1] plot({\x}, {(1-\x)^3});
      \draw [smooth, samples=100, domain=0:1] plot({\x}, {3*\x * (1-\x)^2});
      \draw [smooth, samples=100, domain=0:1] plot({\x}, {3*(\x)^2 * (1-\x)});
      \draw [smooth, samples=100, domain=0:1] plot({\x}, {\x^3});
  \end{tikzpicture}
  \end{center}
  \caption{Cubic Bernstein polynomials.}
  \label{fig:bern3}
\end{figure}
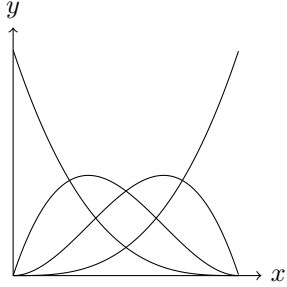

In contrast, suppose that $\mathcal{B}$ is be the basis of piecewise Bernstein polynomials.  On the unit interval, the Bernstein basis of degree $n$
is given by
\begin{equation}
  b^n_i(x) = {n \choose i} x^i (1-x)^{n-i}, \ \ \ 0 \leq i \leq n.
\end{equation}
These polynomials give a nonnegative partition of unity forming a basis for polynomials of degree $n$.  They are readily mapped to any compact interval $[a, b]$, and given their geometric decomposition~\cite{arnold2009geometric}, they can be easily assembled across cells to form $C^0$ piecewise polynomials.  The cubic basis is given in Figure~\ref{fig:bern3}.

The Bernstein polynomials possess a convex hull property -- the graph of
\(
p(x) = \sum_{i=0}^n c_i b^n_i(x)
\)
for $x \in [0, 1]$ 
lies in the convex hull of its control net $\left\{\left( \tfrac{i}{n}, c_i \right)\right\}_{i=0}^n$~\cite{LaiSch07}.  So, if the coefficients $c_i \in I$, then the polynomial can only take values in $I$.
However, there exist polynomials that are uniformly positive on $[0, 1]$ but have negative coefficients in the Bernstein basis, as shown in
Figure~\ref{fig:edgecasebernstein}.
This situation is the opposite of the Lagrange basis -- every member of $\mathcal{W}_h^{\mathcal{B}, I}$ is bounds-constrained, but we only have a proper subset of constrained polynomials.

We define $\mathcal{J}^{[m,M]}_h$ to be those members of $\mathcal{W}_h$ whose degrees of freedom in the chosen basis (either Lagrange or Bernstein) satisfy the bounds constraints.  With the Lagrange basis, we obtain all members of $\mathcal{W}_h$ that satisfy the bounds constraints uniformly, but have an outside approximation since members may violate those bounds between the Lagrange nodes.
With the Bernstein basis, $\mathcal{J}_h^{[m,M]}$ contains only functions uniformly satisfying the bounds constraints, but not all of them.
Constraining the individual $H^1$ approximating space $\mathcal{W}_h$ in this way, we can define $\mathcal{J}_h \subset \mathcal{V}_h$ by constraining the desired fields.
  
This construction, in either the Lagrange or Bernstein basis, allows us to define (discretely) feasible sets.   For example, discretizing the heat equation with nonnegative values gives $\mathcal{K}_h = \mathcal{J}_h^{\mathcal{B}, [0, \infty)}$, while constraining the order parameter $c$ in Cahn-Hilliard gives $\mathcal{K}_h = \mathcal{J}_h^{\mathcal{B}, [-1+\delta,1-\delta]} \times \mathcal{W}_h$, where the choice of $\delta$ is described in Section~\ref{sec:cahn_hilliard}.

After defining the discrete feasible set, $\mathcal{K}_h$, we define
$\mathcal{K}_h^s$ to be its $s$-way product,
Then, we replace the variational equation~\eqref{eq:abstractgalerkin} with a variational inequality seeking $Y_h \in \mathcal{K}_h^s$ such that
\begin{equation}
  \label{eq:abstdiscvi}
  \langle \mathcal{F}(Y_h), V_h - Y_h \rangle \geq 0, \ \ \ V_h \in \mathcal{K}_h^s.
\end{equation}
Then, we use these stage values to form the constrained version of the collocation polynomial, $u(t^n + \tau k)$, and take $y^{n+1} = u(t^{n+1})$ as the update, or just $y^{n+1} = Y_{h,s}$ in the stiffly accurate case.

By construction,~\eqref{eq:abstdiscvi} can only produce solutions for which the discretized Runge-Kutta stages satisfy the bounds constraint (in some sense).
Moreover, for stiffly accurate collocation schemes (in particular, RadauIIA methods) this also guarantees that $y^{n+1}_h$ will satisfy the bounds constraints.
For non-stiffly accurate collocation schemes such as Gauss-Legendre, however, this is not actually the case.  
With $\{ \ell_i \}_{i=0}^s$ the Lagrange polynomials associated with 0 and the collocation nodes, a simple argument shows that every other $\ell_i$ will take negative values at 1.

In this case, it would be possible to fuse together the stage problem~\eqref{eq:abstdiscvi} with the update formula~\eqref{eq:discreteupdatestage} into an $(s+1) \times (s+1)$ variational inequality, including the bounds constraints on the stages and new values $y^{n+1}$.
However, we can obtain an even stronger method that imposes bounds constraints uniformly in time by reformulating the collocation polynomial in terms of its Bernstein coefficients, and we turn now to this case.

\subsection{Bernstein form of the collocating polynomial}
To this point, we have worked with the the collocation polynomial through its value at time step $t^n$ and at the discrete time levels $t^n + \mathbf{c}_i k$.
Now, we show how to reformulate the polynomial in terms of its Bernstein representation, whence the variational inequality analogous to~\eqref{eq:abstdiscvi} will imply uniform-in-time bounds constraints.

In continuous collocation with $\mathbf{c}_1\neq 0$, the previous value $y^n$ with the stage values $Y_j$ define the Lagrange form of the collocation polynomial.
As before, we put $\mathbf{c}_0 = 0$ and let $\{\ell_i\}_{i=0}^s$ be the Lagrange polynomials associated with these nodes and let $\overline{Y}$ be defined as in~\eqref{eq:ybar} and the collocating polynomial by~\eqref{eq:interp}.
Equivalently, we may write this polynomial in the Bernstein basis by
\begin{equation}
  u(t^n + k \tau) = \sum_{j=0}^s \overline{Z}_j b_j(\tau).
\end{equation}

We note that $\overline{Z}_0 = \overline{Y}_0 = y^n$ since $b_0(0) = 1$ and $b_i(0) = 0$ for $1 \leq i \leq s$.  We write the vector of
Bernstein coefficients as
\begin{equation}
  \overline{Z} = \begin{bmatrix} y^n \\ Z
  \end{bmatrix}.
\end{equation}

Evaluating this Bernstein form at each Lagrange node, 
\begin{equation}
  \overline{Y}_i = u(t^n + k \mathbf{c}_i) =
  \sum_{j=0}^s \overline{Y}_j b_j (\mathbf{c}_i),
\end{equation}
so that
\begin{equation}
  \label{eq:bernconvertcont}
  \overline{Y} = \overline{V} \overline{Z},
\end{equation}
where $\overline{V}$ is the Bernstein-Vandermonde matrix
\begin{equation}
  \overline{V}_{ij} = b_j(\mathbf{c}_i), \ \ \ 0 \leq i, j \leq s.
\end{equation}

We can write~\eqref{eq:bernconvertcont} as
\begin{equation}
\begin{bmatrix}
  y^n \\ Y
\end{bmatrix}
= 
\begin{bmatrix} 1 & \mathbf{0}^T \\ \mathbf{v} & V \end{bmatrix}
\begin{bmatrix} y^n \\ Z \end{bmatrix}
= \begin{bmatrix}
  y^n \\
  V Z + y^n \mathbf{v}
  \end{bmatrix}.
\end{equation}

Now, since $Y = VZ + y^n \mathbf{v}$, we can rewrite~\eqref{eq:stagedef} as a modified variational
problem seeking the Bernstein coefficients $Z \in \mathbf{V}^s$  by
\begin{equation}
  \label{eq:stagedefbern}
  \left( \left(V Z\right)_i, v_i \right)  = \left( y^n( 1 - \mathbf{v}_i ), v_i \right)
  + k \sum_{j=1}^s A_{ij} \left\langle
  F_{t^n+\mathbf{c}_j k}\left(\left(V Z\right)_j + y^n \mathbf{v}_j \right), v_i \right\rangle.
\end{equation}

Since only the last Bernstein polynomial is nonzero at $\tau=1$, we can update the solution by evaluating the collocation polynomial with
\begin{equation}
  y^{n+1} = Z_s,
\end{equation}
and this holds regardless of whether the original collocation schemes was stiffly accurate with $\mathbf{c}_s = 1$.

When the Butcher tableau has $\mathbf{c}_1 = 0$, we saw that a more general confluent interpolation problem was required to define the collocating polynomial. In this case, the Bernstein-Vandermonde matrix must be replaced with a confluent Bernstein-Vandermonde matrix. Evaluating~\eqref{eq:confinterp_stages} at $\tau = 0$,
\begin{equation}
u(t^n) = Y_1 = y^n.
\end{equation}
Writing the collocation polynomial in the Bernstein basis,
\begin{equation}
u(t^n + \tau k) = \sum_{j = 0}^s \overline{Z}_j b_j(\tau),
\end{equation}
and evaluating $u$ at the collocation nodes and $u'$ at $0$, we arrive at the following system of equation:
\begin{align}
Y_i &= \sum_{j = 0}^s \overline{Z}_j b_j(\mathbf{c}_i), \qquad i = 1, \dots, s,\label{eq:confluent_stages}\\
\left\langle F_{t^n}(y^n), v \right\rangle &= \left(\sum_{j = 0}^s \overline{Z}_j b_j'(0), v \right)\label{eq:confluent_deriv}, \qquad v\in V,
\end{align}
Note,~\eqref{eq:confluent_deriv} must be cast in a variational sense because $F_{t^n}$ maps $V$ into its dual. By the Riesz Representation theoreom, there exists some element $f_{t^n}(y^n) \in V$ such that $\langle F_{t^n}(y^n), v\rangle = (f_{t^n}(y^n), v)$ for all $v \in V$. So, we may replace~\eqref{eq:confluent_deriv} with 
\begin{equation}\label{eq:riesz_map_deriv}
f_{t^n}(y^n) = \sum_{j = 0}^s \overline{Z}_j b_j'(0).
\end{equation} 
Thus, we must solve the system
\begin{equation}
\begin{bmatrix}  b_0(0) & b_1(0) & \dots &b_s(0)\\ b_0'(0) & b_1'(0) & \dots &b_s'(0)\\  b_0(\mathbf{c}_2) & b_1(\mathbf{c}_2) & \dots &b_s(\mathbf{c}_2)\\
\vdots & \vdots & & \vdots \\
b_0(\mathbf{c}_s) & b_1(\mathbf{c}_s) & \dots & b_s(\mathbf{c}_s)
\end{bmatrix}\overline{Z} = \begin{bmatrix} y^n \\ f_{t^n}(y^n) \\ Y_2 \\ \vdots \\ Y_s\end{bmatrix}.
\end{equation}
Recalling that $b_0(0) = 1, b_i(0) = 0$ for $i >1$, and 
\begin{equation}
b_i'(0) = \begin{cases} -s & \text{ if } i = 0,\\ s &\text{ if } i =1,\\ 0 &\text{ otherwise}\end{cases}, 
\end{equation}
we may write this system as 
\begin{equation}
\begin{bmatrix} 1 & 0 & \dots & 0 \\ 0 & 1 &\dots & 0 \\ b_0(\mathbf{c}_2) & b_1(\mathbf{c}_2) & \dots & b_s(\mathbf{c}_2) \\ \vdots & \vdots & & \vdots \\ b_0(\mathbf{c}_s) & b_1(\mathbf{c}_s) & \dots  &b_s(\mathbf{c}_s)\end{bmatrix} \overline{Z} - \begin{bmatrix} 0 \\ \frac{1}{s}f_{t^n}(y^n) \\ 0 \\ \vdots \\ 0\end{bmatrix} = \begin{bmatrix} y^n\\ y^n\\ Y_2 \\ \vdots \\ Y_s\end{bmatrix}
\end{equation}
We partition this matrix as before to obtain
\begin{equation}\label{eq:confluent_system}
\begin{bmatrix} 1 & \mathbf{0}^T \\ \mathbf{v} & V\end{bmatrix} \begin{bmatrix} y^n \\ Z \end{bmatrix} - \begin{bmatrix} 0\\  w \end{bmatrix} = \begin{bmatrix} y^n \\ Y \end{bmatrix},
\end{equation}
where $\mathbf{v} = [0, b_0(\mathbf{c}_2), \dots, b_0(\mathbf{c}_s)]^T$, $w = \left[\frac{1}{s}f_{t^n}(y^n), 0, \dots, 0\right]^T$ and
\begin{equation}
V = \begin{bmatrix} 1 & 0 & \dots & 0\\ b_1(\mathbf{c}_2) & b_2(\mathbf{c}_2) &  \dots & b_s(\mathbf{c}_2) \\ \vdots  & \vdots & & \vdots \\ b_1(\mathbf{c}_s) & b_2(\mathbf{c}_s) & \dots & b_s(\mathbf{c}_s)\end{bmatrix}.
\end{equation}
Now, substituting~\eqref{eq:confluent_system} into~\eqref{eq:stagedef}, we may solve for the coefficients of the confluent interpolating polynomial written in terms of the Bernstein basis.

In either case, we may write the reformulated system abstractly in the same form as~\eqref{eq:abstractstage} as
\begin{equation}
  \label{eq:abstractbernstage}
  \langle \mathcal{F}^b(Z), V\rangle = 0, \ \ \ V \in \mathcal{V},
\end{equation}
and, like~\eqref{eq:abstdiscvi} we may also force the Bernstein stages to satisfy the bounds constraints by posing the variational equality of
seeking $Z \in \mathcal{K}_h^s$ such that 
\begin{equation}
  \label{eq:abstdiscvibern}
  \langle \mathcal{F}^b(Z), V - Z\rangle \geq 0, \ \ \ V \in \mathcal{K}_h^s,
\end{equation}

If $\mathbf{c}_1 \neq 0$, on $[t^n, t^{n+1}]$, we may represent the solution in either Lagrange or Bernstein form by
\begin{equation}
  y(t^n + \tau k) = \sum_{i=0}^s \overline{Y}_i \ell_i(\tau)
  = \sum_{i=0}^s \overline{Z}_i b^{s+1}_i(\tau),
\end{equation}
while if $c_1 = 0$, we have the Bernstein representation
\begin{equation}
y(t^n + \tau k) = \sum_{i=0}^s \overline{Z}_i b^{s+1}_i(\tau).
\end{equation}
We have actually proposed four different approaches to enforcing bounds constraints.
First, the finite element space may be parameterized with either Lagrange or Bernstein polynomials.  Athough these choices in~\eqref{eq:abstractgalerkin} give the same solution to discrete variational problems, applying bounds constraints to the coefficients in the variational inequality~\eqref{eq:abstdiscvi} give genuinely different results.
Second, writing the collocation problem in either the Lagrange or Bernstein basis again leads to different methods with variational inequalities.

Using the Lagrange basis in both space and time leads to a method that satisfies the bounds constraints at the Lagrange nodes at each discrete collocation stage, although no guarantees hold at other points in space and time.
On the other hand, using the Bernstein basis in both space and time leads to a method that satisfies the bounds constraints uniformly in both variables.
Using the Bernstein polynomials in space with the standard Lagrange form of the collocating polynomial gives a solution that, at each discrete Runge-Kutta stage, satisfies the bounds constraints uniformly in space.
By constrast, using the Lagrange polynomials in space but the Bernstein form of the collocating polynomial leads to solutions that satisfy the bounds constraints at each spatial node uniformly in time.

\section{Numerical results}
\label{sec:numres}

All numerical experiments were completed using the Firedrake Project~\cite{Rathgeber:2016}, a Python package designed for the efficient implementation of finite element methods. The Irksome project~\cite{farrell2021irksome,kirby2024extending} sits on top of Firedrake to efficiently apply Runge--Kutta methods to finite element discretizations, and was used extensively throughout. Firedrake provides high-level access to PETSc~\cite{petsc-user-ref,petsc-efficient} through petsc4py~\cite{dalcin2011}. All nonlinear probems were solved using the reduced space Newton-type VI solver {\ttfamily vinewtonrsls}, as implemented in PETSc.  We set the absolute tolerance to $10^{-8}$, and solved the resulting system with the LU-factorization. The snapshots in section~\ref{sec:cahn_hilliard} were generated using ParaView~\cite{ParaView}.\\

\noindent \textbf{Notation}: We have described four different methods determined by the choice of spatial finite element parameterization (between Lagrange and Bernstein), and the choice of representation of  the collocation polynomial (between Lagrange and Bernstein). 

The following numerical examples will showcase each of these different methods, often with varying spatial polynomial degrees, and different collocation polynomial degrees. The spatial finite element and piecewise polynomial order will be denoted by $\mathcal{P}_s$, where $\mathcal{P}$ will take on the symbol $\mathcal{L}$ or $\mathcal{B}$ to denote the Lagrange and Bernstein finite elements, respectively, and $s$ indicates the polynomial degree. The RadauIIA collocation scheme will be denoted by RIIA, taking an argument of the form $\mathcal{P}_s$ to denote the basis chosen to represent the collocation polynomial and the number of stages. Finally, VP will be appended in the event a variational problem is being solved, and VI will be appended when the variational problem is replaced with a variational inequality. For example, if piecewise-cubic Lagrange finite elements are used in space, timestepping is done with the 2-stage RadauIIA method, the collocation polynomial is represented in the Bernstein basis, and bounds constraints are enforced via a variational inequality, we will denote this method by $\mathcal{L}_3$-RIIA($\mathcal{B}_2$)-VI.

\subsection{Phytoplankton}

As a first numerical example, we consider~\eqref{eq:PhytoSys}, a model of phytoplankton growth on two nutrients.  We choose the same physical parameters and initial conditions as in~\cite{bruggeman2007}. Specifically, let $a = b = K_C = K_N = r_{max} = 1$, 
$e = 0.3$, $C_0 = 29.98$, $N_0 = 9.98$, and $P_0 = D_0 = 0.01$. The exact solution is nonnegative uniformly in time, and satisfies the linear invariants 
$C+P+D = 30$ and $N + P + D = 10$, representing the conservation of total carbon and total nitrogen in the system, respectively. In place of an analytical solution, we compare our approximations against a high-resolution numerical solution computed using RIIA($\mathcal{L}_2$)-VP with timestep $k = 10^{-5}$.

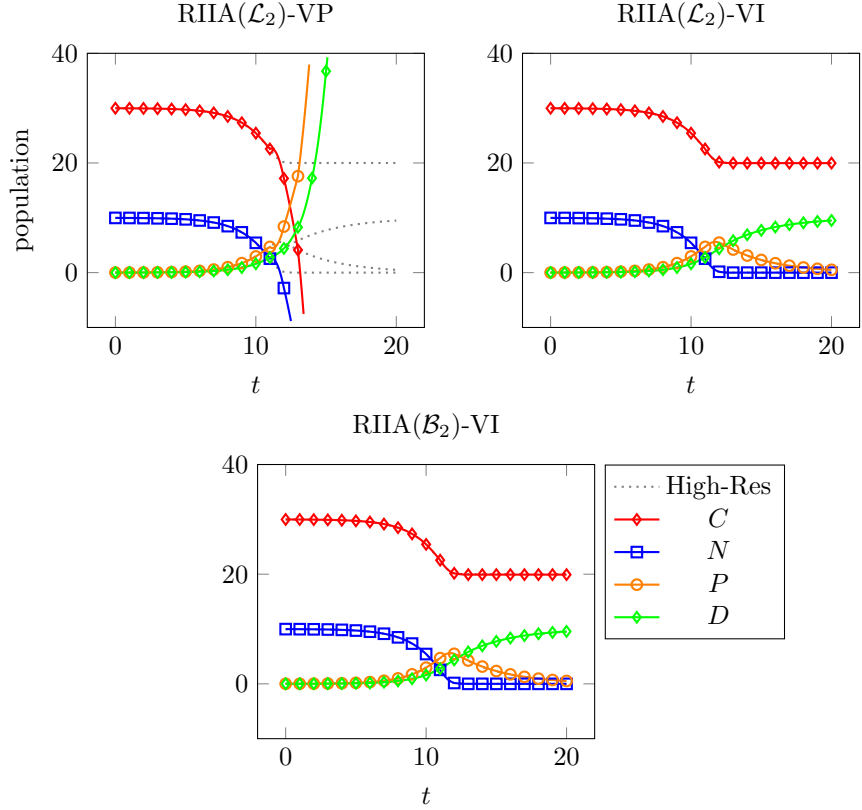
\begin{figure}[ht]
\centering
\begin{subfigure}{0.5\textwidth}
\centering
\begin{tikzpicture}
\centering
\begin{axis}[width=\textwidth, legend style={nodes={scale=1.0, transform shape}}, 
        title={RIIA($\mathcal{L}_2$)-VP},
	ylabel near ticks,
	legend pos=outer north east,
	xlabel=$t$,
	ymin=-10,
	ymax=40,
    	ylabel={population}]
    
    \addplot[gray, dotted, thick] table[y=C, x=time, col sep=comma] {phyto_high_res_filtered.csv};
    \addlegendentry {High-Res}
    \addplot[gray, dotted, thick] table[y=N, x=time, col sep=comma, forget plot] {phyto_high_res_filtered.csv};
    \addplot[gray, dotted, thick] table[y=P, x=time, col sep=comma, forget plot] {phyto_high_res_filtered.csv};
    \addplot[gray, dotted, thick] table[y=D, x=time, col sep=comma, forget plot] {phyto_high_res_filtered.csv};
    
        \addplot[red, mark=diamond, mark repeat=10, restrict y to domain=-10:40, unbounded coords=discard, thick] table[y=C, x=time, col sep=comma] {phyto_L2F_1.0.csv};
    \addlegendentry {$C$}
        \addplot[blue, mark=square, mark repeat=10, restrict y to domain=-10:40, unbounded coords=discard, thick] table[y=N, x=time, col sep=comma] {phyto_L2F_1.0.csv};
    \addlegendentry {$N$}
        \addplot[orange, mark=o, mark repeat=10, restrict y to domain=-10:40, unbounded coords=discard, thick] table[y=P, x=time, col sep=comma] {phyto_L2F_1.0.csv};
    \addlegendentry {$P$}
        \addplot[green, mark=diamond, mark repeat=10, restrict y to domain=-10:40, unbounded coords=discard, thick] table[y=D, x=time, col sep=comma] {phyto_L2F_1.0.csv};
    \addlegendentry {$D$}
\legend{}
\end{axis}
\end{tikzpicture}
\end{subfigure}%
\begin{subfigure}{0.5\textwidth}
\centering
\begin{tikzpicture}
\centering
\begin{axis}[width=\textwidth, legend style={nodes={scale=1.0, transform shape}}, 
        title={RIIA($\mathcal{L}_2$)-VI},
	ylabel near ticks,
	legend pos=outer north east,
	xlabel=$t$,
	ymin=-10,
	ymax=40]
    
    \addplot[gray, dotted, thick] table[y=C, x=time, col sep=comma] {phyto_high_res_filtered.csv};
    \addlegendentry {High-Res}
    \addplot[gray, dotted, thick] table[y=N, x=time, col sep=comma, forget plot] {phyto_high_res_filtered.csv};
    \addplot[gray, dotted, thick] table[y=P, x=time, col sep=comma, forget plot] {phyto_high_res_filtered.csv};
    \addplot[gray, dotted, thick] table[y=D, x=time, col sep=comma, forget plot] {phyto_high_res_filtered.csv};
    
        \addplot[red, mark=diamond, mark repeat=10, thick] table[y=C, x=time, col sep=comma] {phyto_L2T_1.0.csv};
    \addlegendentry {$C$}
        \addplot[blue, mark=square, mark repeat=10, thick] table[y=N, x=time, col sep=comma] {phyto_L2T_1.0.csv};
    \addlegendentry {$N$}
        \addplot[orange, mark=o, mark repeat=10,  thick] table[y=P, x=time, col sep=comma] {phyto_L2T_1.0.csv};
    \addlegendentry {$P$}
        \addplot[green, mark=diamond, mark repeat=10, thick] table[y=D, x=time, col sep=comma] {phyto_L2T_1.0.csv};
    \addlegendentry {$D$}
\legend{}
\end{axis}
\end{tikzpicture}
\end{subfigure}%

\begin{subfigure}{0.5\textwidth}
\centering
\begin{tikzpicture}
\centering
\begin{axis}[width=\textwidth, legend style={nodes={scale=1.0, transform shape}}, 
        title={RIIA($\mathcal{B}_2$)-VI},
	ylabel near ticks,
	legend pos=outer north east,
	xlabel=$t$,
	ymin=-10,
	ymax=40]
    
    \addplot[gray, dotted, thick] table[y=C, x=time, col sep=comma] {phyto_high_res_filtered.csv};
    \addlegendentry {High-Res}
    \addplot[gray, dotted, thick] table[y=N, x=time, col sep=comma, forget plot] {phyto_high_res_filtered.csv};
    \addplot[gray, dotted, thick] table[y=P, x=time, col sep=comma, forget plot] {phyto_high_res_filtered.csv};
    \addplot[gray, dotted, thick] table[y=D, x=time, col sep=comma, forget plot] {phyto_high_res_filtered.csv};
    
        \addplot[red, mark=diamond, mark repeat=10, thick] table[y=C, x=time, col sep=comma] {phyto_B2T_1.0.csv};
    \addlegendentry {$C$}
        \addplot[blue, mark=square, mark repeat=10, thick] table[y=N, x=time, col sep=comma] {phyto_B2T_1.0.csv};
    \addlegendentry {$N$}
        \addplot[orange, mark=o, mark repeat=10,  thick] table[y=P, x=time, col sep=comma] {phyto_B2T_1.0.csv};
    \addlegendentry {$P$}
        \addplot[green, mark=diamond, mark repeat=10, thick] table[y=D, x=time, col sep=comma] {phyto_B2T_1.0.csv};
    \addlegendentry {$D$}

\end{axis}
\end{tikzpicture}
\end{subfigure}
\caption{The numerical solution to~\eqref{eq:PhytoSys} integrated using RIIA($\mathcal{L}_2$)-VP, RIIA($\mathcal{L}_2$)-VI, and RIIA($\mathcal{B}_2$)-VI, with $k = 1.0$.}
\label{fig:phyto_simulation}
\end{figure}

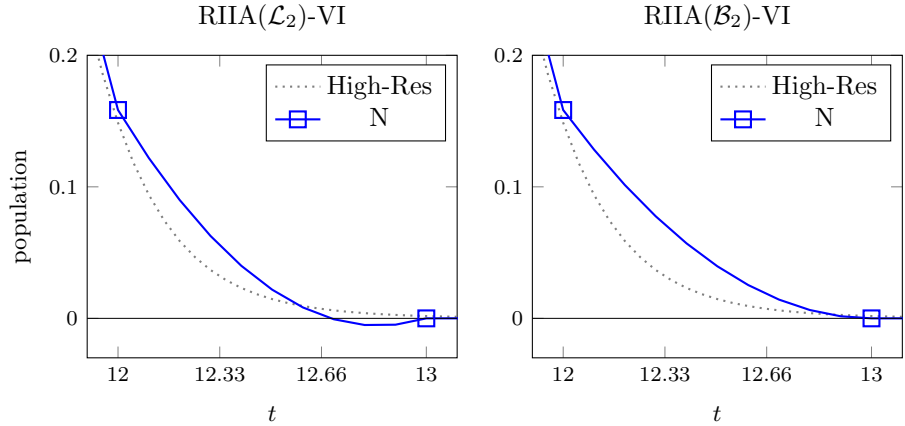
\begin{figure}[ht]
\centering
\begin{subfigure}[t]{0.48\textwidth}
\begin{tikzpicture}
    \begin{axis}[width=\textwidth, legend style={nodes={scale=1.0, transform shape}}, 
        title={RIIA($\mathcal{L}_2$)-VI},
	ylabel near ticks,
	legend pos=north east,
	xlabel=$t$,
	ylabel=population,
	small,
        xmin=11.9,
        xmax=13.1,
        ymin=-0.03,
        ymax=0.2,
        xtick={12.0, 12.33, 12.66, 13.0},
        ytick={0.2, 0.1, 0.0}
    ]
            \addplot[gray, dotted, thick] table[y=N, x=time, col sep=comma] {phyto_high_res_filtered.csv};
        \addlegendentry{High-Res}
        
        \addplot[blue, mark=square, mark repeat=10, mark size=3pt,  thick] table[y=N, x=time, col sep=comma] {phyto_L2T_1.0.csv};
        \addlegendentry {N}

        \draw[ultra thin] (axis cs:\pgfkeysvalueof{/pgfplots/xmin},0) -- (axis cs:\pgfkeysvalueof{/pgfplots/xmax},0);

    \end{axis}
\end{tikzpicture}
\end{subfigure}\hspace{0.04\textwidth}%
\begin{subfigure}[t]{0.48\textwidth}
\centering
\begin{tikzpicture}
    \begin{axis}[width=\textwidth, legend style={nodes={scale=1.0, transform shape}}, 
        title={RIIA($\mathcal{B}_2$)-VI},
	ylabel near ticks,
	legend pos=north east,
	xlabel=$t$,
        small,
        xmin=11.9,
        xmax=13.1,
        ymin=-0.03,
        ymax=0.2,
        xtick={12.0, 12.33, 12.66, 13.0},
        ytick={0.2, 0.1, 0.0}
        ]
        
        \addplot[gray, dotted, thick] table[y=N, x=time, col sep=comma] {phyto_high_res_filtered.csv};
        \addlegendentry {High-Res}
        
        \addplot[blue, mark=square, mark repeat=10, mark size=3pt,  thick] table[y=N, x=time, col sep=comma] {phyto_B2T_1.0.csv};
        \addlegendentry {N}

        \draw[ultra thin] (axis cs:\pgfkeysvalueof{/pgfplots/xmin},0) -- (axis cs:\pgfkeysvalueof{/pgfplots/xmax},0);

    \end{axis}
\end{tikzpicture}
\end{subfigure}
\caption{Bounds violations (left) and uniform nonnegativity (right) of the numerical solution to~\eqref{eq:PhytoSys} integrated using RIIA($\mathcal{L}_2$)-VI and RIIA($\mathcal{B}_2$)-VI, respectively, with $k=1.0$.}
\label{fig:phyto_vios}
\end{figure}

We integrate the system using the 2-stage RadauIIA method. We use a large step size of $k = 1.0$ to approximate the system at final time $20.0$. We utilize the classical Runge-Kutta
method, as well as integrators constrained to preserve nonnegativity using both the Lagrange and Bernstein form of the collocation polynomial. The results are shown
in Figure~\ref{fig:phyto_simulation}.
Without enforcing bounds constraints, the classical Runge-Kutta method predicts a negative stage-value for the nitrogen source. 
This, when coupled with the
other equations, skews the solution of the system leading to instability for all of the populations.  This behavior would resolve with a small enough time step (labelled ``High-Res'' in Figure~\ref{fig:phyto_simulation}).  
When nonnegativity is enforced, however, the numerical solution remains stable
and tracks the exact solution well, even with our large time steps.
Zooming in on the numerical solution computed using the two constrained integrators, we can see the difference between enforcing constraints on
the Lagrange form and the Bernstein form of the collocation polynomial. In Figure~\ref{fig:phyto_vios} (left), we see that when constraints are enforced on the Lagrange form, the
collocation polynomial is constrained at the times corresponding to the collocation nodes, $t = 12.0, t = 12.33$, and $t = 13.0$, but drops negative elsewhere in the interval. However,
when constraints are enforced on the Bernstein form, Figure~\ref{fig:phyto_vios} (right), the collocation polynomial is nonnegative uniformly in time.

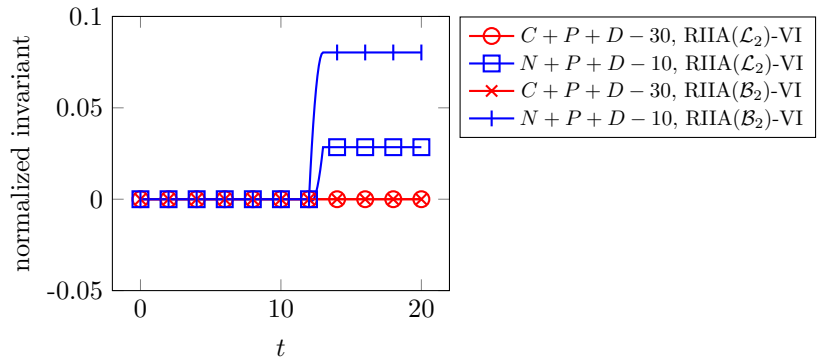
\begin{figure}[ht]
\centering
\begin{tikzpicture}
\centering
\begin{axis}[width=0.5\textwidth, legend style={nodes={scale=0.8, transform shape}}, 
	ylabel near ticks,
	legend pos=outer north east,
	xlabel=$t$,
	ymin=-0.05,
	ymax=0.1,
	ytick={-0.05, 0, 0.05, 0.1},
	yticklabels={-0.05, 0, 0.05, 0.1},
    	ylabel={normalized invariant}]

        \addplot[red, mark=o, mark repeat=20, mark size=3pt, thick] table[y=CPDzero, x=time, col sep=comma] {phyto_L2T_1.0.csv};
    \addlegendentry {$C+P+D - 30$, RIIA($\mathcal{L}_2$)-VI}
    
            \addplot[blue, mark=square, mark repeat=20, mark size=3pt, thick] table[y=NPDzero, x=time, col sep=comma] {phyto_L2T_1.0.csv};
    \addlegendentry {$N+P+D - 10$, RIIA($\mathcal{L}_2$)-VI}
 
         \addplot[red, mark=x, mark repeat=20, mark size=3pt, thick] table[y=CPDzero, x=time, col sep=comma] {phyto_B2T_1.0.csv};
    \addlegendentry {$C+P+D - 30$, RIIA($\mathcal{B}_2$)-VI}
    
            \addplot[blue, mark=|, mark repeat=20, mark size=3pt, thick] table[y=NPDzero, x=time, col sep=comma] {phyto_B2T_1.0.csv};
    \addlegendentry {$N+P+D - 10$, RIIA($\mathcal{B}_2$)-VI}
       
\end{axis}
\end{tikzpicture}
\caption{Linear invariants of~\eqref{eq:PhytoSys} integrated using RIIA($\mathcal{P}_2$)-VI with $k = 1.0$.}
\label{fig:phyto_invariants}
\end{figure}

The solution to~\eqref{eq:PhytoSys} satisfies the linear invariants $C+P+D = 30$ and $N + P + D = 10$. While all classical Runge-Kutta methods will preserve these invariants (even if the solution becomes unstable), the integrators RIIA($\mathcal{P}_2$)-VI, while preserving nonnegativity, have artificially added a small amount of nitrogen to the system. The discrepancy is small, and occurs only where the integrator must act to enforce bounds.
These normalized invariants are plotted in Figure~\ref{fig:phyto_invariants}.
Since the set of feasible solutions that preserves linear invariants remains closed and convex, it should be possible to use a more sophisticated solver that supports equality constraints, but we leave this topic for future study.

\subsection{Heat equation}
\pgfplotstableread[col sep=comma]{heat_conv.csv}\loadedtable

We now consider the heat equation with associated initial and boundary conditions,~\eqref{eq:heat_general}-\eqref{eq:heat_bc}, to study the convergence, performance, and possible constraint violation modes for each of the methods presented.

\subsubsection{Convergence rates}

\begin{figure}[h]
\centering
\begin{subfigure}[t]{0.48\textwidth}
\centering
\begin{tikzpicture}[baseline={(0,0)}]
\begin{axis}[width=\textwidth, legend style={at={(0.5,-0.3)},anchor=north,nodes={scale=0.6, transform shape}, legend columns=2}, 
       title={$\mathcal{P}_1$-RadauIIA(2)},
	xmode=log, 
	ymode=log,
	ylabel near ticks,
	xtick={4,8,16,32, 64, 128}, 
	xticklabels={$2^2$, $2^3$, $2^4$, $2^5$, $2^6$, $2^7$},
	xlabel=N,
    	ylabel={error}
    	]
    \addplot[red,mark=|, mark size=4pt,  thick] table[x=Nspat, y=L1L2TerrorL2, col sep=comma] \loadedtable;
    \addlegendentry[red,mark=oplus,  thick] {$\mathcal{L}_1$-RIIA($\mathcal{L}_2$)-VI};
    
    \addplot[darkgray,mark=o, mark size=4pt, thick] table[x=Nspat, y=L1B2TerrorL2, col sep=comma]\loadedtable;
    \addlegendentry[darkgray,mark=oplus,  thick] {$\mathcal{L}_1$-RIIA($\mathcal{B}_2$)-VI};

    
    \addplot[red,mark=|, mark size=4pt,  thick, dashed, mark options={solid}, forget plot] table[x=Nspat, y=L1L2TerrorH1, col sep=comma] \loadedtable;
    
    \addplot[darkgray,mark=o, mark size=4pt,  thick, dashed,mark options={solid}, forget plot] table[x=Nspat, y=L1B2TerrorH1, col sep=comma]\loadedtable;

    
    \addplot[gray, dashed] table[x=Nspat, y=3Order1, col sep=comma]\loadedtable;
    \addlegendentry[gray, dashed] {$\mathcal{O}(1/N)$};
    
    \addplot[gray] table[x=Nspat, y=Order2, col sep=comma]\loadedtable;
    \addlegendentry[gray] {$\mathcal{O}(1/N^2)$};
\end{axis}

\end{tikzpicture}
\end{subfigure}\hspace{0.04\textwidth}%
\begin{subfigure}[t]{0.48\textwidth}
\centering
\begin{tikzpicture}[baseline={(0,0)}]
\begin{axis}[width=\textwidth, legend style={at={(0.5,-0.3)},anchor=north,nodes={scale=0.6, transform shape}, legend columns=2}, 
        title={$\mathcal{P}_2$-RadauIIA(3)},
	xmode=log, 
	ymode=log,
	ylabel near ticks,
	xtick={4,8,16,32, 64, 128}, 
	xticklabels={$2^2$, $2^3$, $2^4$, $2^5$, $2^6$, $2^7$},
	xlabel=N,
    	ylabel={error}
    	]
    \addplot[red,mark=|, mark size=4pt,  thick] table[x=Nspat, y=L2L3TerrorL2, col sep=comma] \loadedtable;
    \addlegendentry[red,mark=oplus,  thick] {$\mathcal{L}_2$-RIIA($\mathcal{L}_3$)-VI};
    
    \addplot[darkgray,mark=o, mark size=4pt, thick] table[x=Nspat, y=L2B3TerrorL2, col sep=comma] \loadedtable;
    \addlegendentry[darkgray,mark=oplus,  thick] {$\mathcal{L}_2$-RIIA($\mathcal{B}_3$)-VI};
    
    \addplot[brown,mark=triangle, mark size=4pt, thick] table[x=Nspat, y=B2L3TerrorL2, col sep=comma] \loadedtable;
    \addlegendentry[brown,mark=oplus,  thick] {$\mathcal{B}_2$-RIIA($\mathcal{L}_3$)-VI};

    \addplot[blue,mark=square, mark size=4pt,  thick] table[x=Nspat, y=B2B3TerrorL2, col sep=comma] \loadedtable;
    \addlegendentry[blue,mark=oplus,  thick] {$\mathcal{B}_2$-RIIA($\mathcal{B}_3$)-VI};

    
    \addplot[red,mark=|, mark size=4pt,  thick, dashed, mark options={solid}, forget plot] table[x=Nspat, y=L2L3TerrorH1, col sep=comma] \loadedtable;
    
    \addplot[darkgray,mark=o, mark size=4pt,  thick, dashed,mark options={solid}, forget plot] table[x=Nspat, y=L2B3TerrorH1, col sep=comma]\loadedtable;
    
    \addplot[brown,mark=triangle, mark size=4pt,  thick, dashed,mark options={solid}, forget plot] table[x=Nspat, y=B2L3TerrorH1, col sep=comma]\loadedtable;

    \addplot[blue,mark=square, mark size=4pt,  thick, dashed,mark options={solid}, forget plot] table[x=Nspat, y=B2B3TerrorH1, col sep=comma] \loadedtable;
    
    
    \addplot[gray, dashed] table[x=Nspat, y=5Order2, col sep=comma]\loadedtable;
    \addlegendentry[gray, dashed] {$\mathcal{O}(1/N^2)$};    
    
    \addplot[gray] table[x=Nspat, y=Order3, col sep=comma] \loadedtable;
    \addlegendentry[gray] {$\mathcal{O}(1/N^3)$};

\end{axis}
\end{tikzpicture}
\end{subfigure}
\caption{$L^2$ (solid) and $H^1$ (dotted) error in the approximation of $u(1.0)$ using a uniform $N\times N$ mesh and $k = 1/N$.} 
\label{fig:heat_conv_mismatch}
\end{figure}
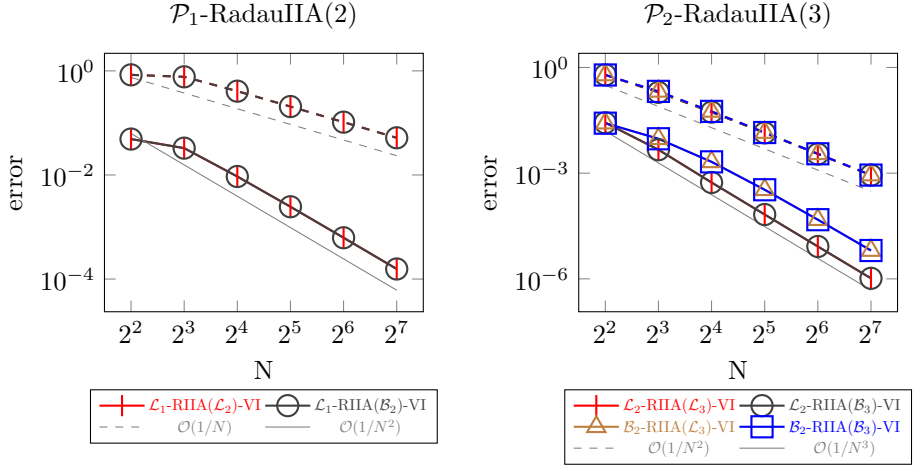

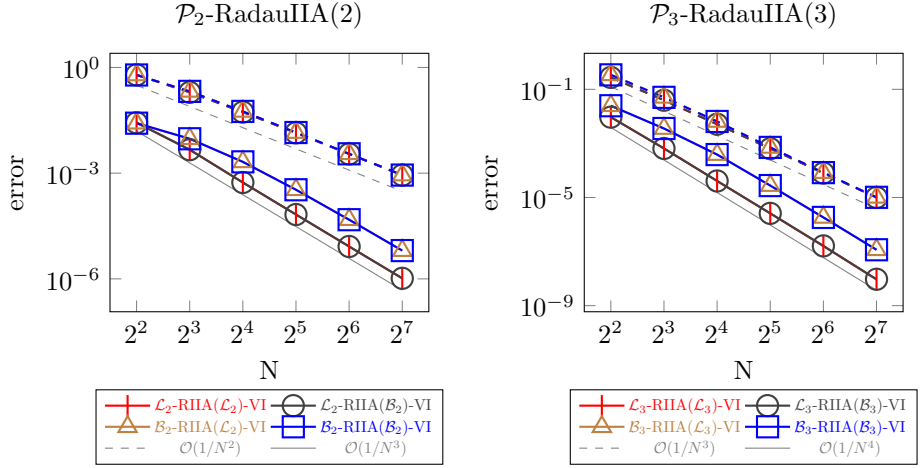
\begin{figure}[ht]
\centering
\begin{subfigure}[t]{0.48\textwidth}
  \centering
\begin{tikzpicture}
\begin{axis}[width=\textwidth, legend style={at={(0.5,-0.3)},anchor=north,nodes={scale=0.6, transform shape}, legend columns=2},
        title={$\mathcal{P}_2$-RadauIIA(2)},
	xmode=log, 
	ymode=log,
	xtick={4,8,16,32, 64, 128}, 
	xticklabels={$2^2$, $2^3$, $2^4$, $2^5$, $2^6$, $2^7$},
	xlabel=N,
    	ylabel={error}]
    \addplot[red,mark=|, mark size=4pt,  thick] table[x=Nspat, y=L2L2TerrorL2, col sep=comma]\loadedtable;
    \addlegendentry[red,mark=oplus,  thick] {$\mathcal{L}_2$-RIIA($\mathcal{L}_2$)-VI};
    
    \addplot[darkgray,mark=o, mark size=4pt, thick] table[x=Nspat, y=L2B2TerrorL2, col sep=comma]\loadedtable;
    \addlegendentry[darkgray,mark=oplus,  thick] {$\mathcal{L}_2$-RIIA($\mathcal{B}_2$)-VI};
    
    \addplot[brown,mark=triangle, mark size=4pt, thick] table[x=Nspat, y=B2L2TerrorL2, col sep=comma]\loadedtable;
    \addlegendentry[brown,mark=oplus,  thick] {$\mathcal{B}_2$-RIIA($\mathcal{L}_2$)-VI};

    \addplot[blue,mark=square, mark size=4pt,  thick] table[x=Nspat, y=B2B2TerrorL2, col sep=comma]\loadedtable;
    \addlegendentry[blue,mark=oplus,  thick] {$\mathcal{B}_2$-RIIA($\mathcal{B}_2$)-VI};

    
    \addplot[red,mark=|, mark size=4pt,  thick, dashed, mark options={solid}, forget plot] table[x=Nspat, y=L2L2TerrorH1, col sep=comma] \loadedtable;
    
    \addplot[darkgray,mark=o, mark size=4pt,  thick, dashed,mark options={solid}, forget plot] table[x=Nspat, y=L2B2TerrorH1, col sep=comma] \loadedtable;
    
    \addplot[brown,mark=triangle, mark size=4pt,  thick, dashed,mark options={solid}, forget plot] table[x=Nspat, y=B2L2TerrorH1, col sep=comma]\loadedtable;

    \addplot[blue,mark=square, mark size=4pt,  thick, dashed,mark options={solid}, forget plot] table[x=Nspat, y=B2B2TerrorH1, col sep=comma] \loadedtable;
    
    
    \addplot[gray, dashed] table[x=Nspat, y=5Order2, col sep=comma] \loadedtable;
    \addlegendentry[gray, dashed] {$\mathcal{O}(1/N^2)$};
    
    \addplot[gray] table[x=Nspat, y=Order3, col sep=comma] \loadedtable;
    \addlegendentry[gray] {$\mathcal{O}(1/N^3)$};
\end{axis}

\end{tikzpicture}
\end{subfigure}\hspace{0.04\textwidth}%
\begin{subfigure}[t]{0.48\textwidth}
  \centering
\begin{tikzpicture}
\centering
\begin{axis}[width=\textwidth, legend style={at={(0.5,-0.3)},anchor=north,nodes={scale=0.6, transform shape}, legend columns=2},
        title={$\mathcal{P}_3$-RadauIIA(3)},
	xmode=log, 
	ymode=log,
	xtick={4,8,16,32, 64, 128}, 
	xticklabels={$2^2$, $2^3$, $2^4$, $2^5$, $2^6$, $2^7$},
	xlabel=N,
    	ylabel={error}]
    	
    \addplot[red,mark=|, mark size=4pt,  thick] table[x=Nspat, y=L3L3TerrorL2, col sep=comma] \loadedtable;
    \addlegendentry[red,mark=oplus,  thick] {$\mathcal{L}_3$-RIIA($\mathcal{L}_3$)-VI};
    
    \addplot[darkgray,mark=o, mark size=4pt, thick] table[x=Nspat, y=L3B3TerrorL2, col sep=comma] \loadedtable;
    \addlegendentry[darkgray,mark=oplus,  thick] {$\mathcal{L}_3$-RIIA($\mathcal{B}_3$)-VI};
    
    \addplot[brown,mark=triangle, mark size=4pt, thick] table[x=Nspat, y=B3L3TerrorL2, col sep=comma] \loadedtable;
    \addlegendentry[brown,mark=oplus,  thick] {$\mathcal{B}_3$-RIIA($\mathcal{L}_3$)-VI};

    \addplot[blue,mark=square, mark size=4pt,  thick] table[x=Nspat, y=B3B3TerrorL2, col sep=comma] \loadedtable;
    \addlegendentry[blue,mark=oplus,  thick] {$\mathcal{B}_3$-RIIA($\mathcal{B}_3$)-VI};

    
    \addplot[red,mark=|, mark size=4pt,  thick, dashed, mark options={solid}, forget plot] table[x=Nspat, y=L3L3TerrorH1, col sep=comma] \loadedtable;
    
    \addplot[darkgray,mark=o, mark size=4pt,  thick, dashed,mark options={solid}, forget plot] table[x=Nspat, y=L3B3TerrorH1, col sep=comma]\loadedtable;
    
    \addplot[brown,mark=triangle, mark size=4pt,  thick, dashed,mark options={solid}, forget plot] table[x=Nspat, y=B3L3TerrorH1, col sep=comma] \loadedtable;

    \addplot[blue,mark=square, mark size=4pt,  thick, dashed,mark options={solid}, forget plot] table[x=Nspat, y=B3B3TerrorH1, col sep=comma]\loadedtable;
    
    
    \addplot[gray, dashed] table[x=Nspat, y=8Order3, col sep=comma]\loadedtable;
    \addlegendentry[gray, dashed] {$\mathcal{O}(1/N^3)$};

    \addplot[gray] table[x=Nspat, y=Order4, col sep=comma]\loadedtable;
    \addlegendentry[gray] {$\mathcal{O}(1/N^4)$};
\end{axis}
\end{tikzpicture}
\end{subfigure}
\caption{$L^2$ (solid) and $H^1$ (dotted) error in the approximation of $u(1.0)$ using a uniform $N\times N$ mesh and $k = 1/N$.} 
\label{fig:heat_conv_match}
\end{figure}

Let $\Omega = [0,1]\times[0,1]$. In order to test the accuracy and performance of the proposed approaches, we pick $f$ and $g$ in \eqref{eq:heat_general} such that the exact solution is given by
\begin{equation}\label{eq:heat_conv_exact}
u = e^{-t}\cos^2(2\pi x)\sin^2(2\pi y).
\end{equation}
We use as the initial condition an appropriately bounds-constrained projection of $u$ onto the spatial finite element space. 
We approximate the solution at $t = 1.0$ using an $N\times N$ square mesh, with each square subdivided into two triangles. 
We use the step size $k = 1 / N$, and advance the system using spatial elements and collocation methods of various polynomial degrees. The $L^2$ and $H^1$ errors 
of the approximation at $t = 1.0$ are plotted in Figures~\ref{fig:heat_conv_mismatch} and \ref{fig:heat_conv_match} for matching and mismatching spatial and temporal polynomial degrees, respectively. We observe optimal convergence rates as $k$ decreases. It should be noted that, for this example, the main difference in the accuracy of the approximation is in the choice 
of spatial finite element.  Since the feasible set for the Bernstein basis contains strictly fewer functions than that for Lagrange, we expect it to have a worse best approximation but hope for the same order of accuracy.
For tests with only linear finite elements in space, we report on the choice of the Lagrange basis in space since it and Bernstein coincide in this case.

\subsubsection{Nonlinear iterations}

Turning to solver performance, we again choose (\ref{eq:heat_conv_exact}) as the exact solution of (\ref{eq:heat_general}), and approximate the solution using various spatial and temporal polynomial degrees. We divide the unit square into an 
$N\times N$ square mesh, and subdivide each square into two triangles. We use either Lagrange or Bernstein finite elements in space and timestep the system using both RadauIIA(2) and RadauIIA(3) with step size $k = 1 / N$. Note that for 
each of the associated unconstrained methods, the resulting system is linear, so each will require only a single Newton iteration per time step. The average number of Newton iterations per timestep is recorded in Table~\ref{tab:heat_equation_aveIts}. We see satisfactory performance for all of the methods tested.  While enforcing bounds constraints 
requires the solution of a more challenging constrained optimization problem, the number of additional Newton iterations per time step is small, and does not grow with the size of the mesh.
Although deploying efficient linear solvers on the reduced space Jacobian remains an open question to us, we view these low iteration counts of the Newton VI solver as a very positive empirical result.

\begin{table}[ht]
\small
\begin{center}
\resizebox{\textwidth}{!}{%
\pgfplotstabletypeset[
	every head row/.style={before row=\hline, after row=\hline}, 
	display columns/0/.style={column name=\rotatebox{-90}{N}},
	assign column name/.style={/pgfplots/table/column name={\rotatebox{90}{#1}}},
	every last row/.style={after row=\hline},
	every first column/.style={column type/.add={|}{}},
	every last column/.style={column type/.add={}{|}},
	columns={Nspat, 
						L1L2TaveIts,
			L1B2TaveIts,
						L2L2TaveIts,
			L2B2TaveIts,
			B2L2TaveIts, 
			B2B2TaveIts,
						L2L3TaveIts,
			L2B3TaveIts,
			B2L3TaveIts, 
			B2B3TaveIts,
						L3L3TaveIts,
			L3B3TaveIts,
			B3L3TaveIts, 
			B3B3TaveIts
			 },
	precision=1,
	fixed zerofill,
	columns/Nspat/.style={column name={$N$}, column type/.add={}{|}, 	precision=0, fixed zerofill},
		columns/L1L2TaveIts/.style={column name={ $\mathcal{L}_1$-RIIA($\mathcal{L}_2$)-VI}},
	columns/L1B2TaveIts/.style={column name={ $\mathcal{L}_1$-RIIA($\mathcal{B}_2$)-VI}, column type/.add={}{|}},
		columns/L2L2TaveIts/.style={column name={ $\mathcal{L}_2$-RIIA($\mathcal{L}_2$)-VI}},
	columns/L2B2TaveIts/.style={column name={ $\mathcal{L}_2$-RIIA($\mathcal{B}_2$)-VI\;}},
	columns/B2L2TaveIts/.style={column name={ $\mathcal{B}_2$-RIIA($\mathcal{L}_2$)-VI}},
	columns/B2B2TaveIts/.style={column name={ $\mathcal{B}_2$-RIIA($\mathcal{B}_2$)-VI}, column type/.add={}{|}},
		columns/L2L3TaveIts/.style={column name={ $\mathcal{L}_2$-RIIA($\mathcal{L}_3$)-VI}},
	columns/L2B3TaveIts/.style={column name={ $\mathcal{L}_2$-RIIA($\mathcal{B}_3$)-VI}},
	columns/B2L3TaveIts/.style={column name={ $\mathcal{B}_2$-RIIA($\mathcal{L}_3$)-VI}},
	columns/B2B3TaveIts/.style={column name={ $\mathcal{B}_2$-RIIA($\mathcal{B}_3$)-VI}, column type/.add={}{|}},
		columns/L3L3TaveIts/.style={column name={ $\mathcal{L}_3$-RIIA($\mathcal{L}_3$)-VI}},
	columns/L3B3TaveIts/.style={column name={ $\mathcal{L}_3$-RIIA($\mathcal{B}_3$)-VI}},
	columns/B3L3TaveIts/.style={column name={ $\mathcal{B}_3$-RIIA($\mathcal{L}_3$)-VI}},
	columns/B3B3TaveIts/.style={column name={ $\mathcal{B}_3$-RIIA($\mathcal{B}_3$)-VI}}
	]\loadedtable
}
\end{center}

\caption{Average number of nonlinear iterations required to integrate~\eqref{eq:heat_general} over $[0, 1]$.}
\label{tab:heat_equation_aveIts}
\end{table}

\subsubsection{Constraint violation modes}

We now give an example to highlight the different ways in which each method of enforcing bounds constraints may fail. Any Runge-Kutta method which is equivalent to a collocation method produces a continuous approximation to the system at hand by means of the collocation polynomial.
The differences between some of the methods presented may only be realized when examining the solution between the discrete timesteps and collocation nodes. In particular, if Bernstein finite elements are used in space, the system will be uniformly bounds-constrained 
at the discrete times, and at the collocation nodes, regardless of the temporal polynomial basis chosen.  Between the collocation nodes, however, using the Lagrange form of the collocation polynomial will not necessarily guarantee a bounds-constrained approximation, while using the Bernstein form will guarantee uniform bounds constraint adherence.
To illustrate this, consider again the heat equation. Choose the initial condition, forcing function, and boundary conditions of~\eqref{eq:heat_general} such that the exact solution is given by 
\begin{equation}\label{eq:heat_vios_exact}
\frac{1}{4}\left(1 - \tanh\left(\frac{0.15 - \sqrt{(x - 0.5)^2 + (y - 0.5)^2)}}{0.015}\right)\right) \left(1 + \tanh\left(75t - 6.0\right)\right).
\end{equation}
Note that the exact solution is bounded below by $0$ uniformly in time and space.

\begin{figure}[ht]
\centering
\begin{subfigure}{.2\textwidth}-
\centering
\rotatebox{90}{\tiny $\mathcal{L}_2$-RIIA($\mathcal{L}_2$)-VP}%
\includegraphics[width=\textwidth] {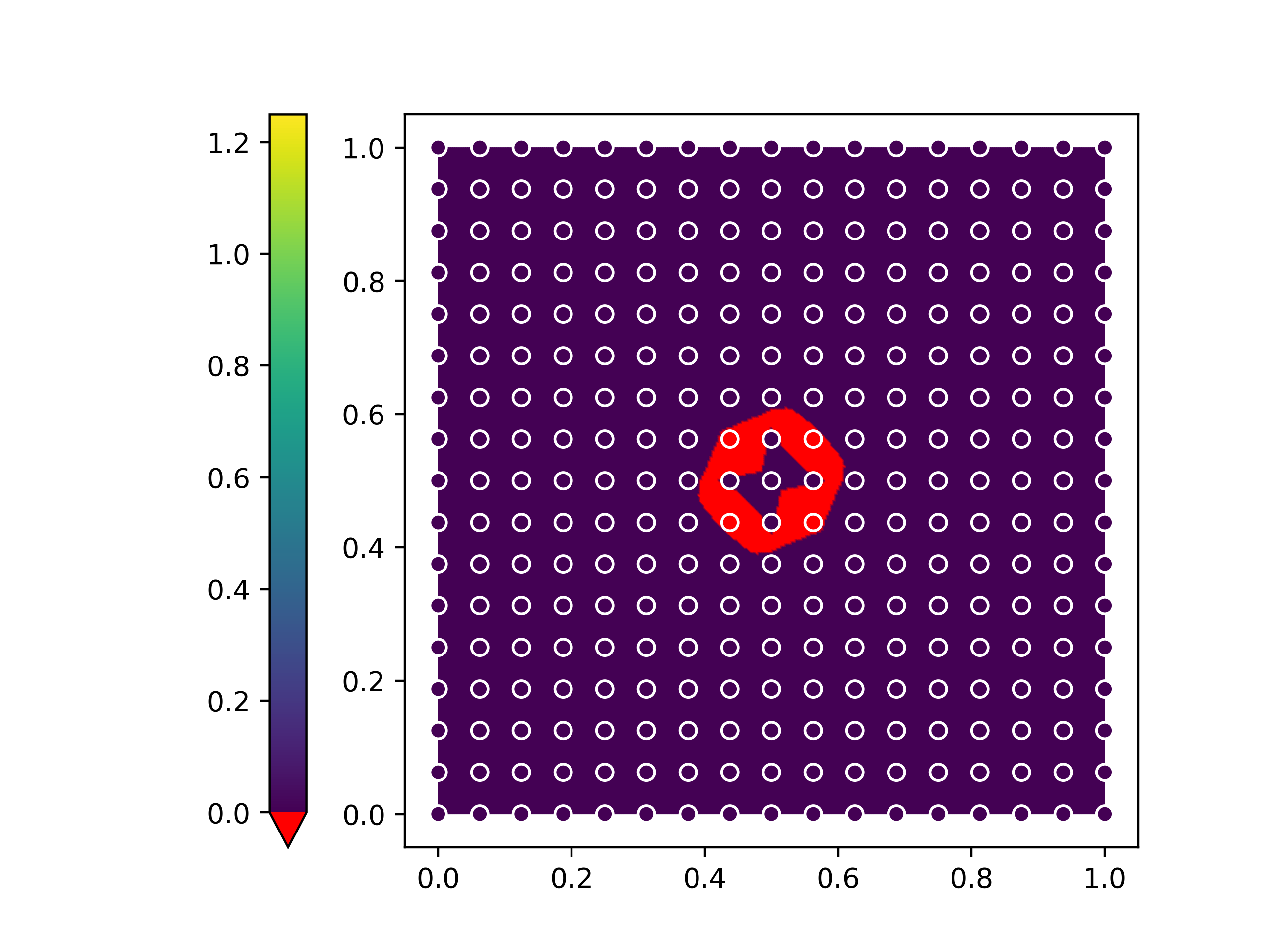}
\end{subfigure}%
\begin{subfigure}{.2\textwidth}
\centering
\includegraphics[width=\textwidth] {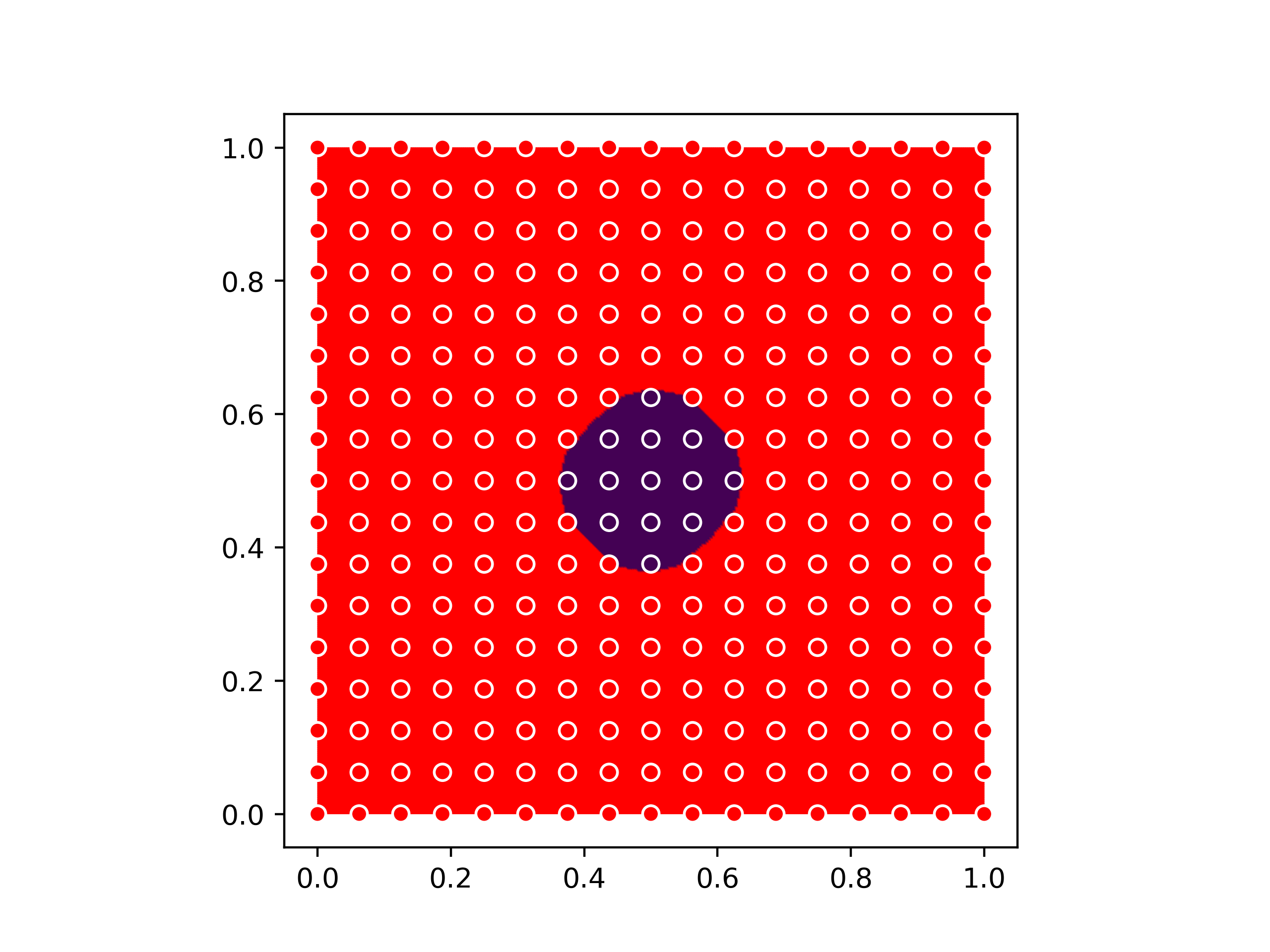}
\end{subfigure}%
\begin{subfigure}{.2\textwidth}
\centering
\includegraphics[width=\textwidth] {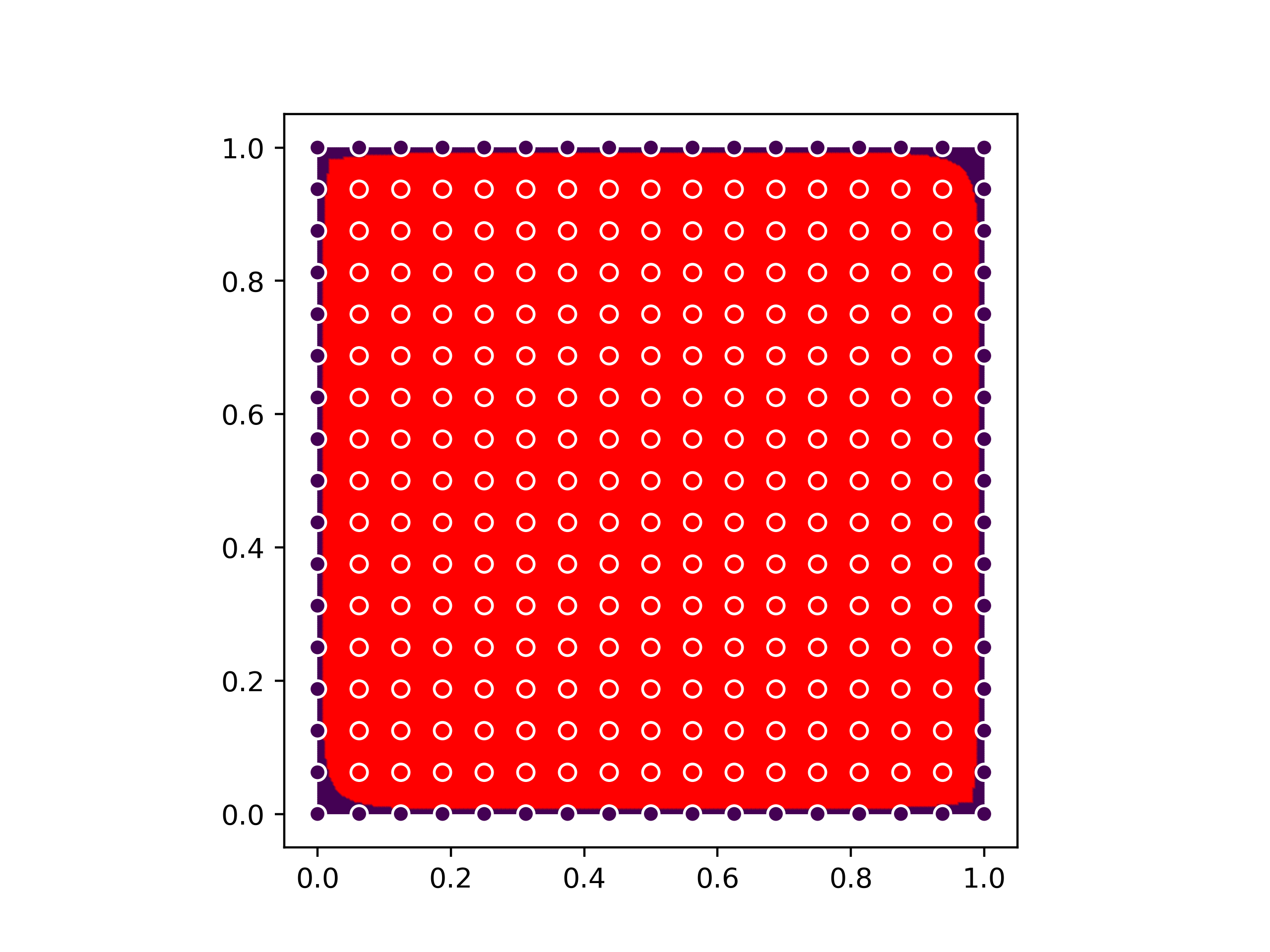}
\end{subfigure}%
\begin{subfigure}{.2\textwidth}
\centering
\includegraphics[width=\textwidth] {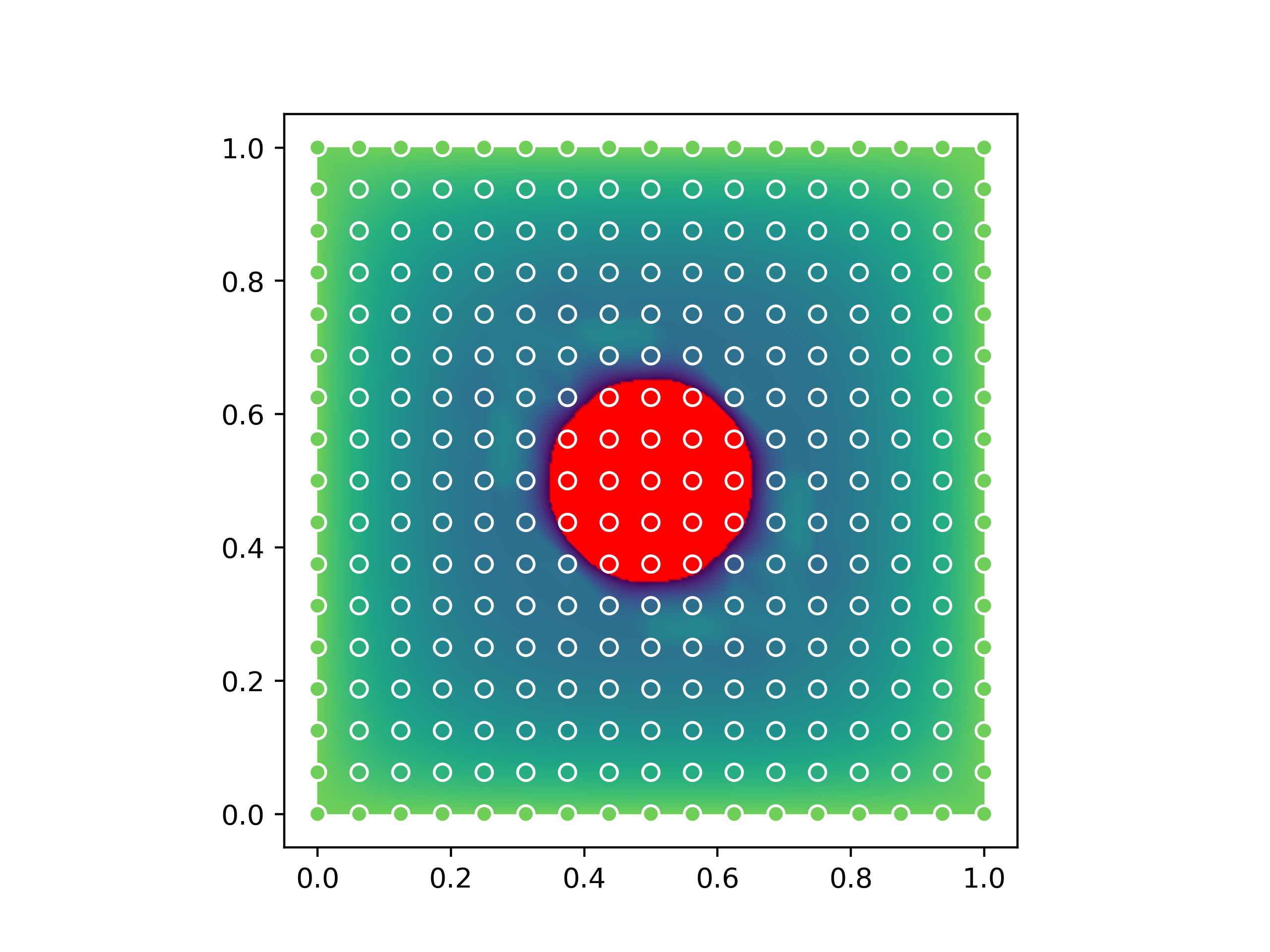}
\end{subfigure}%
\begin{subfigure}{.2\textwidth}
\centering
\includegraphics[width=\textwidth] {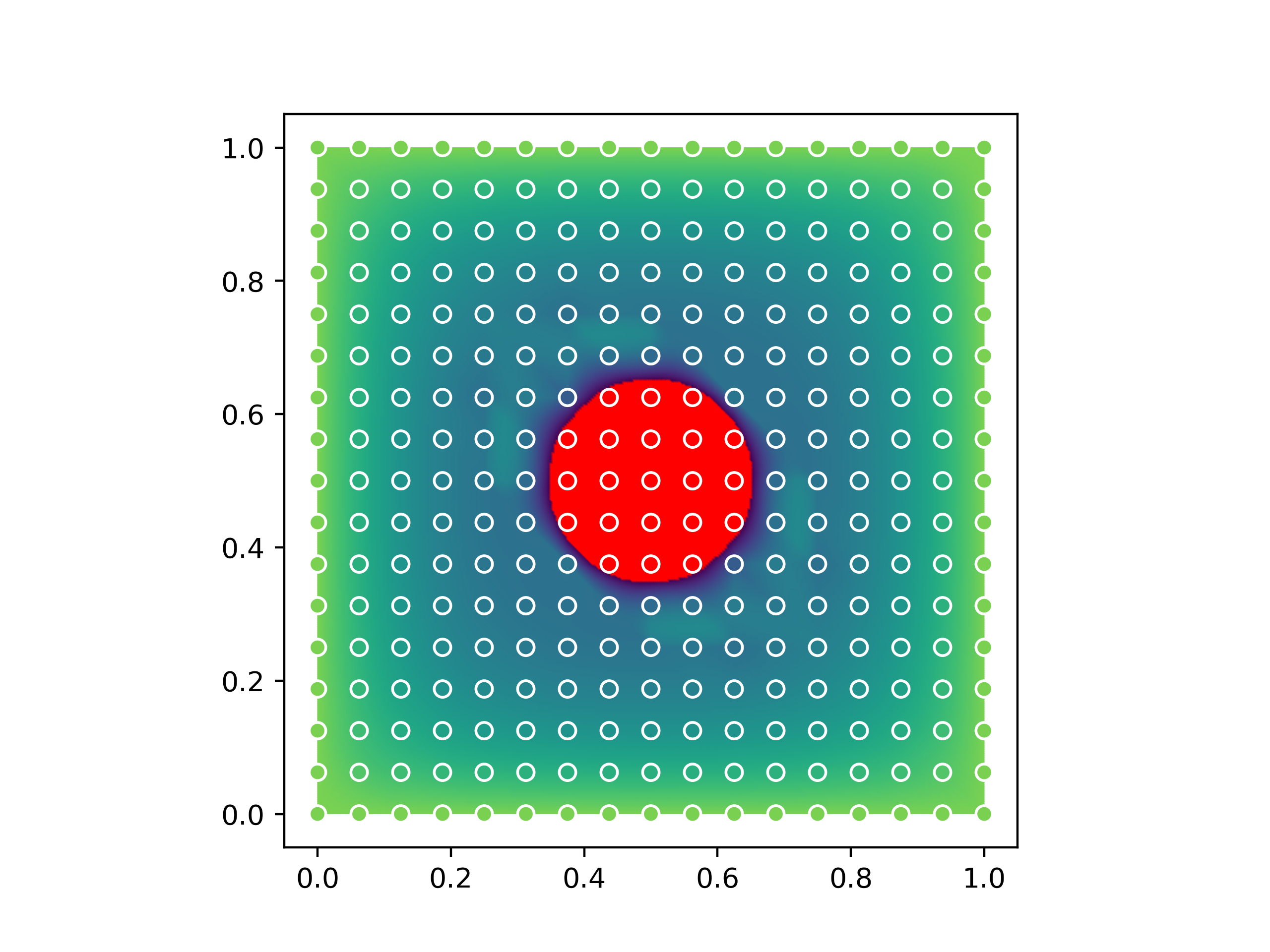}
\end{subfigure}%

\begin{subfigure}{.2\textwidth}
\centering
\rotatebox{90}{\tiny $\mathcal{L}_2$-RIIA($\mathcal{L}_2$)-VI}%
\includegraphics[width=\textwidth] {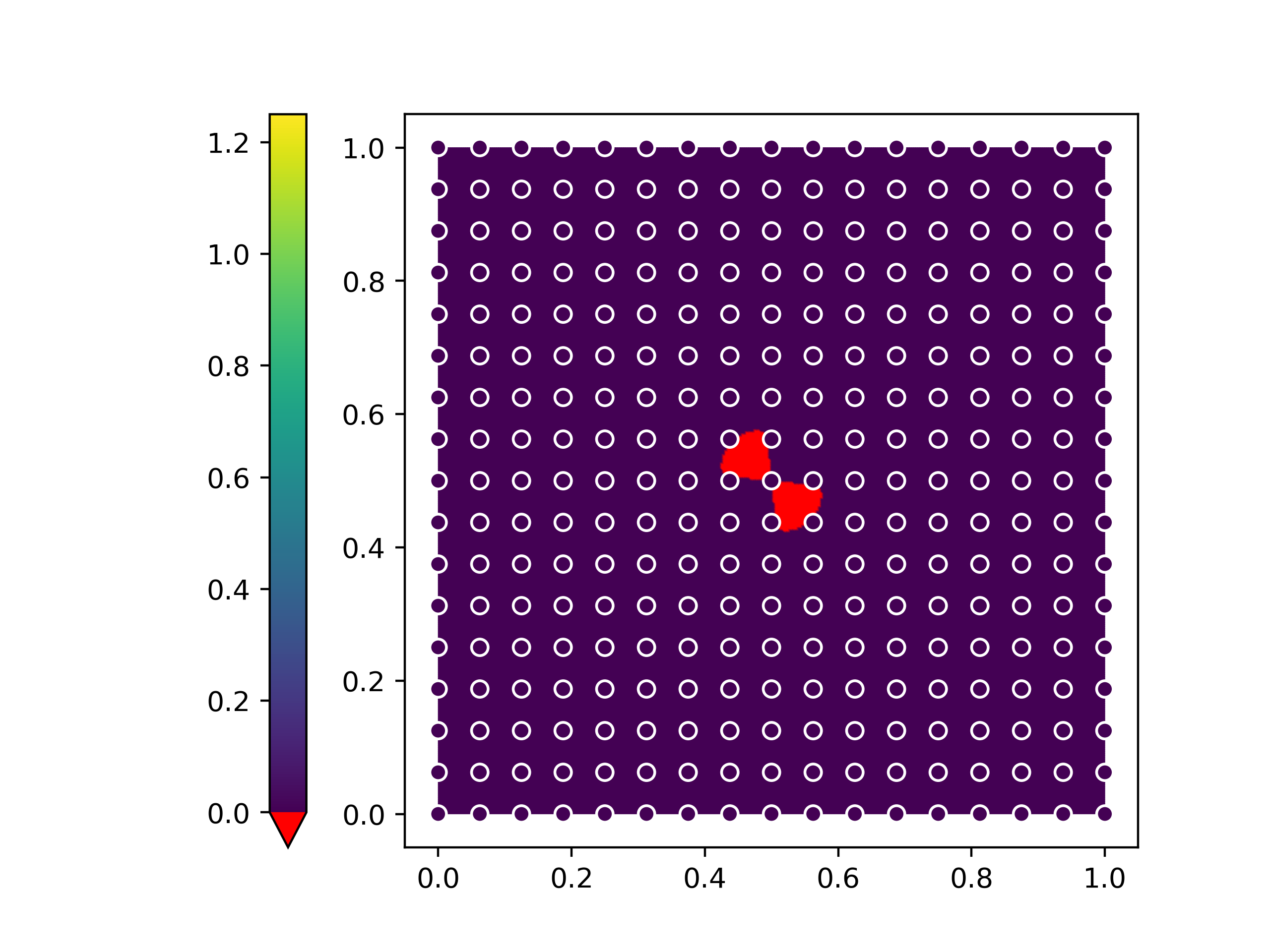}
\end{subfigure}%
\begin{subfigure}{.2\textwidth}
\centering
\includegraphics[width=\textwidth] {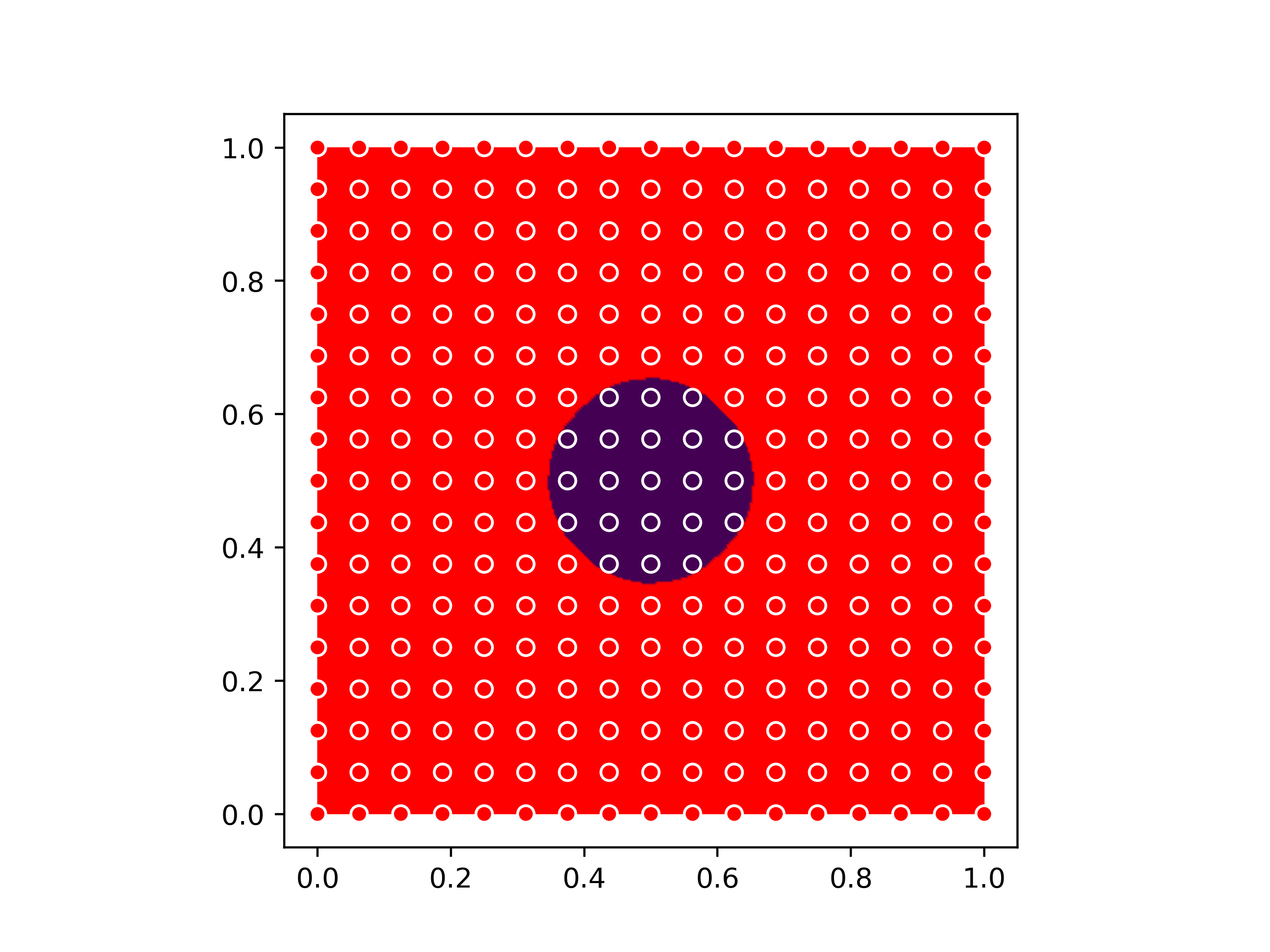}
\end{subfigure}%
\begin{subfigure}{.2\textwidth}
\centering
\includegraphics[width=\textwidth] {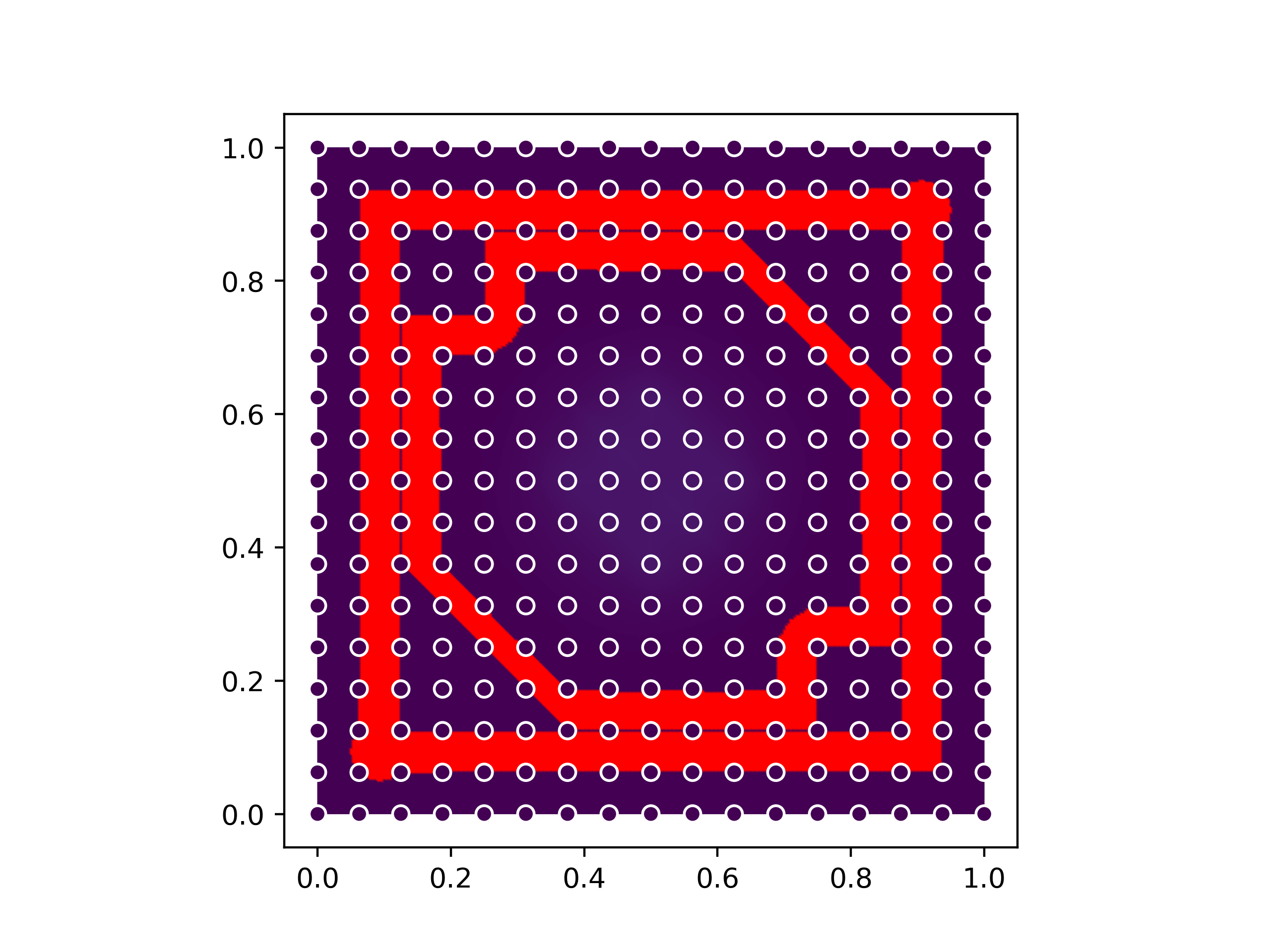}
\end{subfigure}%
\begin{subfigure}{.2\textwidth}
\centering
\includegraphics[width=\textwidth] {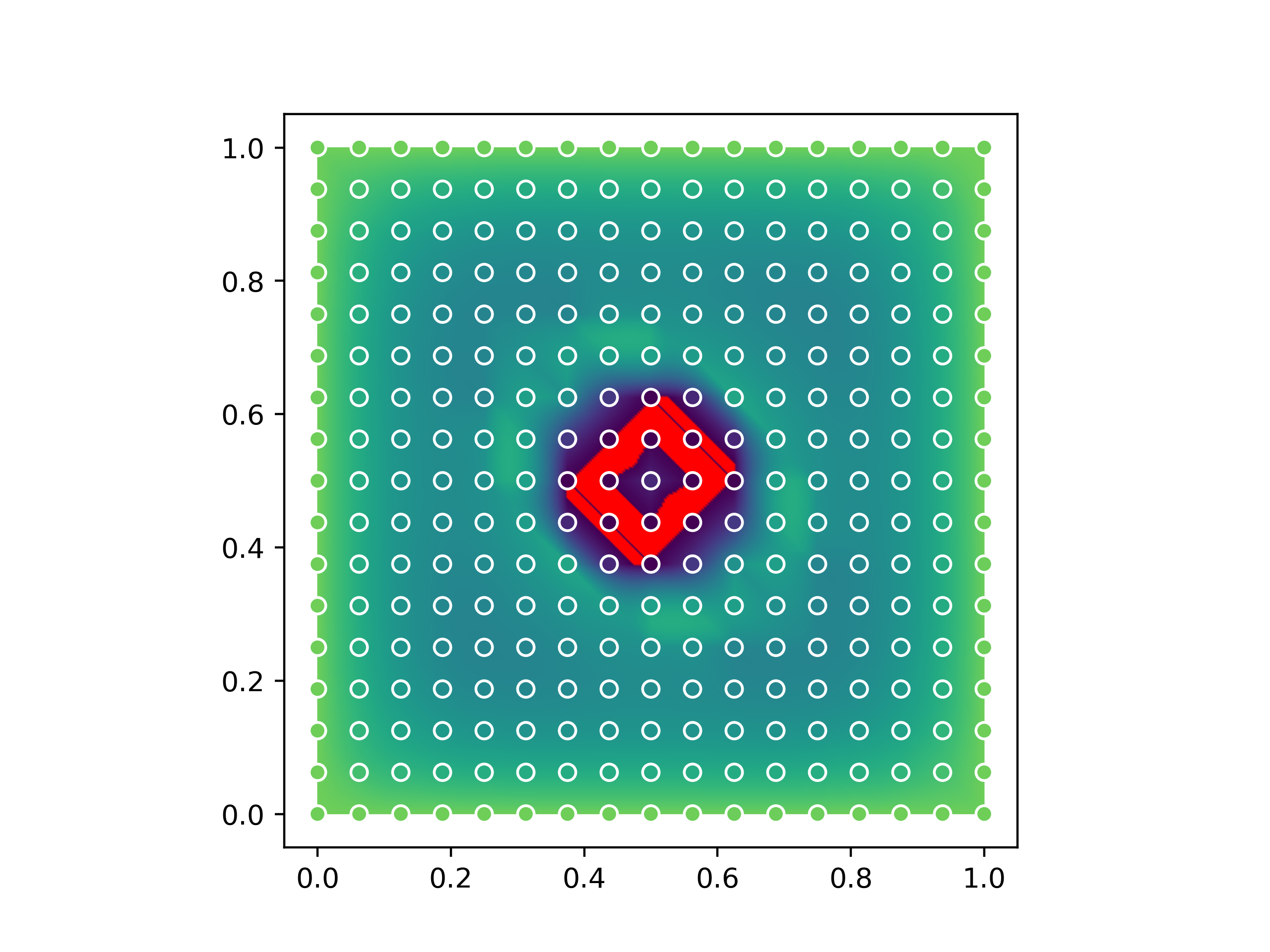}
\end{subfigure}%
\begin{subfigure}{.2\textwidth}
\centering
\includegraphics[width=\textwidth] {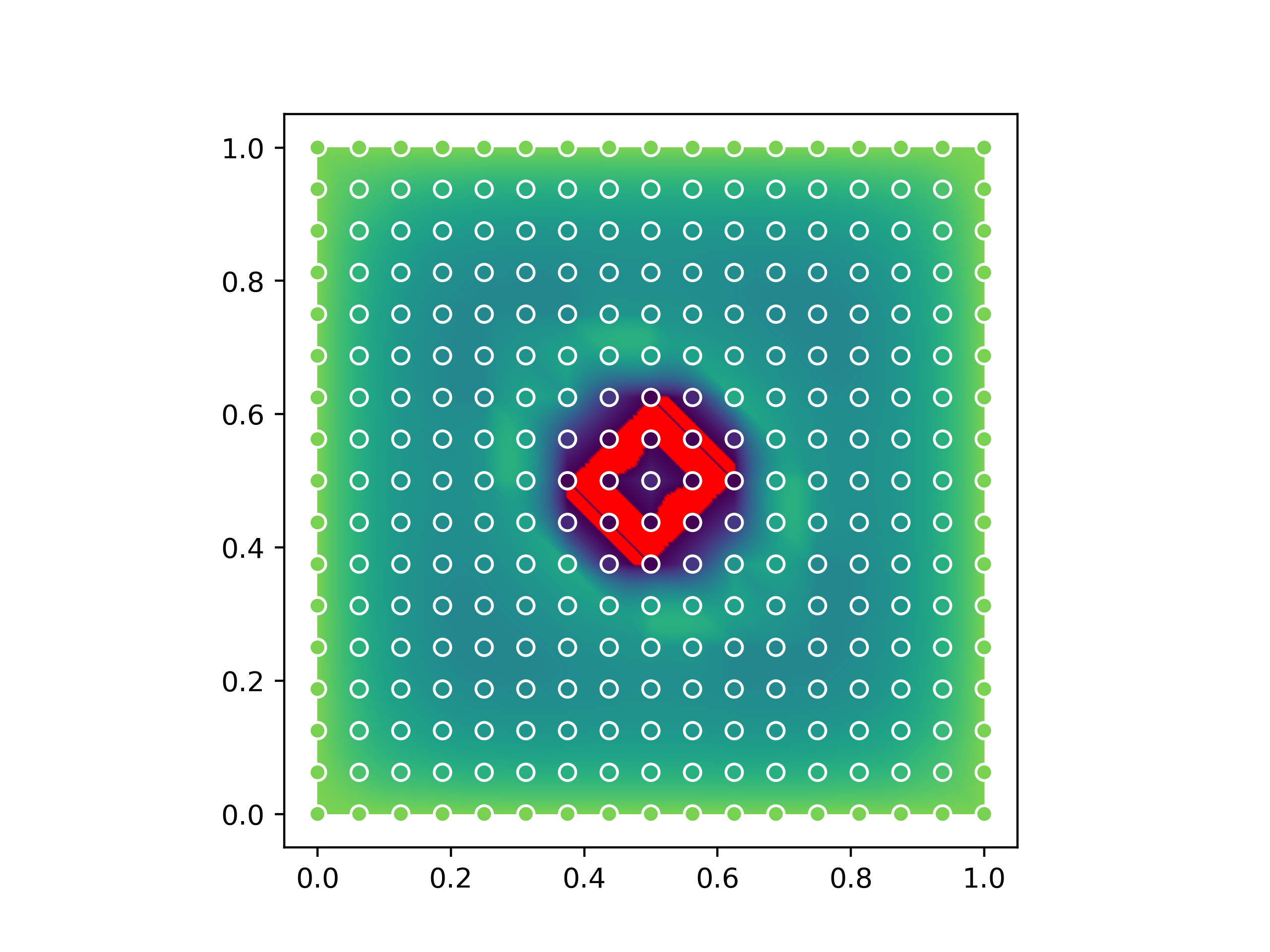}
\end{subfigure}%

\begin{subfigure}{.2\textwidth}
\centering
\rotatebox{90}{\tiny $\mathcal{B}_2$-RIIA($\mathcal{L}_2$)-VI}%
\includegraphics[width=\textwidth] {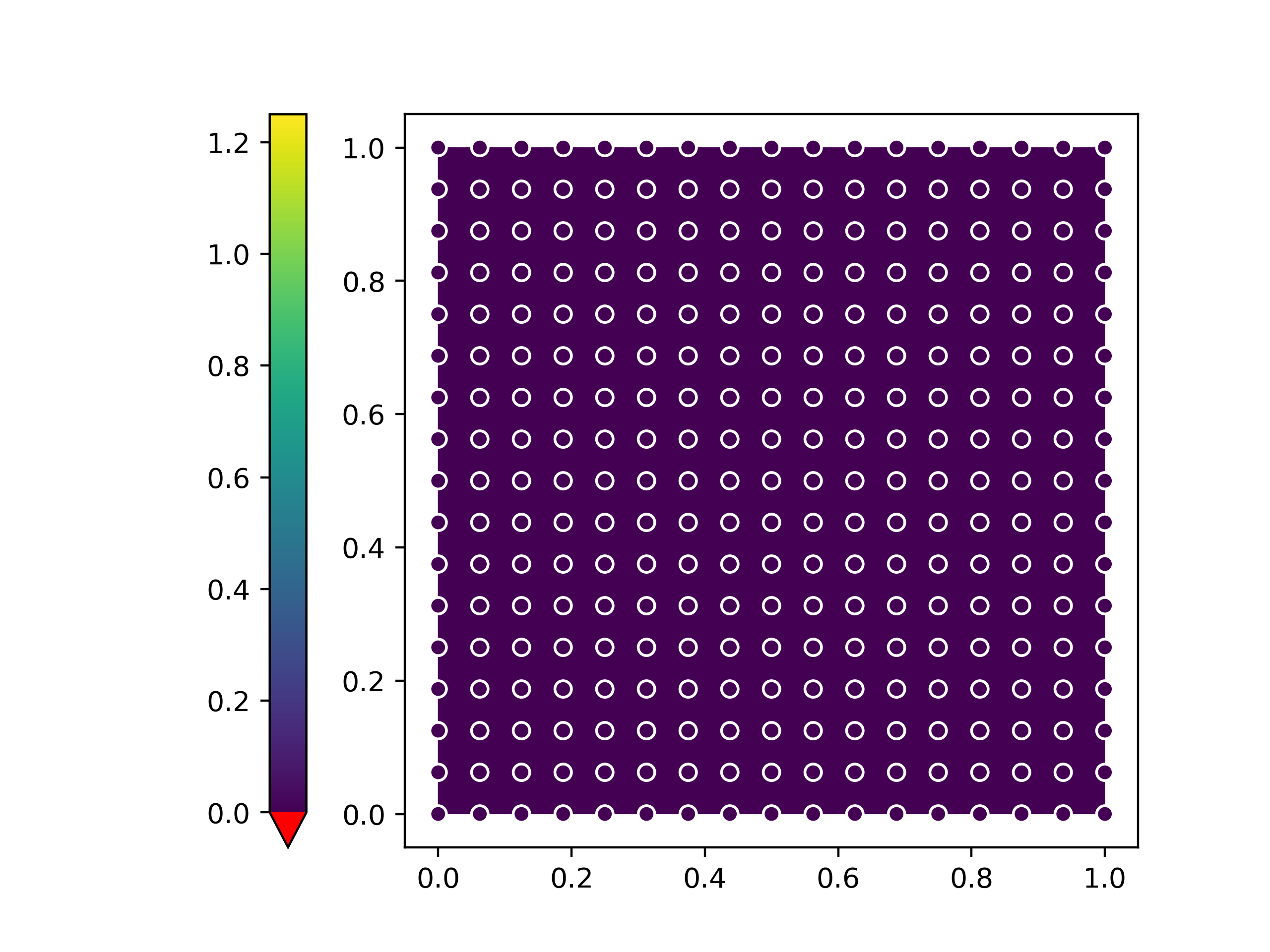}
\end{subfigure}%
\begin{subfigure}{.2\textwidth}
\centering
\includegraphics[width=\textwidth] {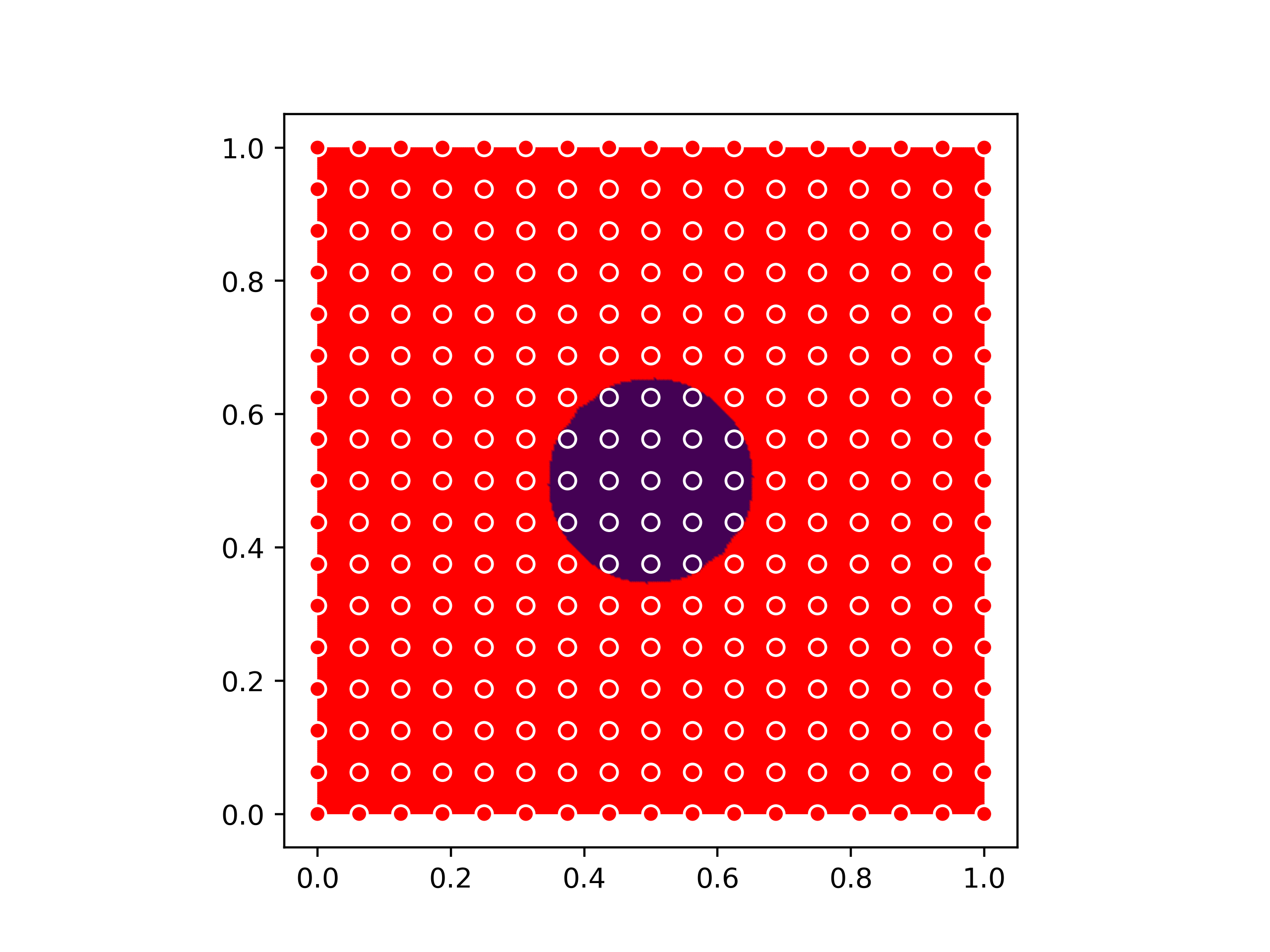}
\end{subfigure}%
\begin{subfigure}{.2\textwidth}
\centering
\includegraphics[width=\textwidth] {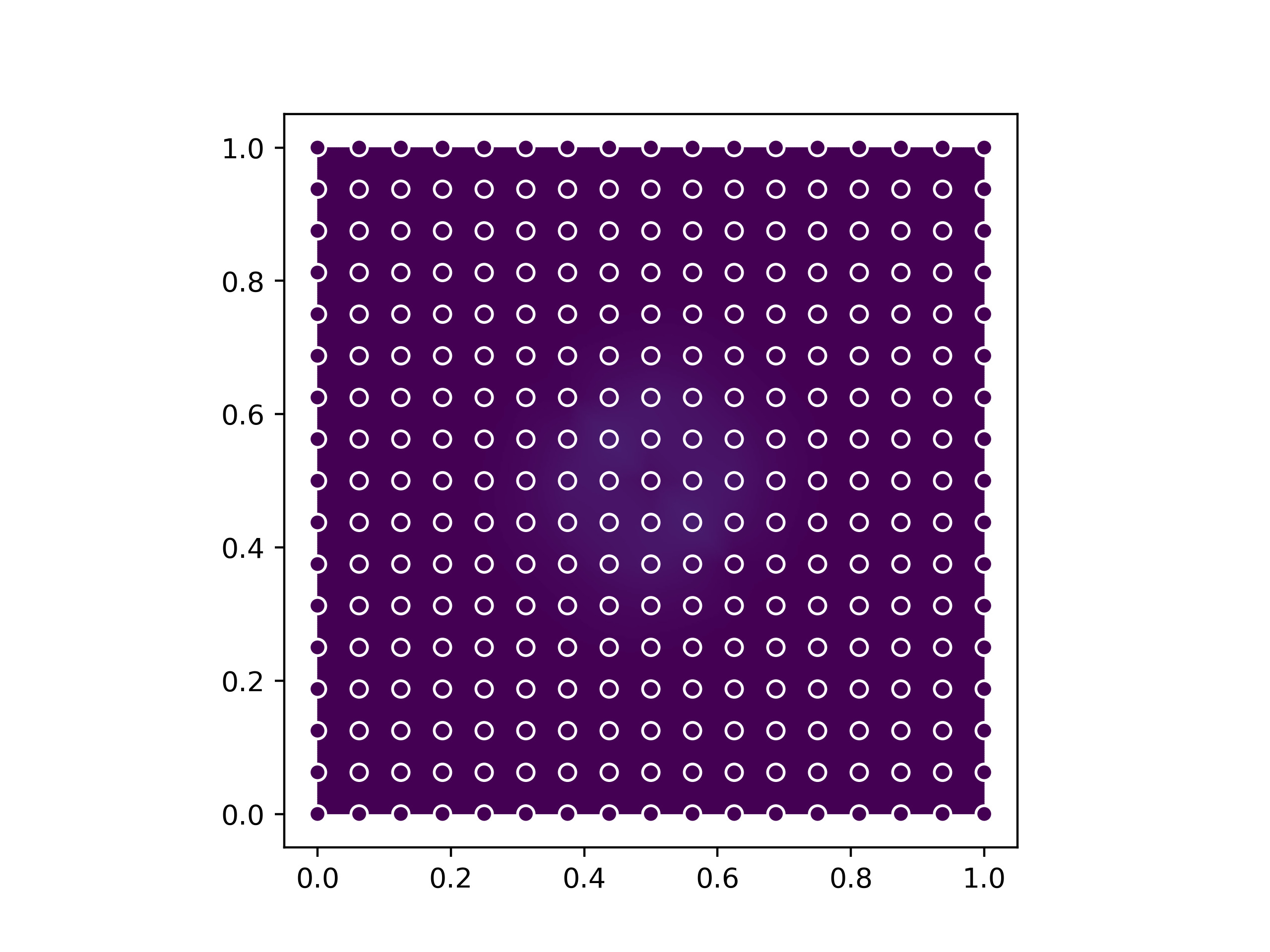}
\end{subfigure}%
\begin{subfigure}{.2\textwidth}
\centering
\includegraphics[width=\textwidth] {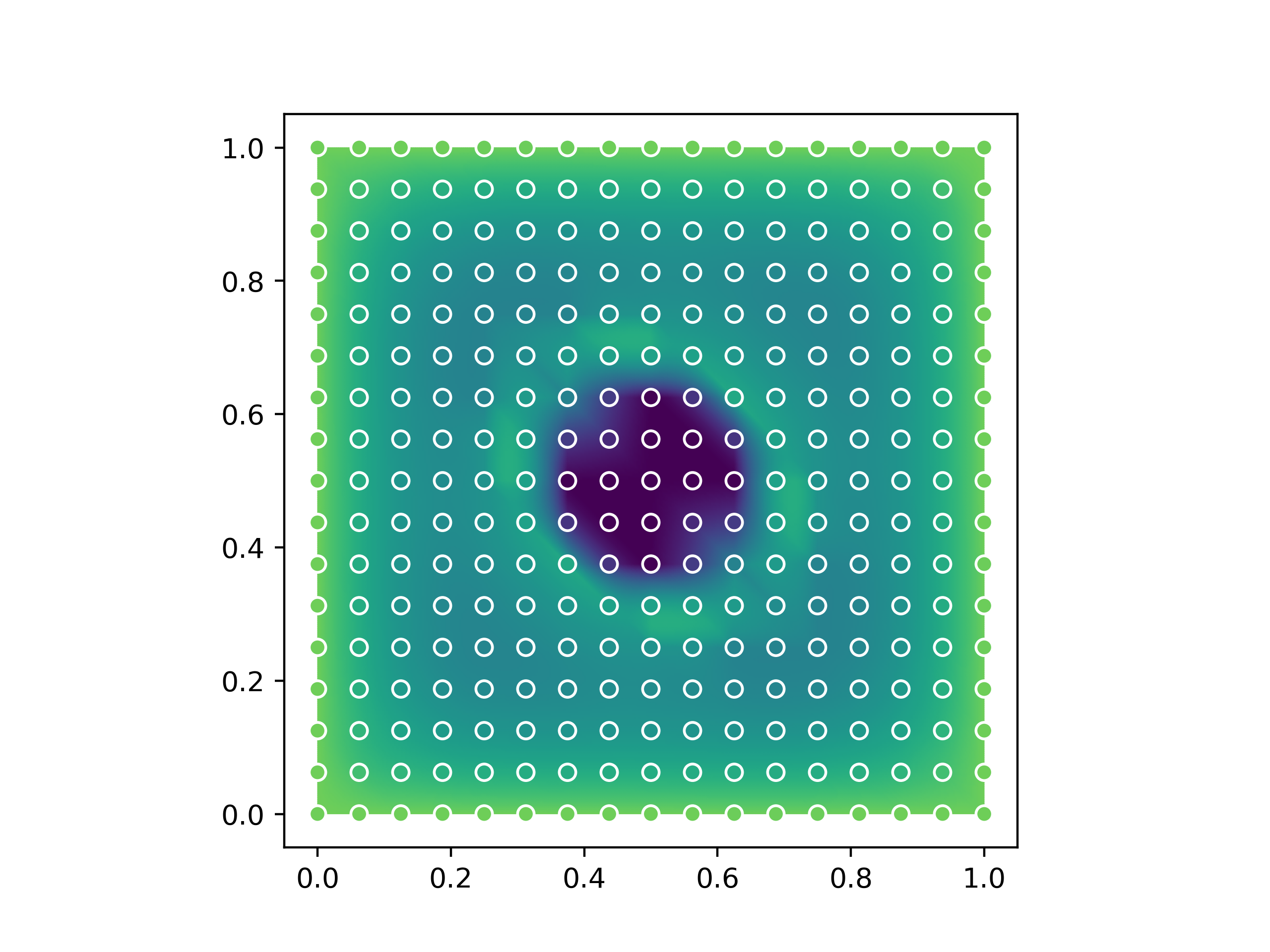}
\end{subfigure}%
\begin{subfigure}{.2\textwidth}
\centering
\includegraphics[width=\textwidth] {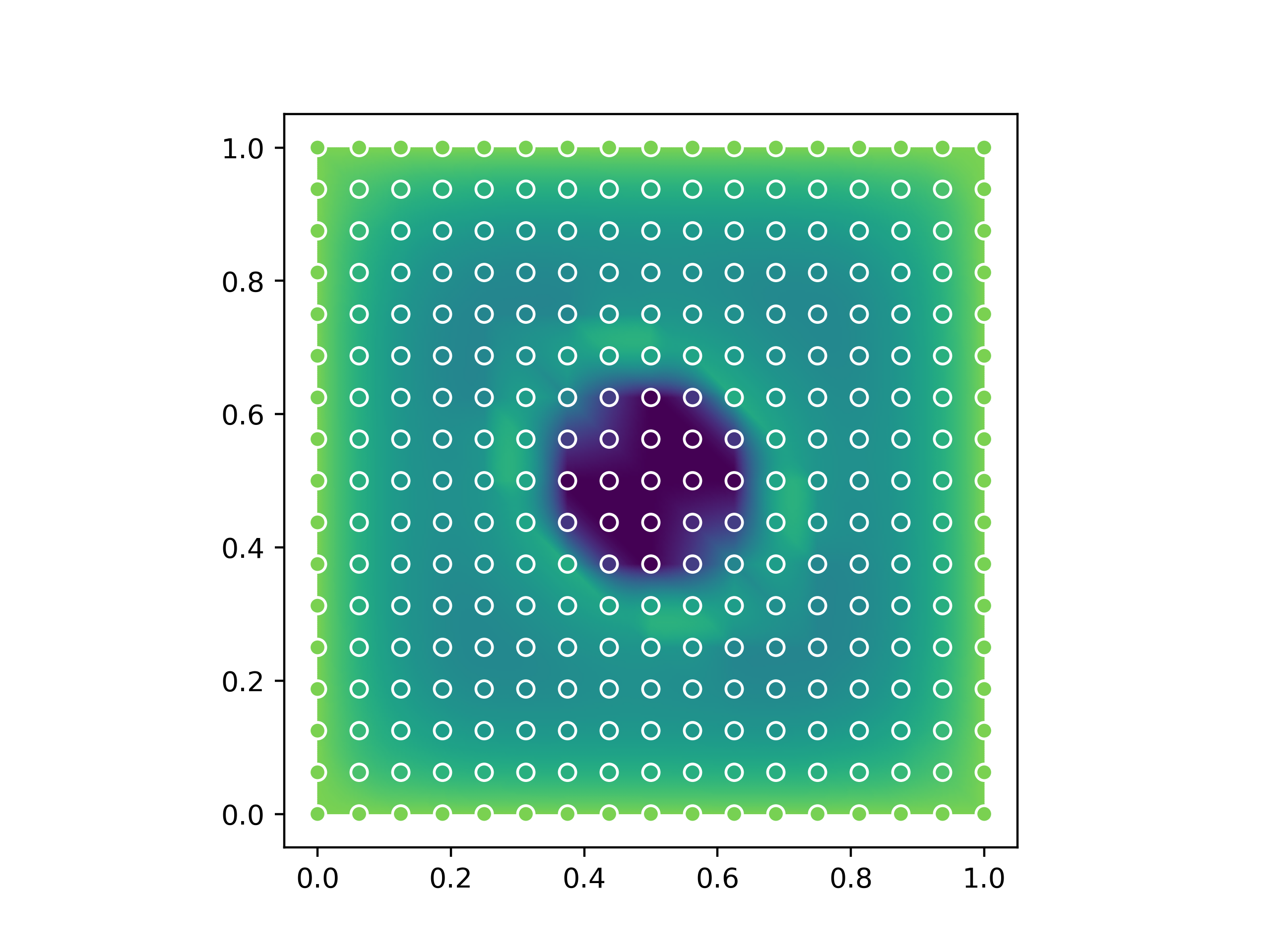}
\end{subfigure}%

\begin{subfigure}{.2\textwidth}
\centering
\rotatebox{90}{\tiny $\mathcal{L}_2$-RIIA($\mathcal{B}_2$)-VI}%
\includegraphics[width=\textwidth] {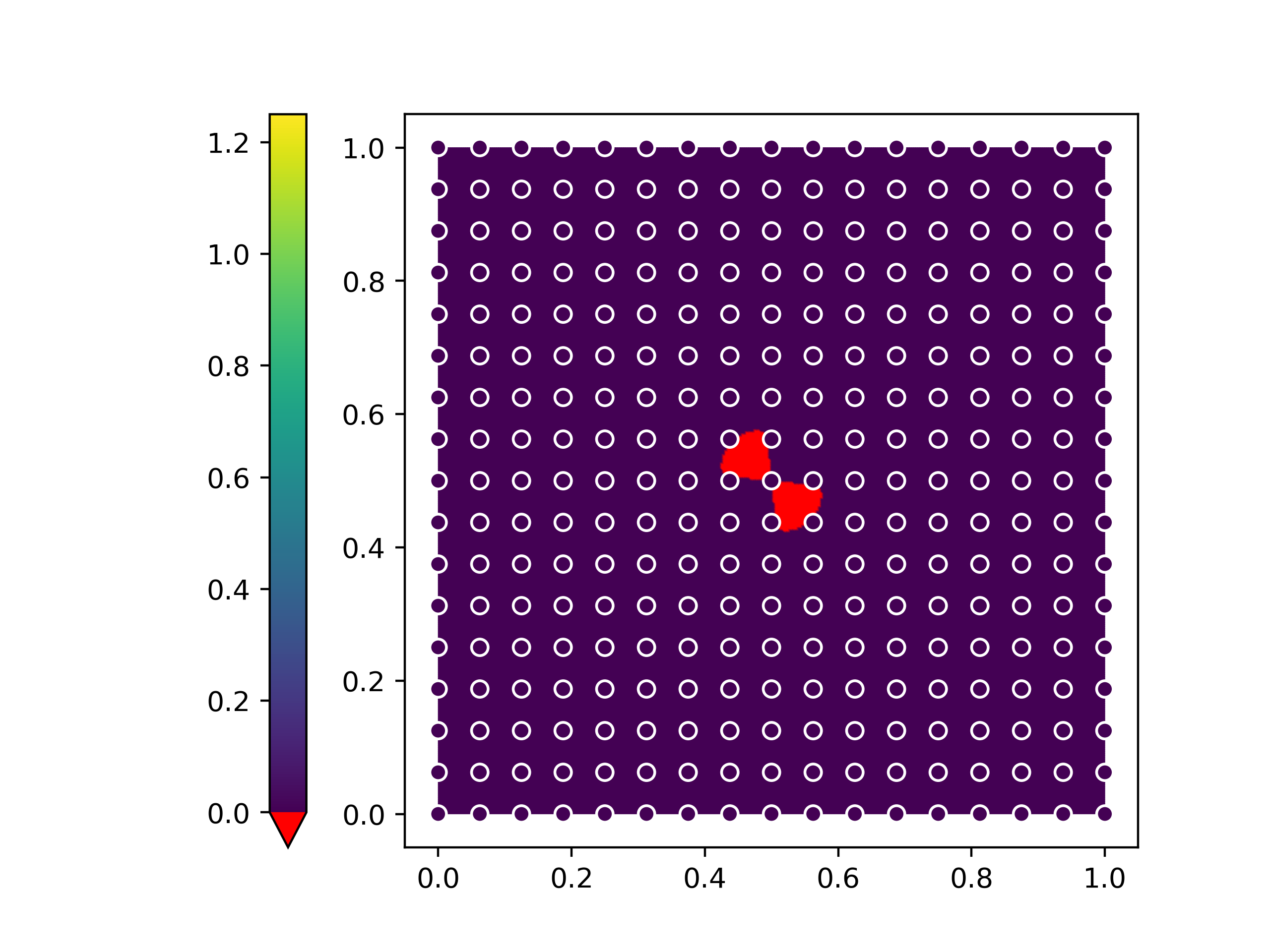}
\end{subfigure}%
\begin{subfigure}{.2\textwidth}
\centering
\includegraphics[width=\textwidth] {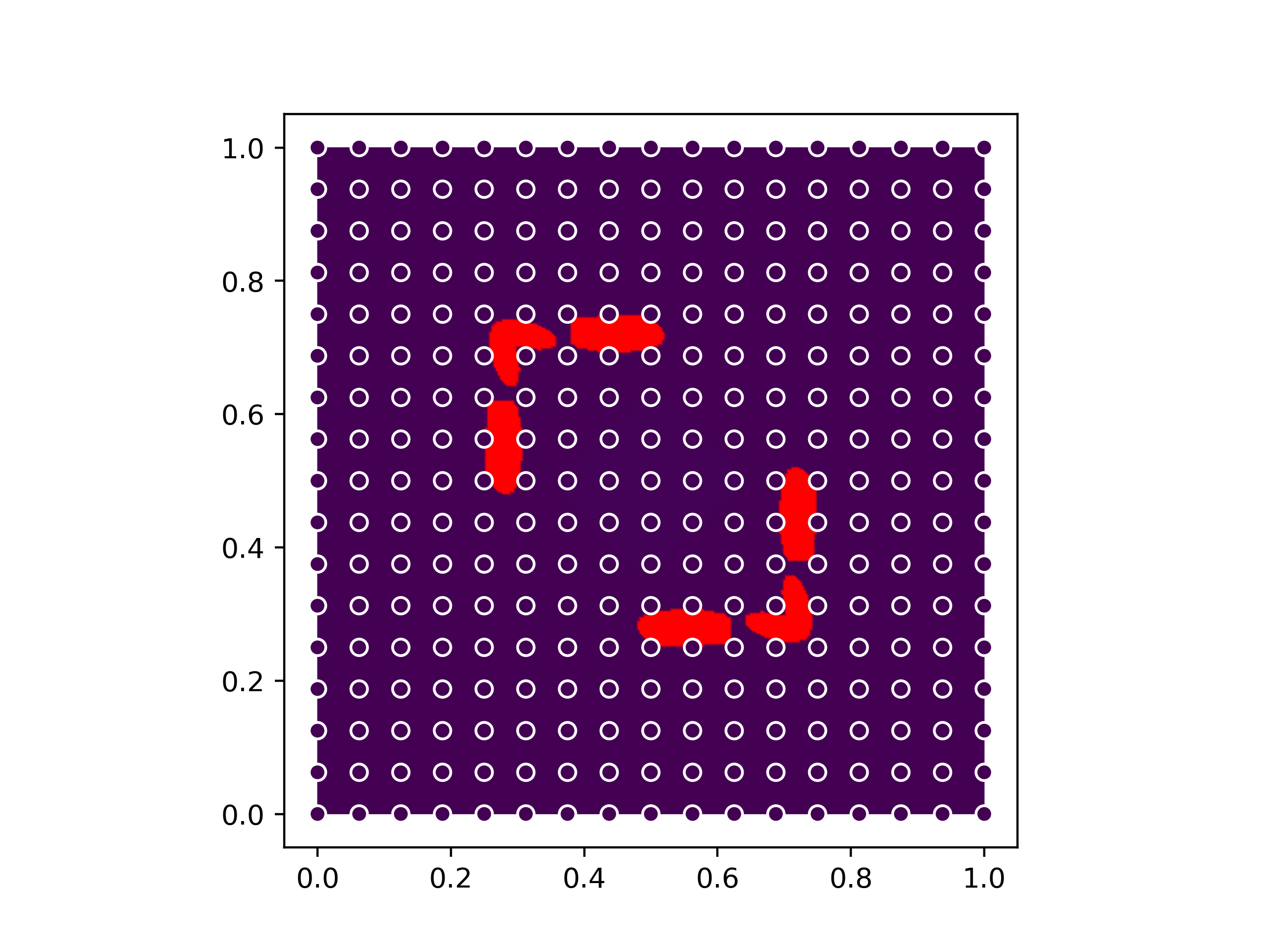}
\end{subfigure}%
\begin{subfigure}{.2\textwidth}
\centering
\includegraphics[width=\textwidth] {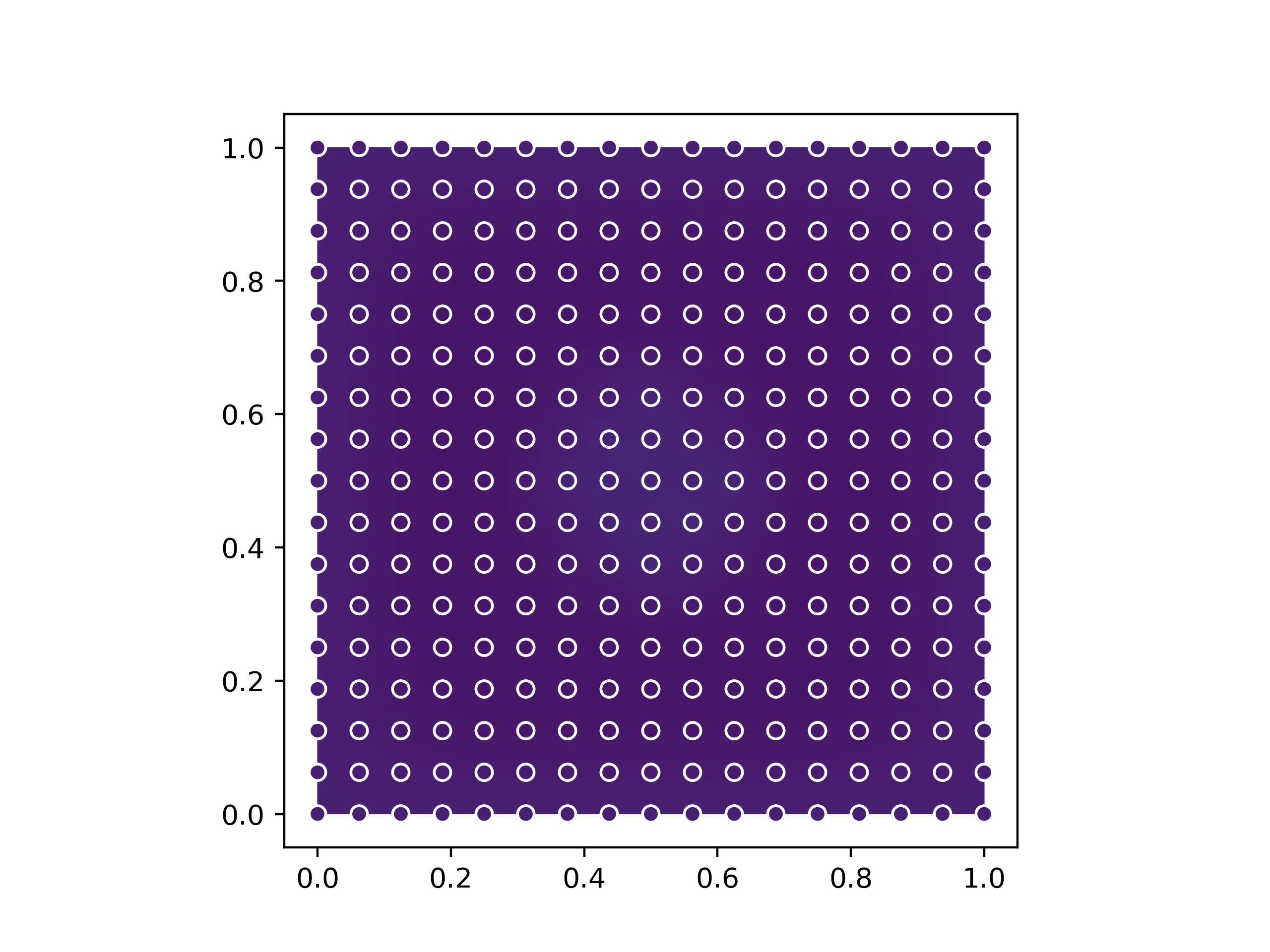}
\end{subfigure}%
\begin{subfigure}{.2\textwidth}
\centering
\includegraphics[width=\textwidth] {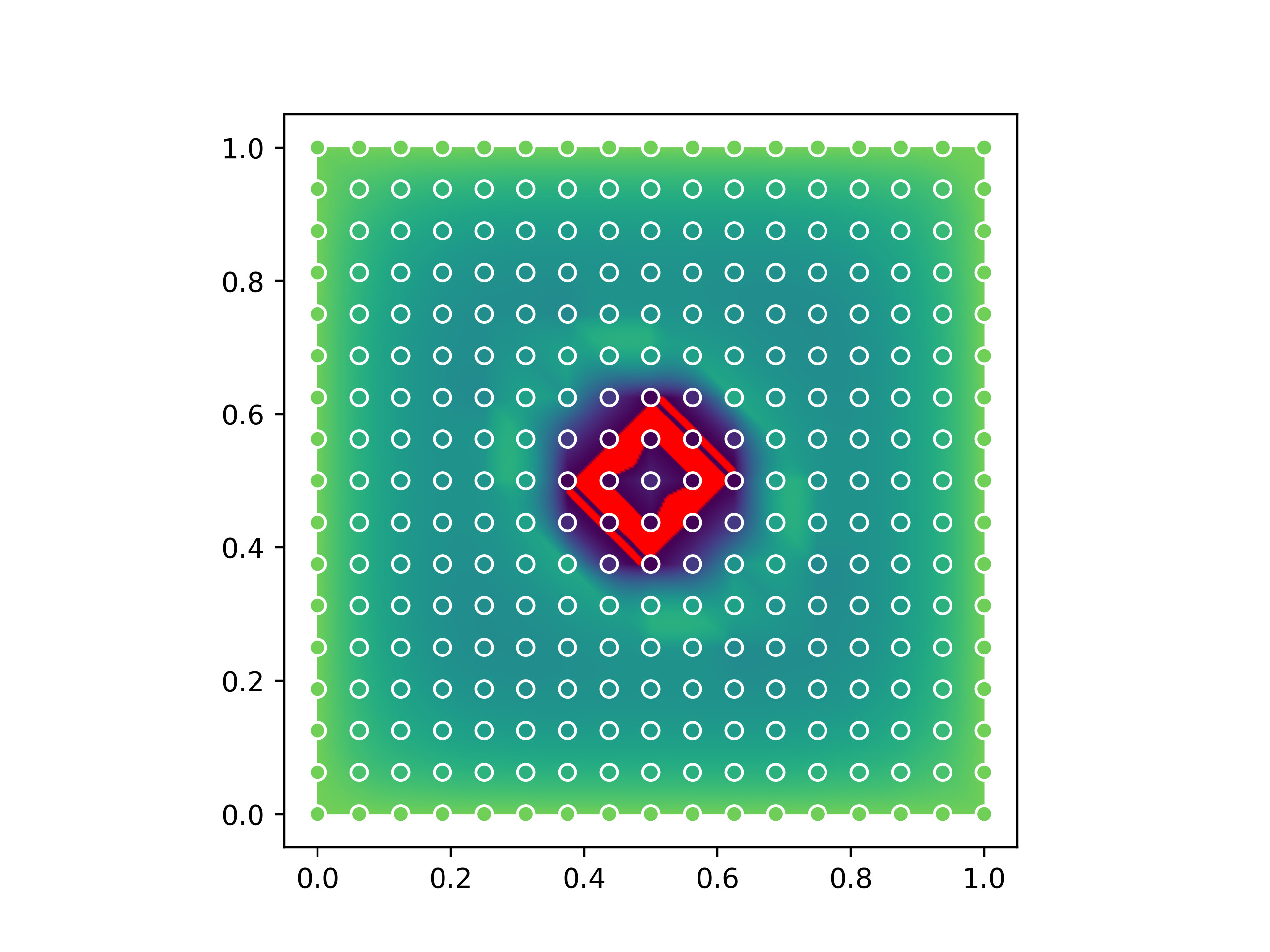}
\end{subfigure}%
\begin{subfigure}{.2\textwidth}
\centering
\includegraphics[width=\textwidth] {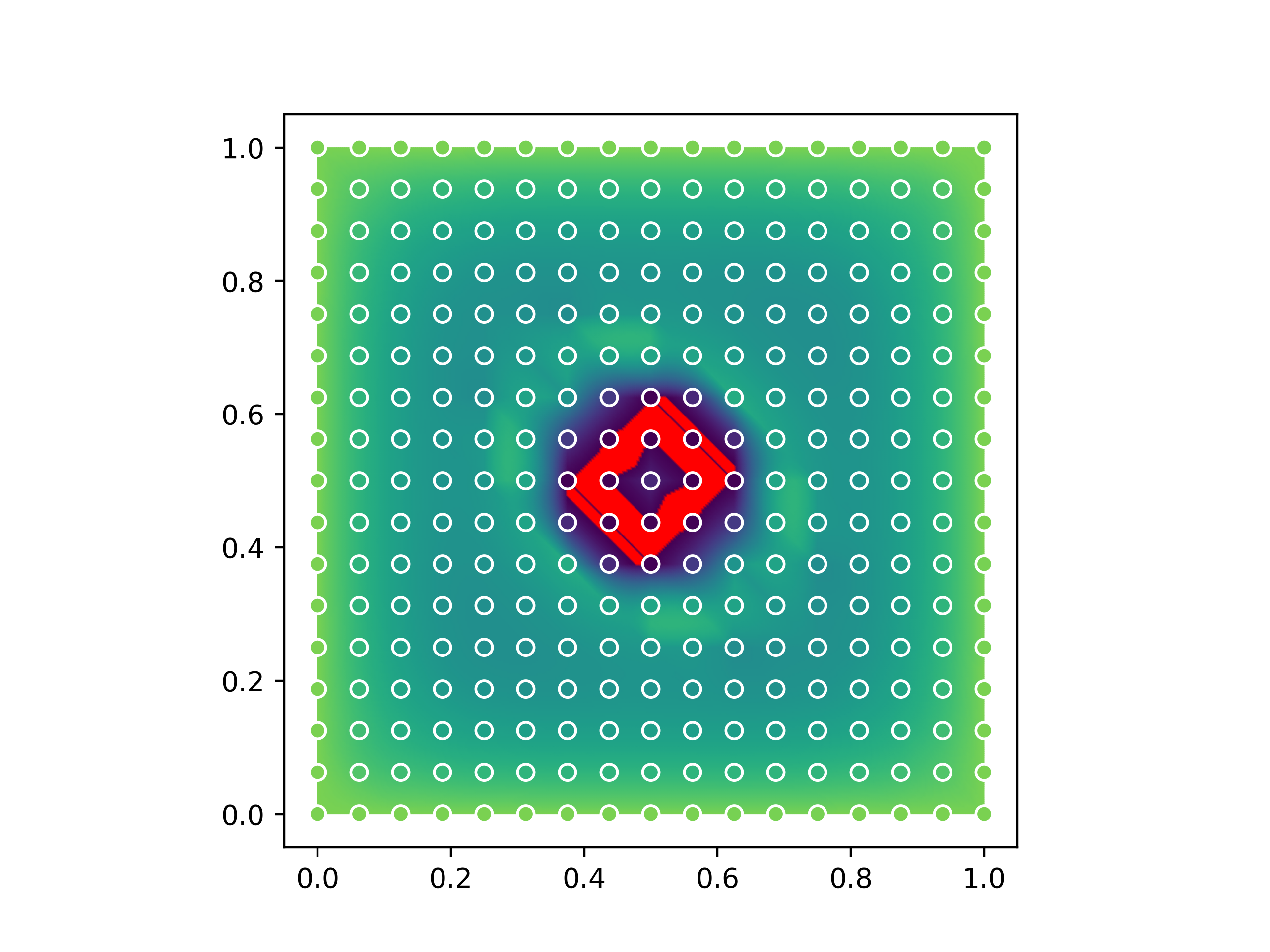}
\end{subfigure}%

\begin{subfigure}{.2\textwidth}
\centering
\rotatebox{90}{\tiny $\mathcal{B}_2$-RIIA($\mathcal{B}_2$)-VI}%
\includegraphics[width=\textwidth] {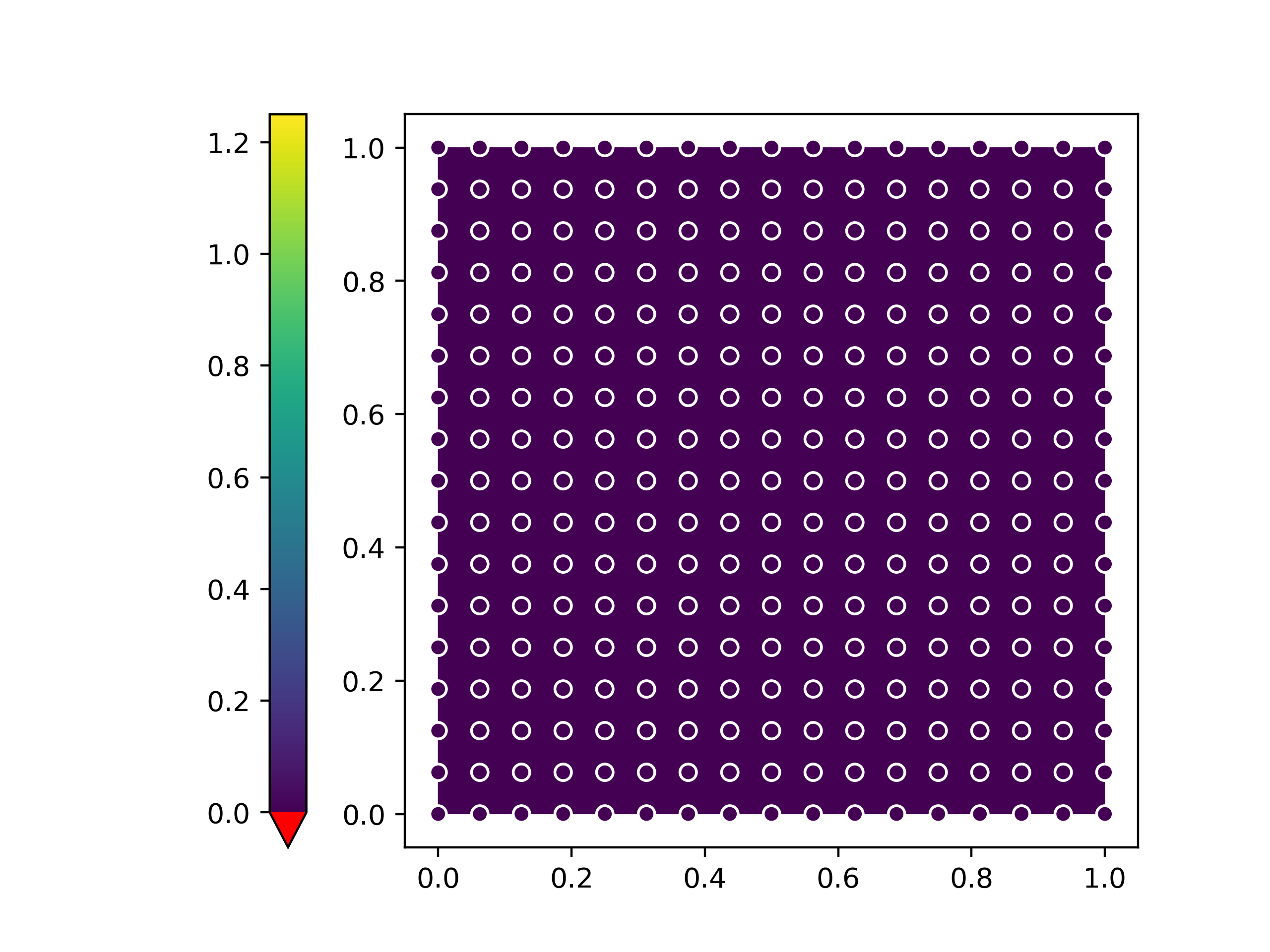}
\subcaption{Initial Time}
\end{subfigure}%
\begin{subfigure}{.2\textwidth}
\centering
\includegraphics[width=\textwidth] {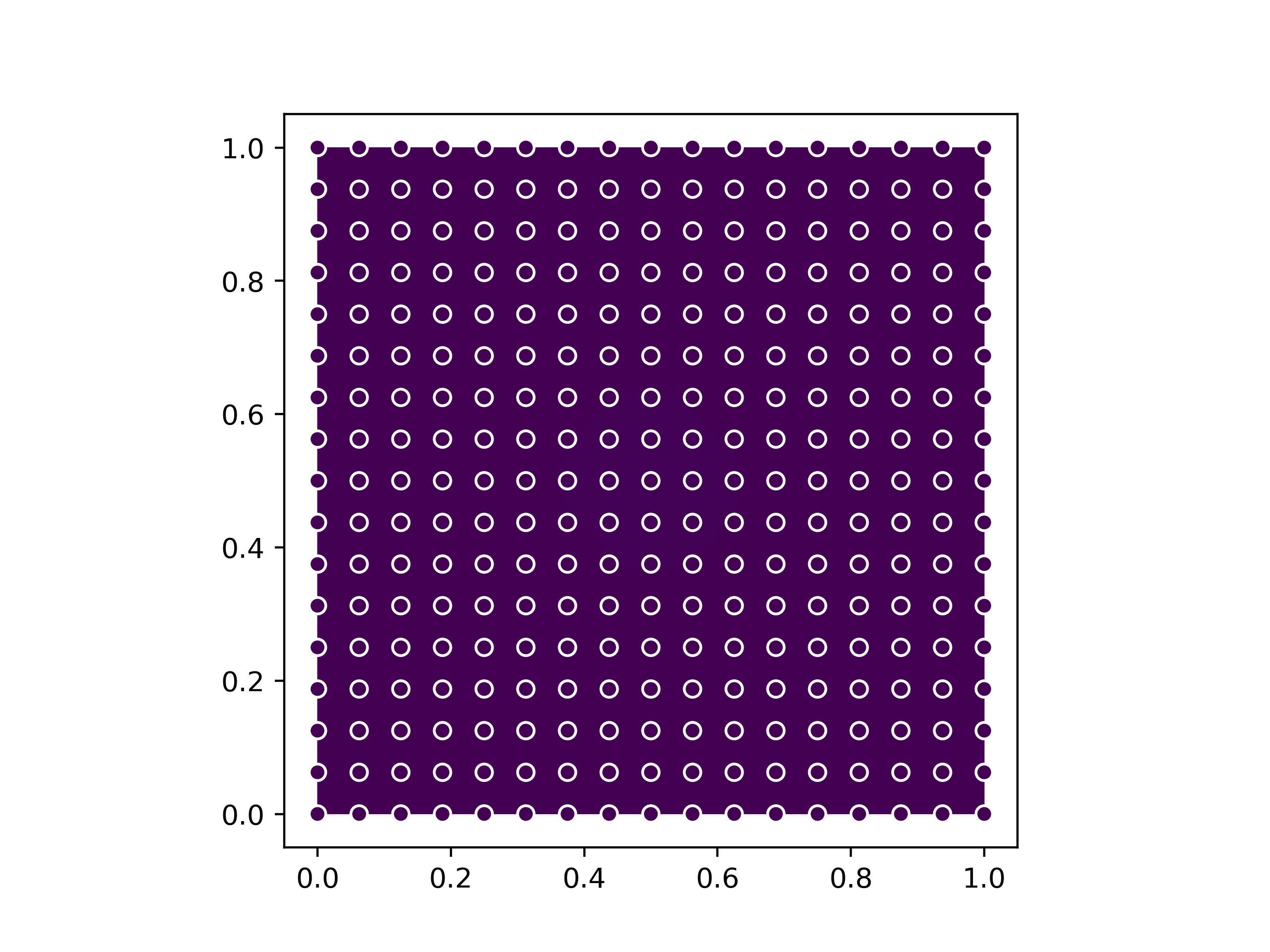}
\subcaption{$t = 0.01k$}
\end{subfigure}%
\begin{subfigure}{.2\textwidth}
\centering
\includegraphics[width=\textwidth] {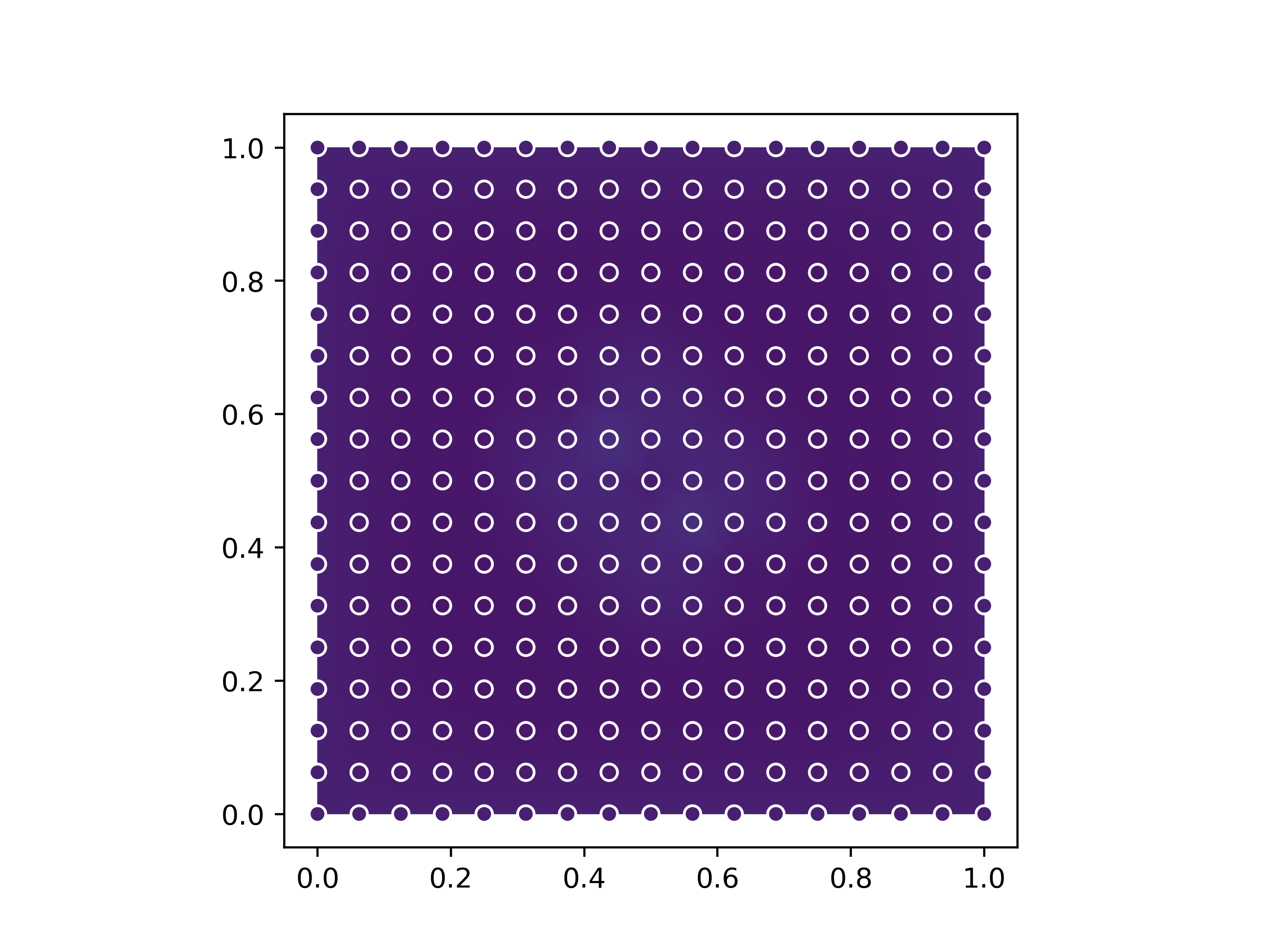}
\subcaption{$t = c_1k = \frac{1}{3}k$}
\end{subfigure}%
\begin{subfigure}{.2\textwidth}
\centering
\includegraphics[width=\textwidth] {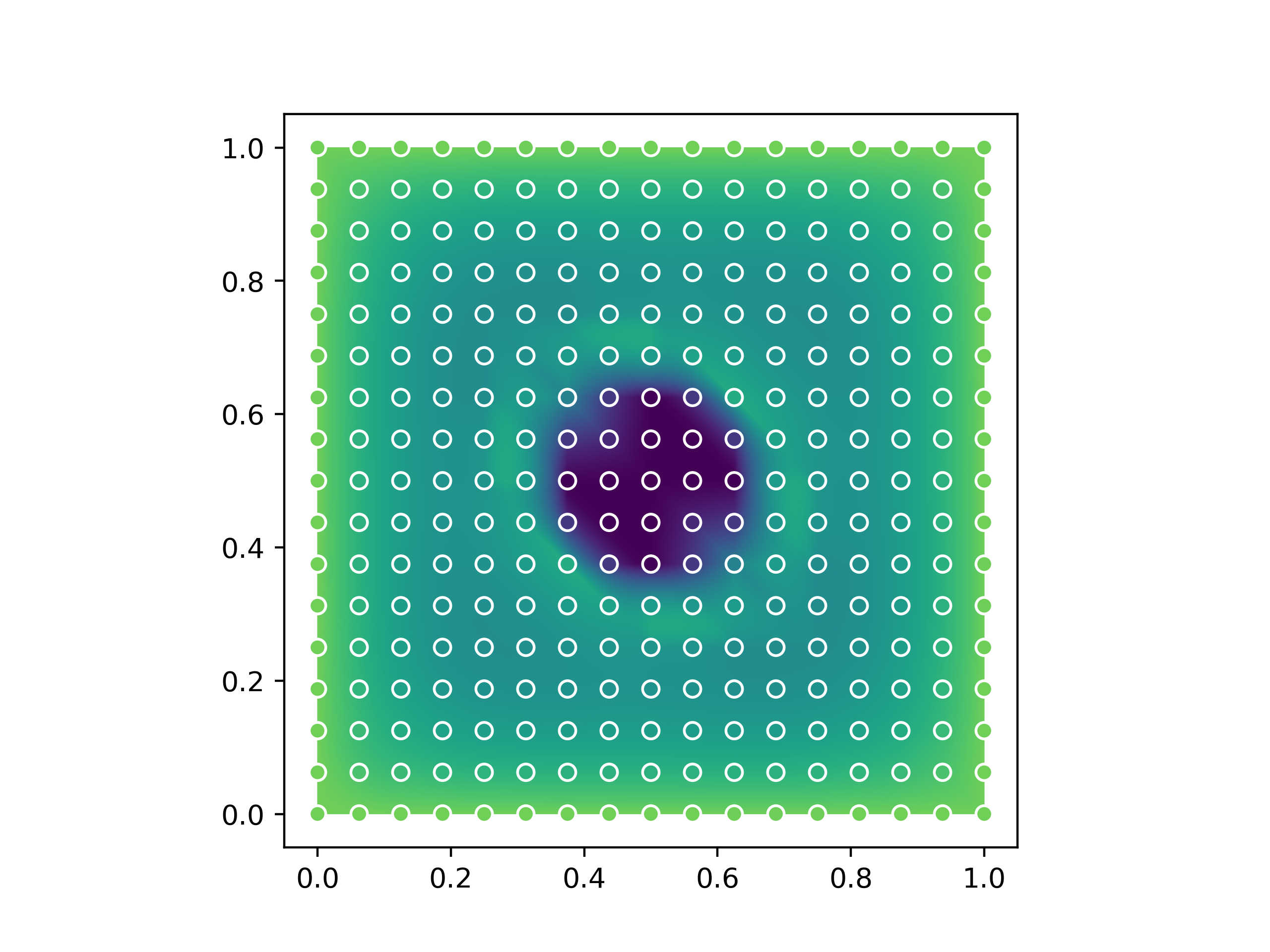}
\subcaption{$t = 0.99k$}
\end{subfigure}%
\begin{subfigure}{.2\textwidth}
\centering
\includegraphics[width=\textwidth] {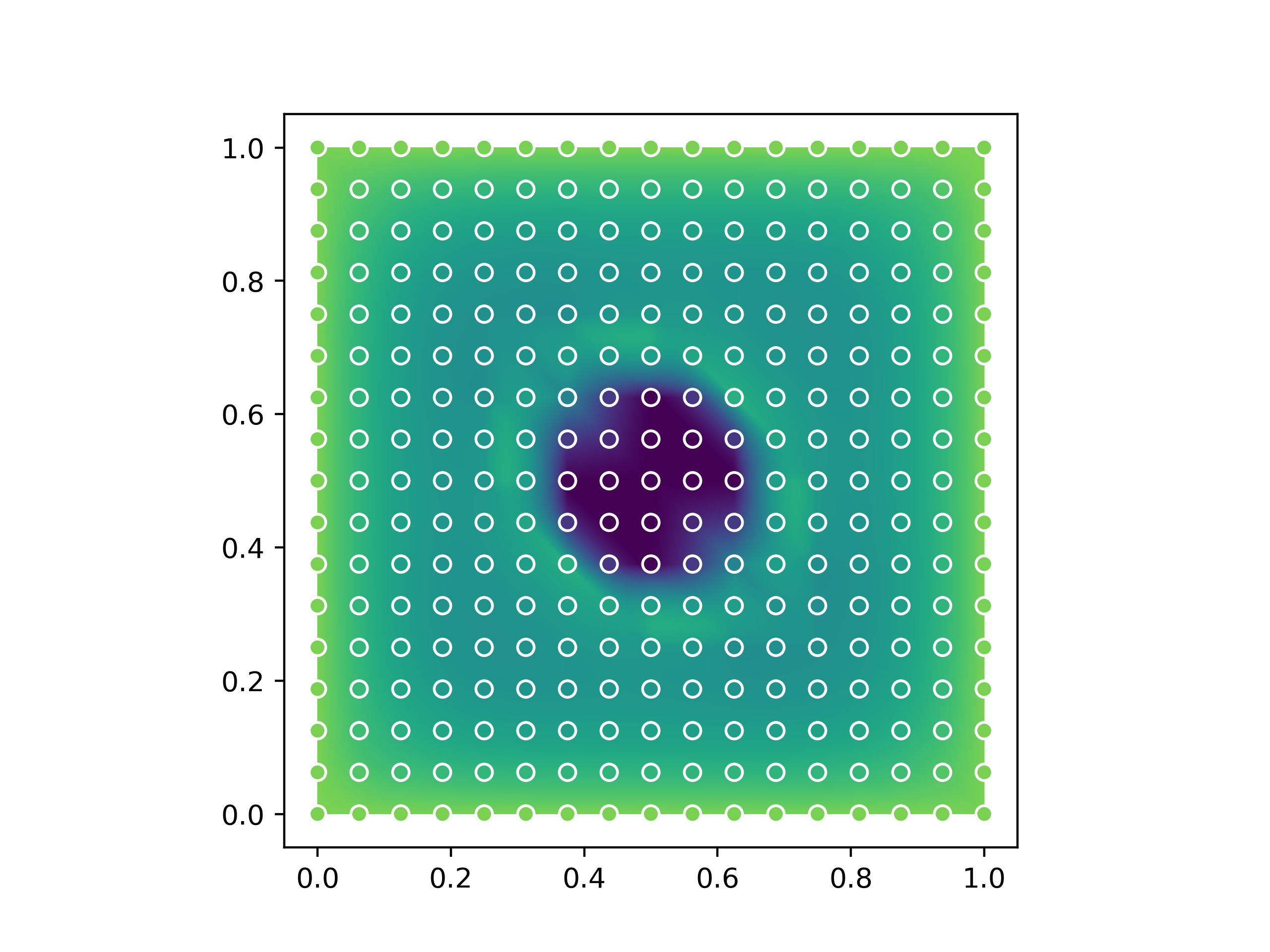}
\subcaption{$t = c_2k = k$}
\end{subfigure}%

\caption{Snapshots of the numerical solution to~\eqref{eq:heat_general} (chosen with exact solution~\eqref{eq:heat_vios_exact}).}
\label{fig:heat_vios_snapshots}
\end{figure}

We partition the unit square into a uniform triangular mesh consisting of $128$ elements. Quadratic finite elements are used in space, and timestepping is done with RadauIIA(2) with step size $k = 1 / 8$. 
The collocation nodes of the RadauIIA method are $\frac{1}{3}$ and $1$. As such, all of the constrained timestepping schemes proposed will guarantee constraint adherence at the spatial degrees of freedom at the times $0, \frac{1}{3}k$, and $k$. If Bernstein finite elements are used, the constraints will be uniform at these times, while the use of Lagrange finite elements will not necessarily result in a uniform-in-space approximation. At off-collocation-node times, constraints at the degrees of freedom are guaranteed only if the Bernstein form of the collocation polynomial is used. In Figure~\ref{fig:heat_vios_snapshots} we plot snapshots of the approximation at the initial time, the collocation nodes, and two intermediate points, $t = 0, 0.01k, \frac{1}{3}k, 0.99k, k$. The circles plotted represent the degrees of freedom of the spatial function space. They are included to highlight the difference between constraint violations at the degrees of freedom, and those only between them. The color red is used to indicate regions of constraint violation.   

With the exception of the constraint violations, qualitatively similar results are seen at the final time using each of the methods shown.  First, we see that using constrained Lagrange polynomials in space and time results in an approximation with spatial violations both on and off of the collocation nodes, and at the spatial nodes between the collocation nodes in time. Choosing the Lagrange Basis in space and Bernstein in time, we see spatial violations at collocation nodes, but the value at the nodes remains positive uniformly in time. When the Bernstein basis is chosen in space, there are no spatial violations at the collocation nodes in time, regardless of the choice of representation for the collocation polynomial. The third and fifth rows of Figure~\ref{fig:heat_vios_snapshots}, then, highlight the difference between constraining the collocation
 polynomial when cast in the Bernstein and Lagrange bases. Lagrange in time allows off-collocation node violations, while Bernstein does not.\newline

\noindent\textbf{Remark:} The differences between each of these methods become vital when the underlying Runge-Kutta method is not stiffly accurate. If the method lacks a collocation node at the right endpoint (Gauss-Legendre, for instance), the bounds constraints will not necessarily be enforced at the discrete times. This defeats the purpose of enforcing bounds constraints at all, and can greatly impact the quality of the numerical approximation.  

\subsection{Cahn-Hilliard equation}\label{sec:cahn_hilliard}

As a final example, we consider the Cahn-Hilliard model of phase separation in a binary fluid~\cite{cahn1958free}. We pose the problem with the Flory-Huggins type logarithmic potential~\eqref{eq:CH_logarithmic_potential} on a domain $\Omega \subset \mathbb{R}^d$ for $d=2,3$. The Cahn-Hilliard equation is the conserved $H^{-1}$ gradient flow the energy functional
\begin{equation}\label{eq:CH_energy_functional}
E(c) = \int_\Omega \left(F(c) + \frac{\epsilon^2}{2}|\nabla c|^2\right)dx,
\end{equation}
i.e.,
\begin{equation}
\begin{split}
\frac{\partial c}{\partial t} &= \nabla \cdot(M(c)\nabla \mu), \\
\mu &= F'(c) - \epsilon^2 \Delta c,
\end{split}
\end{equation}
\cite{chen2019positivity}. In the logarithmic potential, $\theta_0$ and $\theta_c$ are positive constants with $\theta_0<\theta_c$, representing the absolute temperature of the system and the critical temperature of phase separation, respectively.

The logarithmic potential~\eqref{eq:CH_logarithmic_potential} causes phase separation between the two components of the mixture. The states will tend towards the binodal points $c = c^*$ and $c = -c^*$, for some $c^* \in (-1, 1)$. In particular, $c^*$ is the positive root of $\frac{1}{s}\ln\left(\frac{1 + s}{1-s}\right) = \frac{2\theta_c}{\theta_0}$, \cite{copetti1992}. The distance between $\theta_0$ and $\theta_c$ determines the value of $c^*$. When $\theta_c$ is only slightly larger than $\theta_0$, $c^*$ is well inside the interior of $(-1, 1)$, and one can use a typical finite element method. However, the energy functional~\eqref{eq:CH_energy_functional} makes sense only for $c$ satisfying $-1<c<1$. As the distance between $\theta_0$ and $\theta_c$ grows, $c^*$ and $-c^*$ tend toward the singular points of the potential. When this occurs, overshoots and undershoots become a significant problem. Sampling of the potential outside of $(-1, 1)$ can produce nonphysical imaginary component and/or divergence in the nonlinear solver.

Prior work on the logarithmic Cahn-Hilliard equation has focused primarily on finite difference and linear finite element schemes, while some extensions into higher order methods are also available (e.g. \cite{goudenege2012high} analyzes a p-FEM appraoch). Typically, either a problem-specific bounds-preserving spatial discretization is proposed, and then paired with a known time stepping scheme (for which bounds-preservation is then established), or the constants $\theta$ and $\theta_c$ are chosen so that the singularities of the potential are avoided in the computation. In many works the bounds-preservation of the time stepping scheme is established through a time step restriction which is dependent on the separation between $\theta_0$ and $\theta_c$, or through adaptive schemes \cite{copetti1992, liu2023positivity}. Some works have avoided these restrictions through a careful spatial discretization and analysis of existing timestepping techniques, e.g., \cite{chen2019positivity}. Additionally, some works choose to regularize the logarithm near the singularities, and instead solve this regularized problem. Here, we apply the bounds-preserving method $\mathcal{B}_2$-RIIA($\mathcal{B}_2$)-VI to approximate the solution to~\eqref{eq:cahn_hilliard_general} with the temperatures $\theta$ and $\theta_c$ sufficiently separated so as to cause issues numerically if bounds are not enforced. We obtain a solution which avoids the singularities of the logarithm independent of the mesh diameter and temporal step size, and allows for the high-order approximation of the solution in both space and time.
We also avoid regularizing the problem, thus arriving at a solution to the variational inequality, not a regularized version thereof.

Towards a numerical approximation, we take $\mathcal{V} = H^1(\Omega) \times H^1(\Omega)$,
and we write the system for the order parameter $c:\Omega \rightarrow [-1,1]$ and chemical potential $\mu : \Omega \rightarrow \mathbb{R}$ as
\begin{equation}\label{eq:cahn_hilliard_general}
  \begin{split}
    \left(\frac{\partial c}{\partial t}, v \right) &= - M\left(\nabla \mu, \nabla v\right), \\
    \left(\mu, w \right) &= \left( F'(c), w\right) + \epsilon^2\left( \nabla c, \nabla w \right) \\
  \end{split}
\end{equation}
for all $v, w \in H^1(\Omega)$ and all $0 < t \leq T$. The system is closed with boundary conditions, which we take as either homogeneous Neumann conditions on both $c$ and $\mu$, or periodic boundary conditions.

In the following examples, we take $\Omega = [0, 1]^2$. We use an $100\times 100$ square mesh, with each square divided into two triangles. We represent both the concentration and the chemical potential as piecewise quadratic functions on the mesh, represented in the Bernstein basis. We advance the system in time using the 2-stage RadauIIA collocation method, with the collocation polynomial represented in the Bernstein basis.

To avoid sampling the logarithms outside of $(-1, 1)$, bounds must be enforced on the numerical solution. To do this, we choose some $\delta_{\text{b}}$ such that $c^* \in (-1 + \delta_{\text{b}}, 1 - \delta_{\text{b}})$, and then require that $c\in [-1 + \delta_{\text{b}}, 1 - \delta_{\text{b}}]$ uniformly in space and time.

Finally, we must regularize the linear solver which is used to advance the nonlinear solve at each timestep. While the variational inequality will ensure that the solution will satisfy the bounds, the linear solves at each step have no such guarantee. As in \cite{chen2019positivity}, we choose some $\delta_{\text{reg}} \in (0, 0.25)$ and define 
\begin{equation*}
\ln_\text{reg}(s) = \begin{cases} \ln(s) &\text{ if } s > \delta_{\text{reg}}\\ \ln(\delta) + \frac{s - \delta_\text{reg}}{\delta_{\text{reg}}} &\text{ if } s \leq \delta_{\text{reg}} \end{cases}
\end{equation*}
We then use the regularized logarithm to compute the Jacobian used to solve the nonlinear system. While this is likely to slow the convergence of the sytem, it allows for the solution to be found without overflow in the linear solves at each step.

For the following examples, we use the parameters $\theta_0 = 2.0$, $\theta_c = 3.5$, $\epsilon = 0.01$, $k =0.0001$, $M = 1$, $\delta_{reg} = 0.001$ and $\delta_b = 10^{-8}$.  In particular, $\theta_0$ and $\theta_c$ are chosen such that, if bounds are not enforced on the numerical solution, the nonlinear solve fails.

First, we approximate the solution at time $t = 1.0$ subject to the initial condition
\begin{equation}\label{eq:CH_ex6_init}
c(0) =  \frac{1}{4} \sin^2 (2\pi x) \sin^2 (2\pi y) \sin(12\pi x)\sin(12\pi y).
\end{equation}

\begin{figure}[ht]
\centering
\begin{subfigure}{0.25\textwidth}
\centering
{\includegraphics[trim={10cm 3cm 19cm 2cm},clip, width=\textwidth] {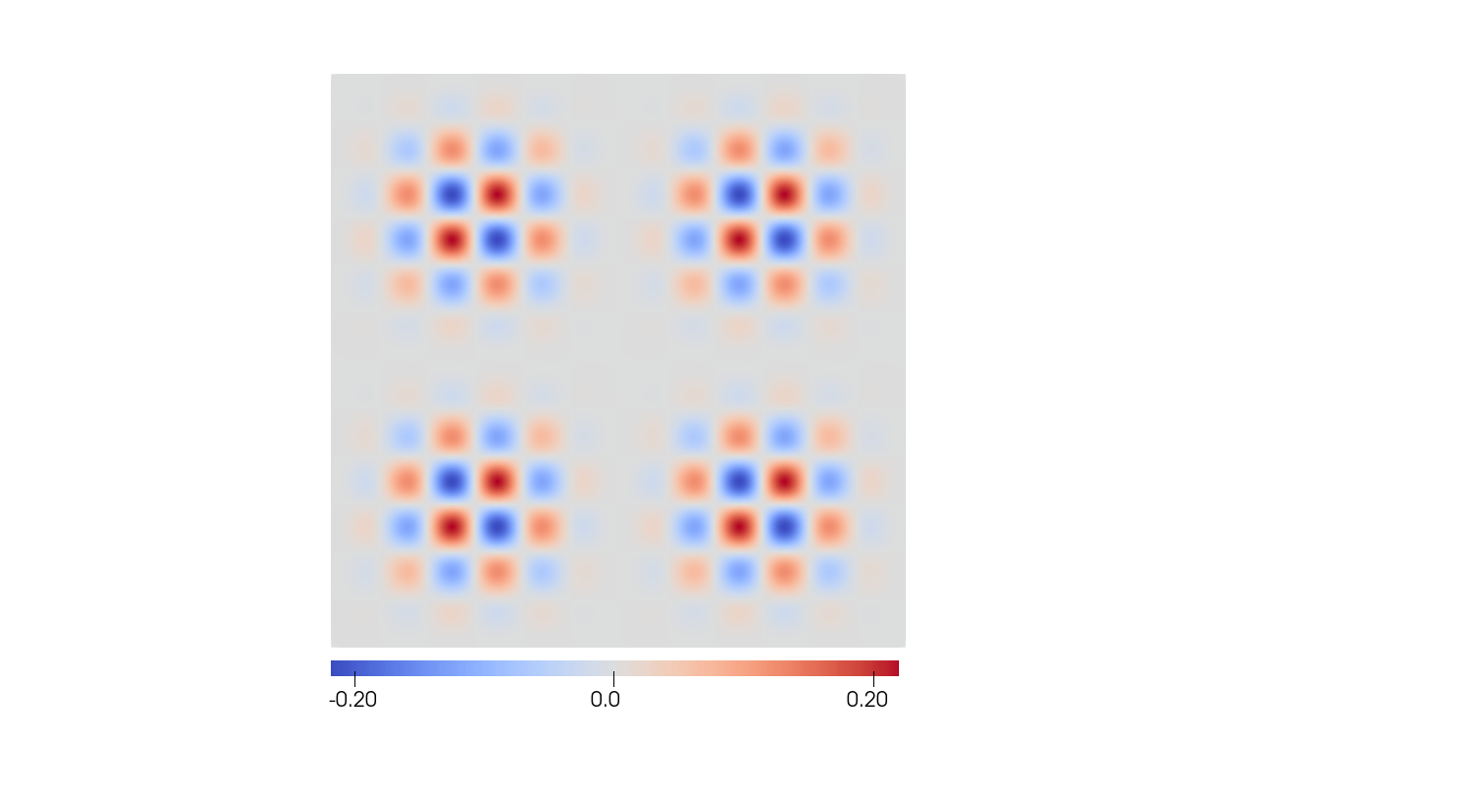}}
\subcaption{$t = 0.0$}
\end{subfigure}%
\begin{subfigure}{0.25\textwidth}
\centering
{\includegraphics[trim={10cm 3cm 19cm 2cm},clip, width=\textwidth] {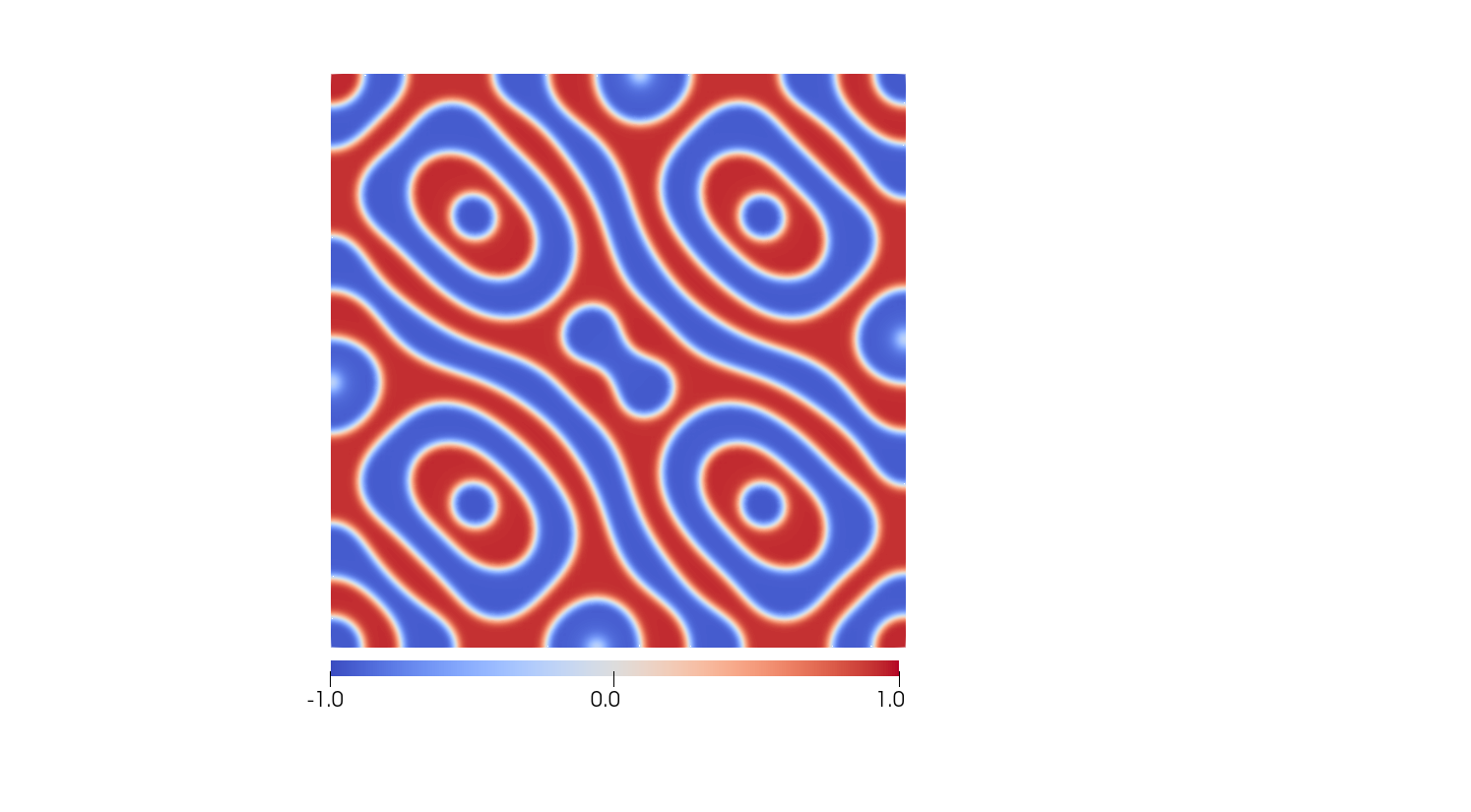}}
\subcaption{$t = 0.01$}
\end{subfigure}%
\begin{subfigure}{0.25\textwidth}
\centering
{\includegraphics[trim={10cm 3cm 19cm 2cm},clip, width=\textwidth] {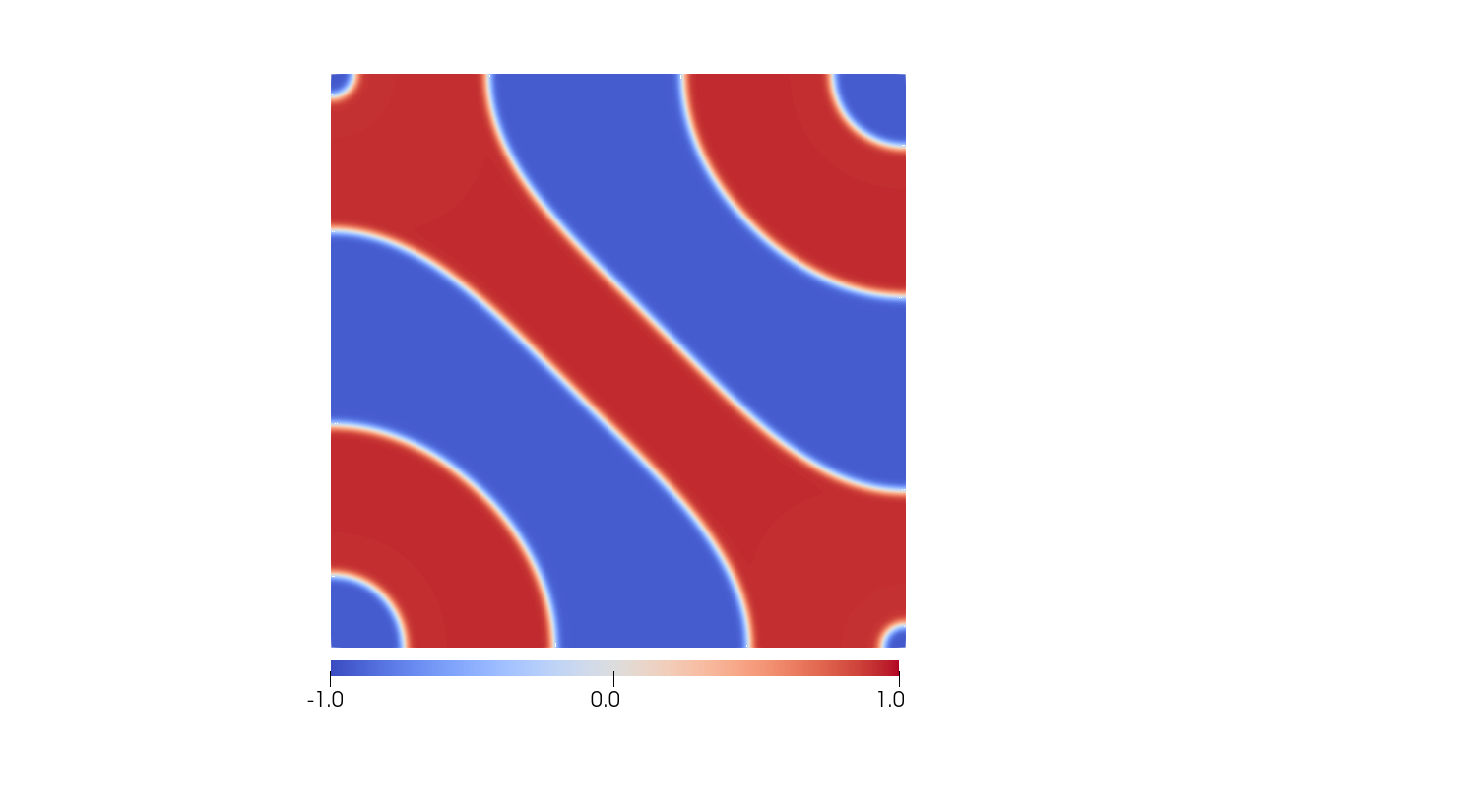}}
\subcaption{$t = 0.5$}
\end{subfigure}%
\begin{subfigure}{0.25\textwidth}
\centering
{\includegraphics[trim={10cm 3cm 19cm 2cm},clip,width=\textwidth] {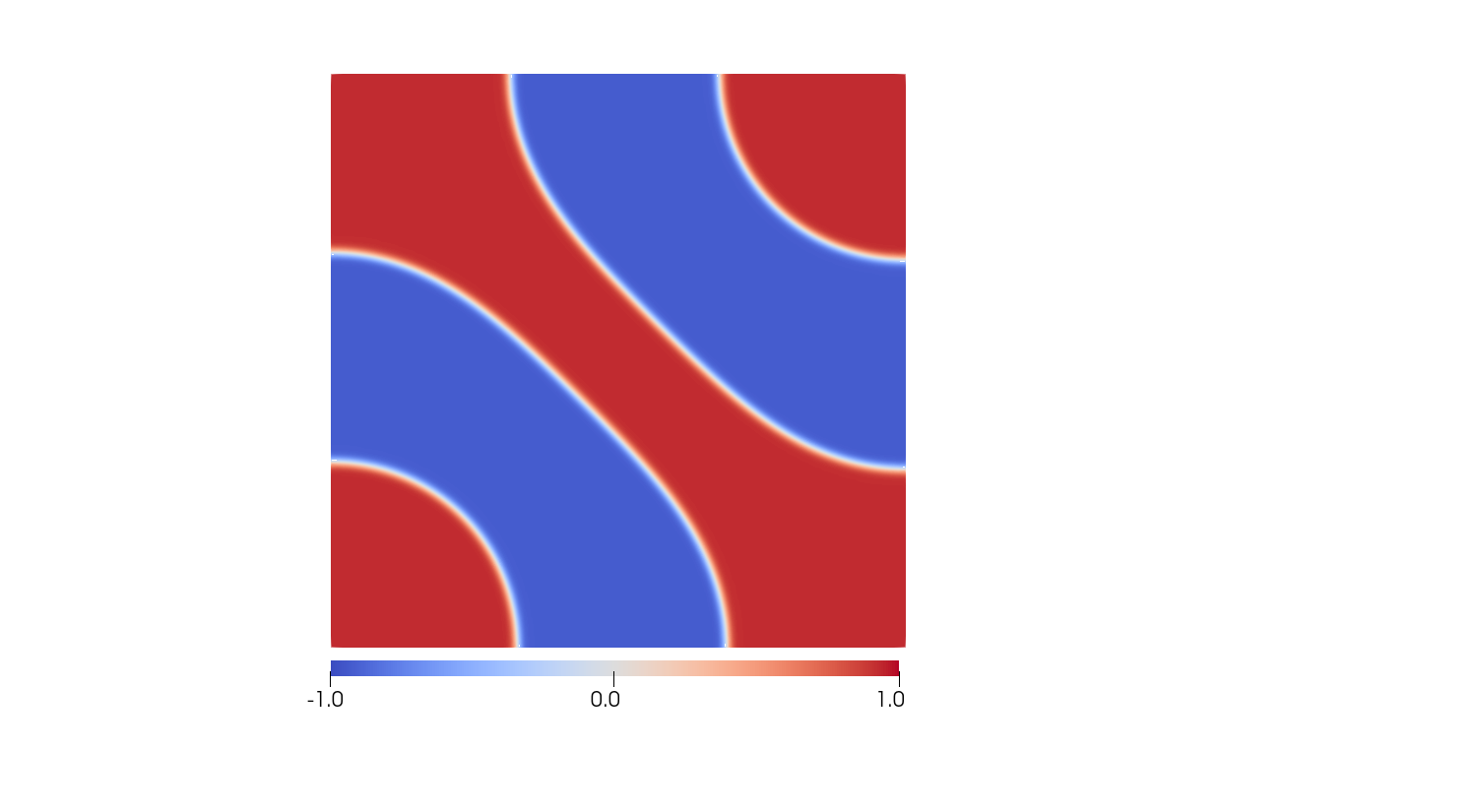}}
\subcaption{$t = 1.0$}
\end{subfigure}%
\caption{Snapshots of the numerical solution to \eqref{eq:cahn_hilliard_general} subject to homogeneous Neumann boundary conditions and the initial concentration~\eqref{eq:CH_ex6_init} at $t =0, 0.01, 0.5, 1.0$.}
\label{fig:CH_ex6_neumann_snaps}
\end{figure}

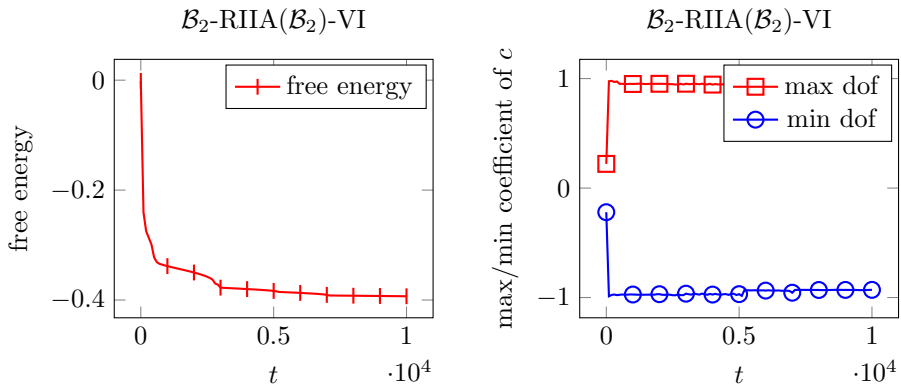
\begin{figure}[ht]

\begin{subfigure}{0.48\textwidth}
\centering
\begin{tikzpicture}
\centering
\begin{axis}[width=\textwidth, legend style={nodes={scale=1.0, transform shape}}, 
        title={$\mathcal{B}_2$-RIIA($\mathcal{B}_2$)-VI},
	ylabel near ticks,
	xlabel=$t$,
    	ylabel={free energy}]
    
    \addplot[red,mark=|, mark size=3pt, mark repeat=10,  each nth point=100, filter discard warning=false, thick] table[y=energy, col sep=comma] {CH6N_data.csv};
    \addlegendentry {free energy};  
\end{axis}
\end{tikzpicture}
\end{subfigure}\hspace{0.04\textwidth}%
\begin{subfigure}{0.48\textwidth}
\centering
\begin{tikzpicture}
\centering
\begin{axis}[width=\textwidth, legend style={nodes={scale=1.0, transform shape}}, 
        title={$\mathcal{B}_2$-RIIA($\mathcal{B}_2$)-VI},
	ylabel near ticks,
	xlabel=$t$,
    	ylabel={max/min coefficient of $c$}]
    
    \addplot[red,mark=square, mark size=3pt, mark repeat=10,  each nth point=100, filter discard warning=false,  thick] table[y=max_dof, col sep=comma] {CH6N_data.csv};
    \addlegendentry {max dof};
    
    \addplot[blue,mark=o, mark size=3pt, mark repeat=10,  each nth point=100, filter discard warning=false,  thick] table[y=min_dof, col sep=comma] {CH6N_data.csv};
    \addlegendentry {min dof};
\end{axis}

\end{tikzpicture}
\end{subfigure}%
\caption{The free energy (left) and maximum and minimum coefficients (right) of the numerical solution to~\eqref{eq:cahn_hilliard_general} subject to homogeneous Neumann boundary conditions and the initial concentration~\eqref{eq:CH_ex6_init}.}
\label{fig:CH_ex6_neumann_data}
\end{figure}

Snapshots of the evolution, subject to homogeneous Neumann boundary conditions, are shown in Figure~\ref{fig:CH_ex6_neumann_snaps}. The total free energy of the system is often used to monitor the quality of a numerical approximation to the Cahn-Hilliard system. 
If all is well, the total free energy should be monotonically decreasing. We let the simulation run to final time $t = 1.0$ and compute the free energy at each timestep. The free energy and maximum/minimum coefficients of the solution are shown in Figure~\ref{fig:CH_ex6_neumann_data}.

\begin{figure}[ht]
\centering
\begin{subfigure}{0.25\textwidth}
\centering
{\includegraphics[trim={10cm 3cm 19cm 2cm},clip, width=\textwidth] {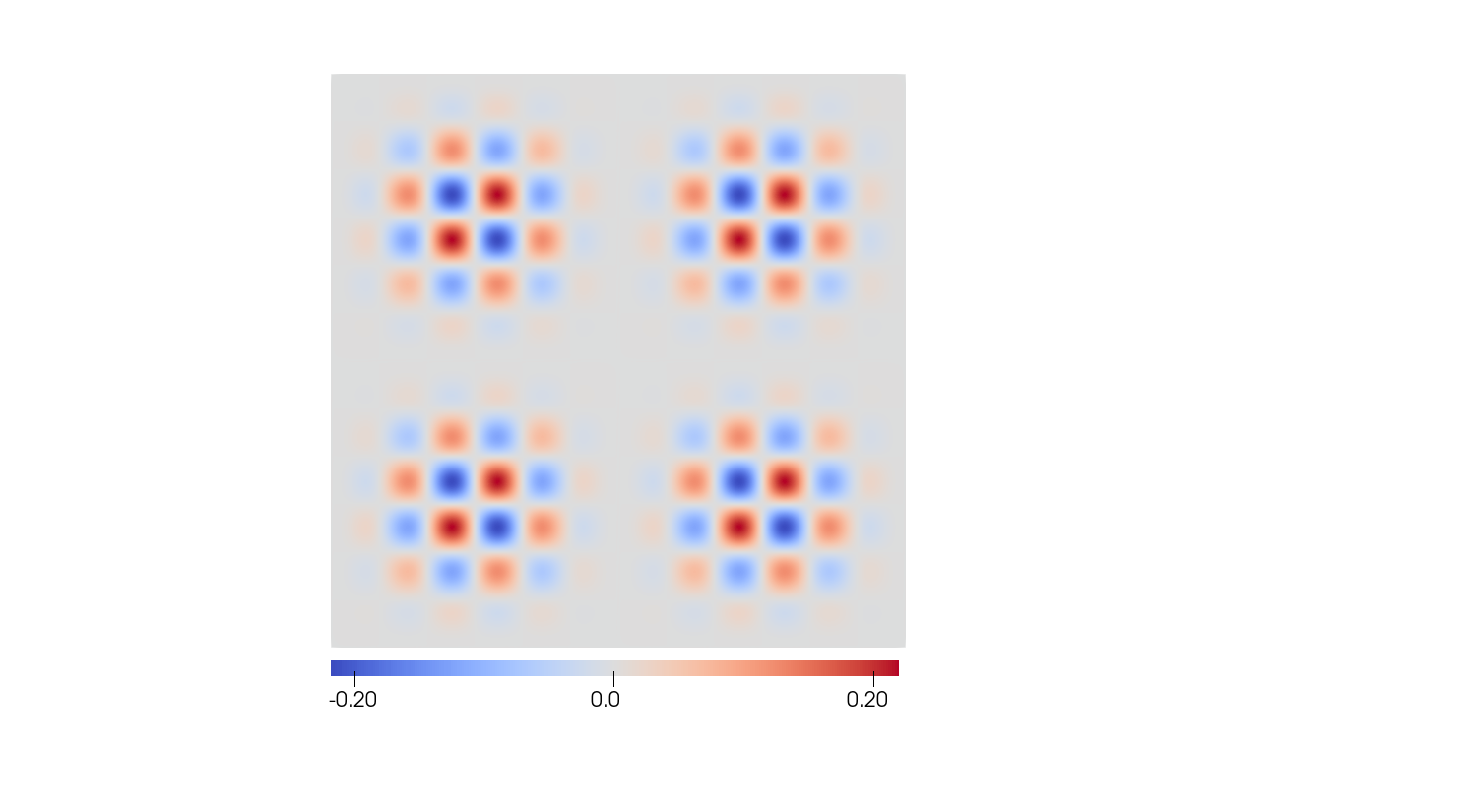}}
\subcaption{$t = 0.0$}
\end{subfigure}%
\begin{subfigure}{0.25\textwidth}
\centering
{\includegraphics[trim={10cm 3cm 19cm 2cm},clip,  width=\textwidth] {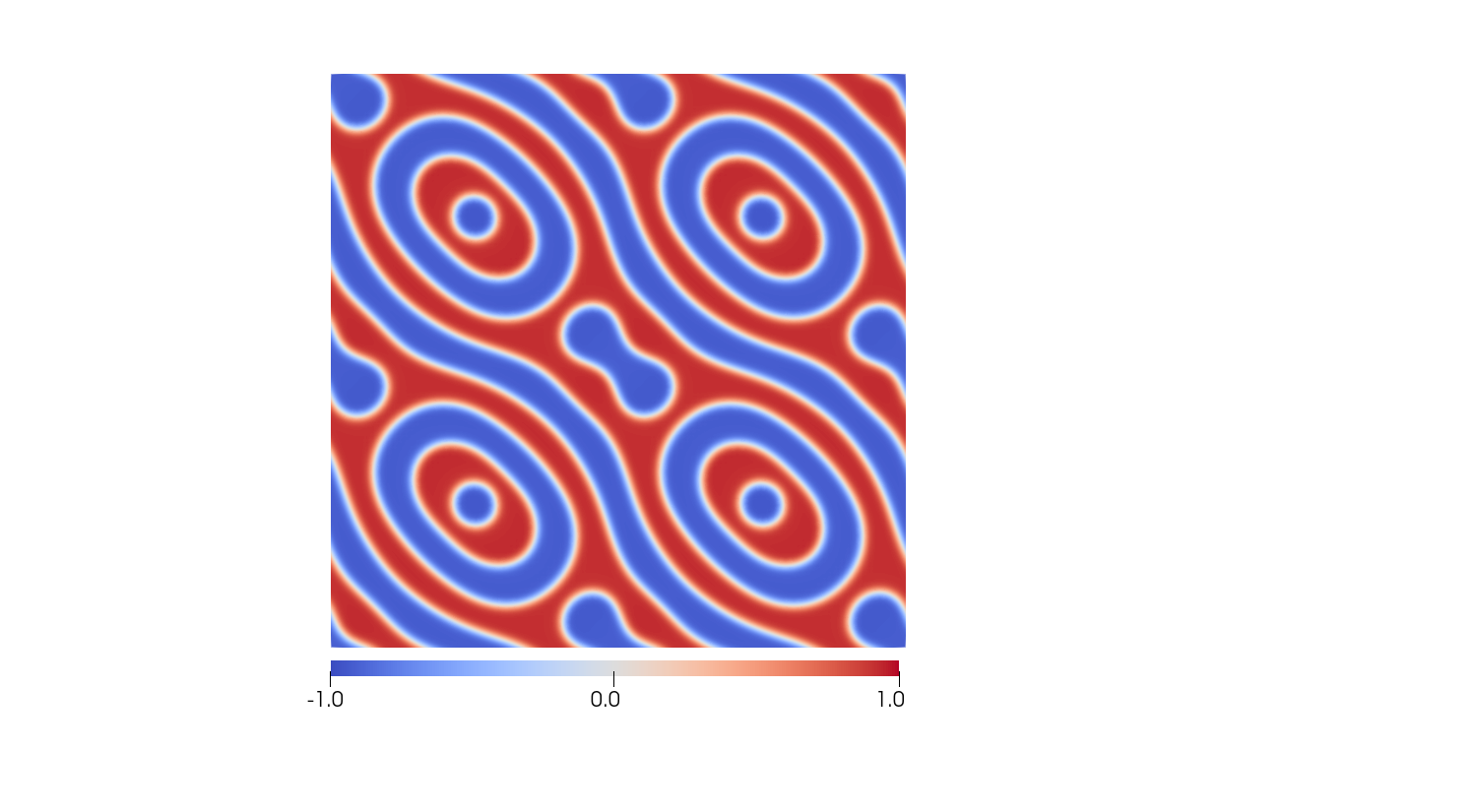}}
\subcaption{$t = 0.01$}
\end{subfigure}%
\begin{subfigure}{0.25\textwidth}
\centering
{\includegraphics[trim={10cm 3cm 19cm 2cm},clip, width=\textwidth] {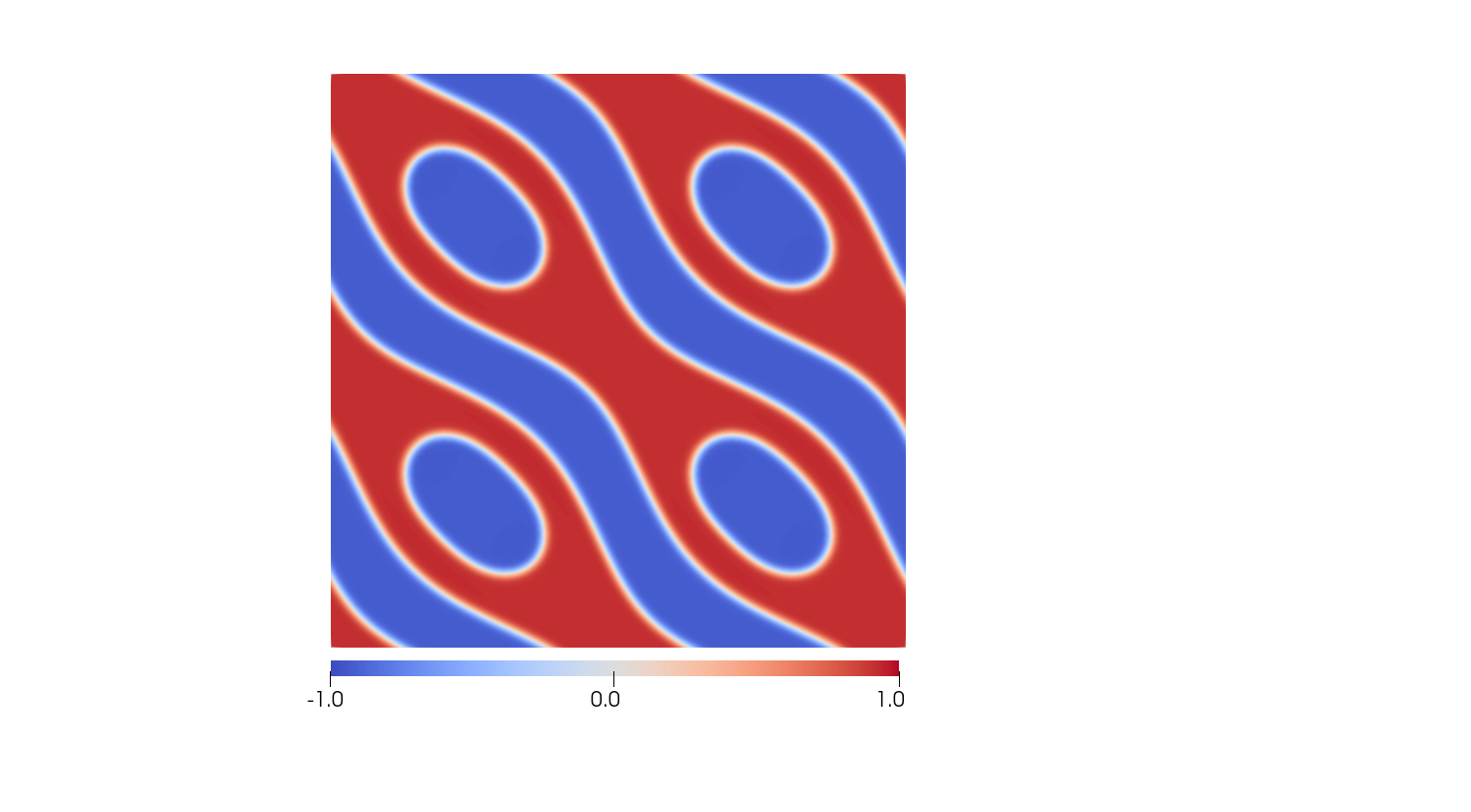}}
\subcaption{$t = 0.06$}
\end{subfigure}%
\begin{subfigure}{0.25\textwidth}
\centering
{\includegraphics[trim={10cm 3cm 19cm 2cm},clip,  width=\textwidth] {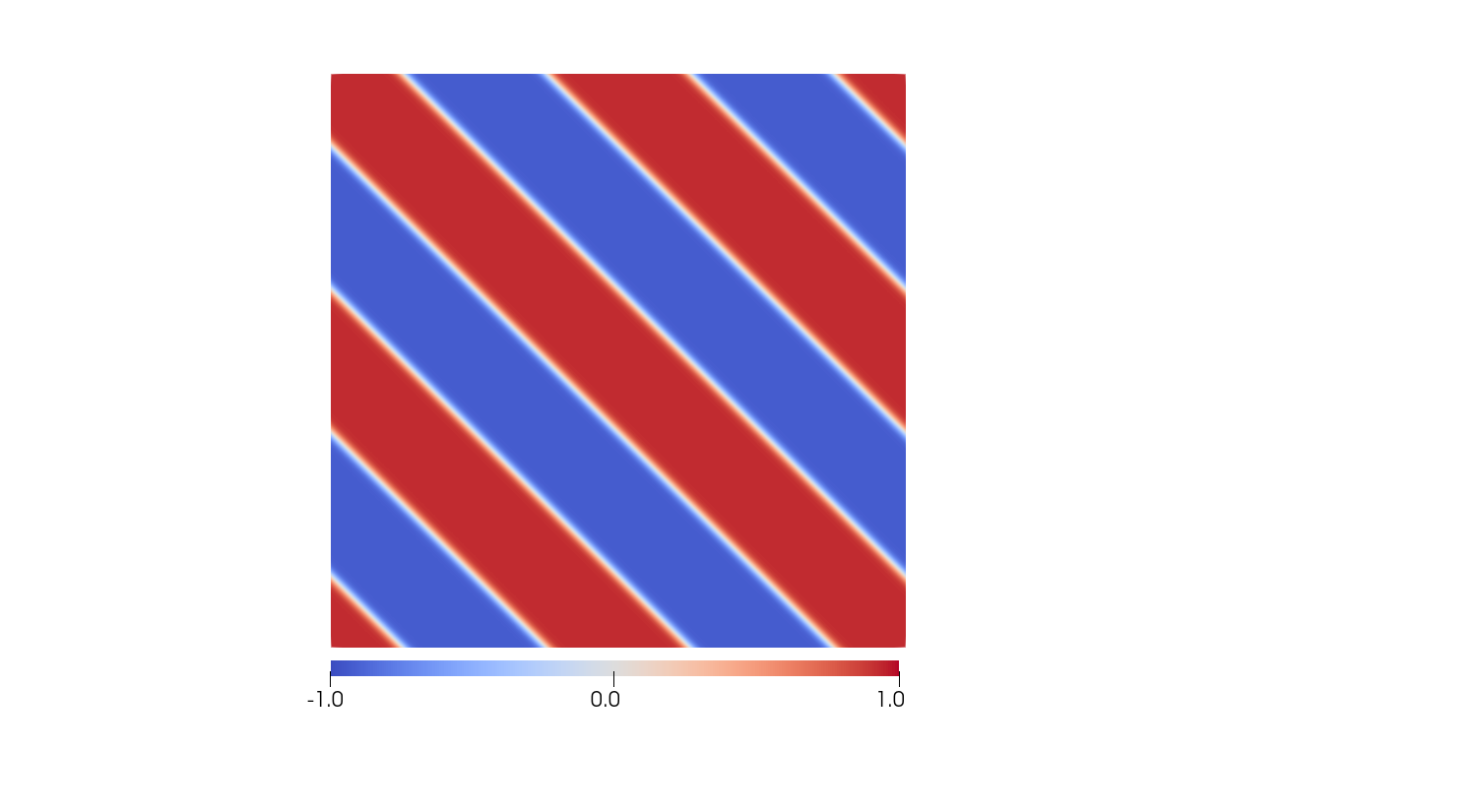}}
\subcaption{$t = 1.0$}
\end{subfigure}%
\caption{Snapshots of the numerical solution to \eqref{eq:cahn_hilliard_general} subject to periodic boundary conditions and the initial concentration~\eqref{eq:CH_ex6_init} at $t =0, 0.01, 0.06, 1.0$.}
\label{fig:CH_ex6_per_snaps}
\end{figure}

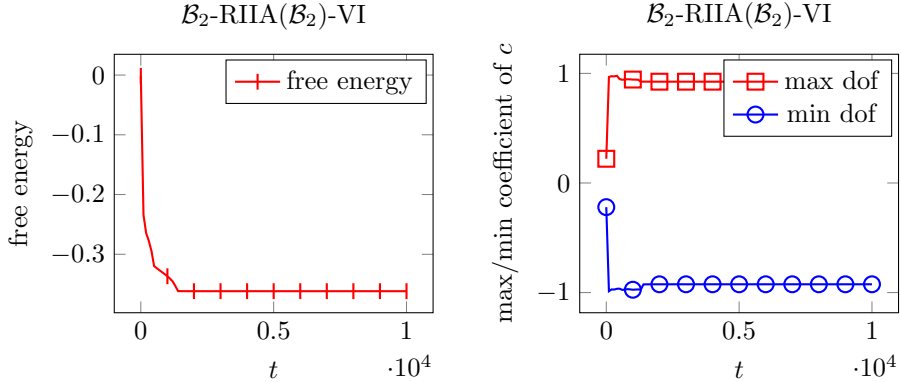
\begin{figure}[ht]

\begin{subfigure}{0.48\textwidth}
\centering
\begin{tikzpicture}
\centering
\begin{axis}[legend style={nodes={scale=1.0, transform shape}},  width=\textwidth, 
        title={$\mathcal{B}_2$-RIIA($\mathcal{B}_2$)-VI},
	ylabel near ticks,
	xlabel=$t$,
    	ylabel={free energy}]
    	
    \addplot[red,mark=|, mark size=3pt, mark repeat=10,  each nth point=100, filter discard warning=false, thick] table[y=energy, col sep=comma] {CH6P_data.csv};
    \addlegendentry {free energy};
\end{axis}
\end{tikzpicture}
\end{subfigure}\hspace{0.04\textwidth}%
\begin{subfigure}{0.48\textwidth}
\centering
\begin{tikzpicture}
\centering
\begin{axis}[legend style={nodes={scale=1.0, transform shape}},  width=\textwidth, 
        title={$\mathcal{B}_2$-RIIA($\mathcal{B}_2$)-VI},
	ylabel near ticks,
	xlabel=$t$,
    	ylabel={max/min coefficient of $c$}]
    	
    \addplot[red,mark=square, mark size=3pt, mark repeat=10,  each nth point=100, filter discard warning=false,  thick] table[y=max_dof, col sep=comma] {CH6P_data.csv};
    \addlegendentry {max dof};
    
    \addplot[blue,mark=o, mark size=3pt, mark repeat=10,  each nth point=100, filter discard warning=false,  thick] table[y=min_dof, col sep=comma] {CH6P_data.csv};
    \addlegendentry {min dof};
\end{axis}

\end{tikzpicture}
\end{subfigure}%
\caption{The free energy (left) and maximum and minimum coefficients (right) of the numerical solution to~\eqref{eq:cahn_hilliard_general} subject to periodic boundary conditions and the initial concentration~\eqref{eq:CH_ex6_init}.}
\label{fig:CH_ex6_per_data}
\end{figure}

With periodic boundary conditions, we obtain similar results. Snapshots are shown in Figure~\ref{fig:CH_ex6_per_snaps} and the free energy and maximum/minimum coefficients are shown in Figure~\ref{fig:CH_ex6_per_data}

Next, we allow the system to evolve given the initial condition
\begin{equation}\label{eq:CH_ex4_init}
c(0) = \frac{1}{10}(2\cdot\text{rand}(x, y) - 1),
\end{equation}
where $\text{rand}(x, y)$ is a random number between $0$ and $1$. 

\begin{figure}[ht]
\centering
\begin{subfigure}{0.25\textwidth}
\centering
{\includegraphics[trim={10cm 3cm 19cm 2cm},clip,  width=\textwidth] {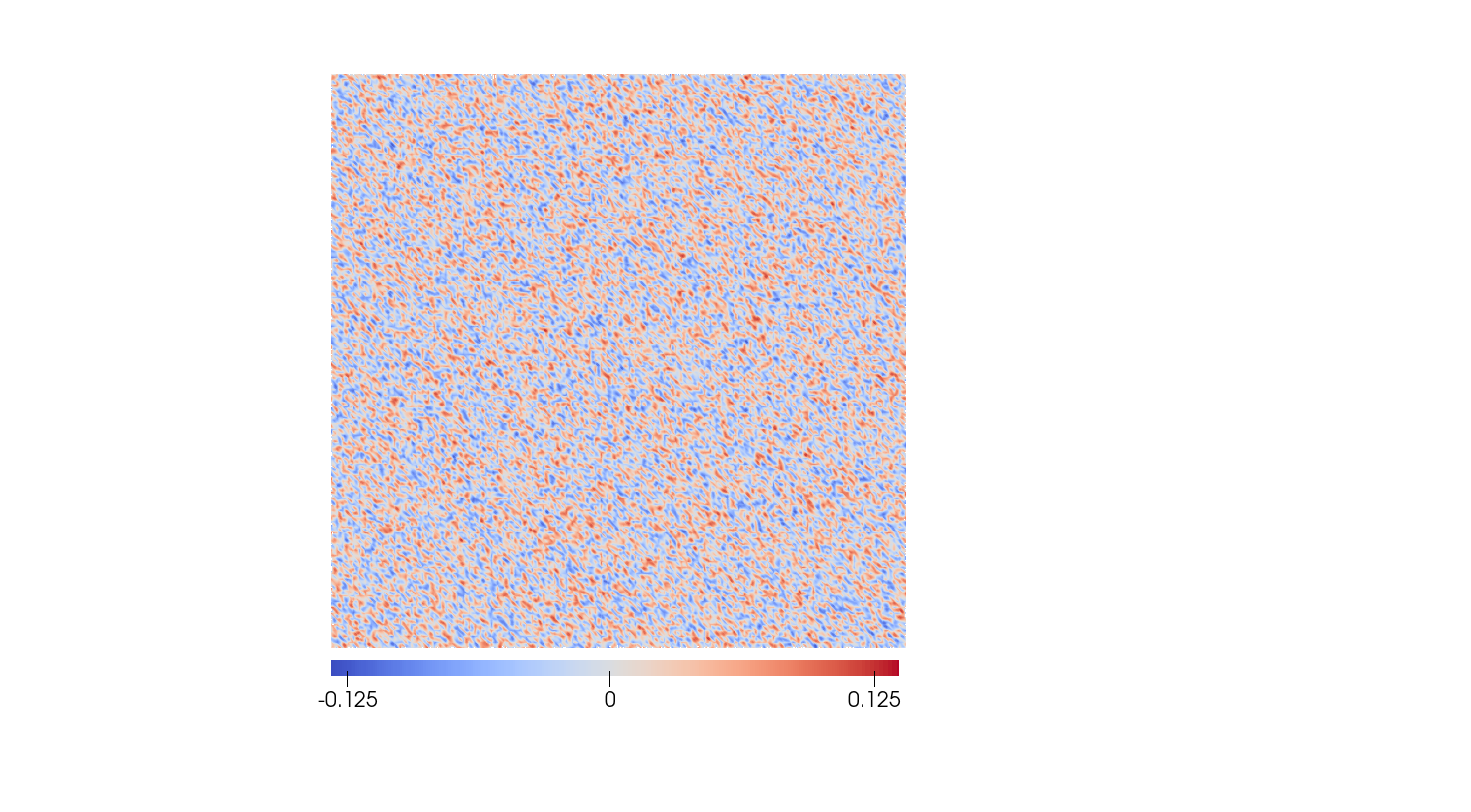}}
\subcaption{$t = 0.0$}
\end{subfigure}%
\begin{subfigure}{0.25\textwidth}
\centering
{\includegraphics[trim={10cm 3cm 19cm 2cm},clip,  width=\textwidth] {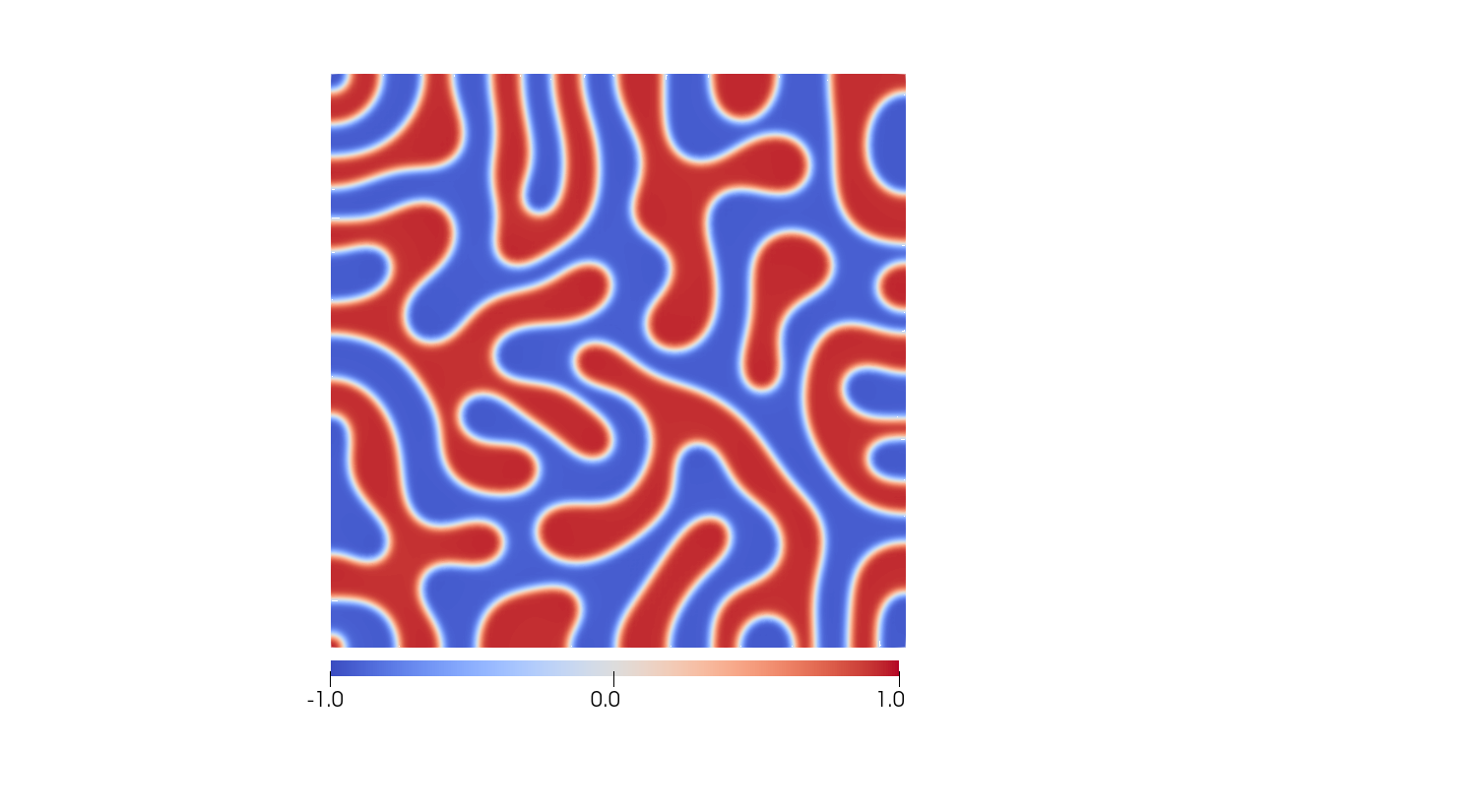}}
\subcaption{$t = 0.01$}
\end{subfigure}%
\begin{subfigure}{0.25\textwidth}
\centering
{\includegraphics[trim={10cm 3cm 19cm 2cm},clip,  width=\textwidth] {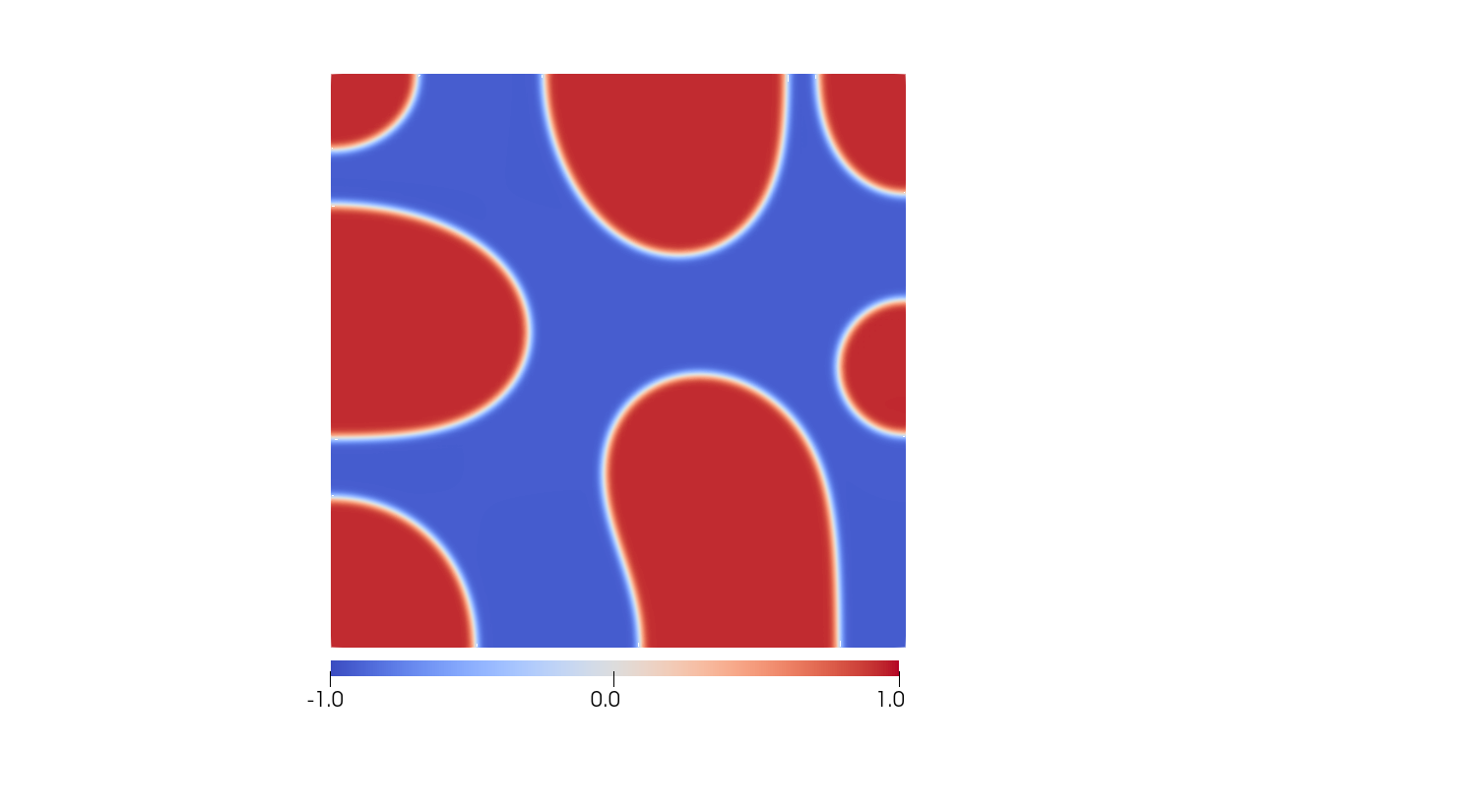}}
\subcaption{$t = 0.5$}
\end{subfigure}%
\begin{subfigure}{0.25\textwidth}
\centering
{\includegraphics[trim={10cm 3cm 19cm 2cm},clip,  width=\textwidth] {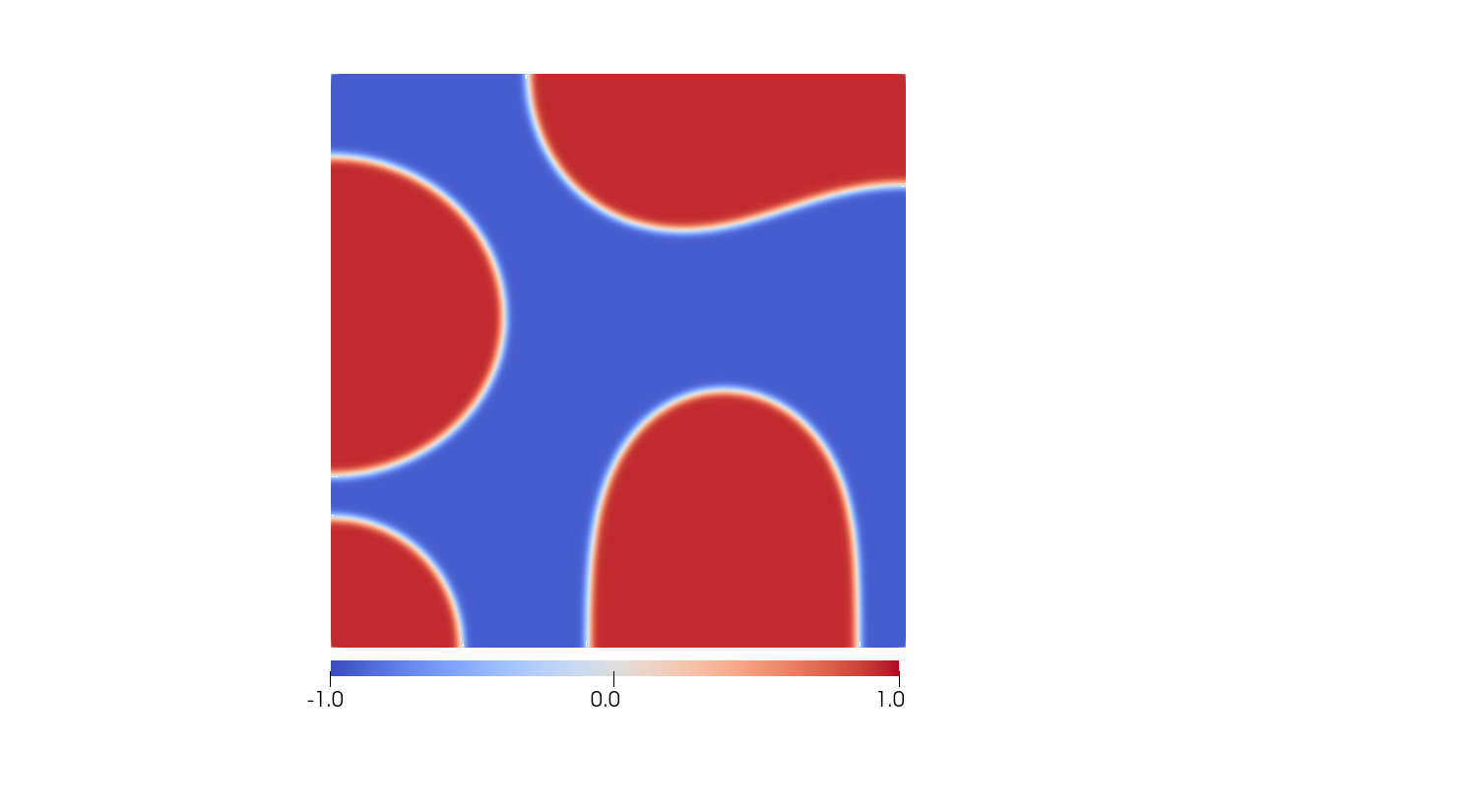}}
\subcaption{$t = 1.0$}
\end{subfigure}%
\caption{Snapshots of the numerical solution to \eqref{eq:cahn_hilliard_general} subject to homogeneous Neumann boundary conditions and the initial concentration~\eqref{eq:CH_ex4_init} at $t =0, 0.01, 0.5, 1.0$.}
\label{fig:CH_ex4_neumann_snaps}
\end{figure}

\begin{figure}[ht]

\begin{subfigure}{0.48\textwidth}
\centering
\begin{tikzpicture}
\centering
\begin{axis}[legend style={nodes={scale=1.0, transform shape}},  width=\textwidth, 
        title={$\mathcal{B}_2$-RIIA($\mathcal{B}_2$)-VI},
	ylabel near ticks,
	xlabel=$t$,
    	ylabel={free energy}]
    
    \addplot[red,mark=|, mark size=3pt, mark repeat=10,  each nth point=100, filter discard warning=false, thick] table[y=energy, col sep=comma] {CH4N_data.csv};
    \addlegendentry {free energy};
\end{axis}
\end{tikzpicture}
\end{subfigure}\hspace{0.04\textwidth}%
\begin{subfigure}{0.48\textwidth}
\centering
\begin{tikzpicture}
\centering
\begin{axis}[legend style={nodes={scale=1.0, transform shape}},  width=\textwidth, 
        title={$\mathcal{B}_2$-RIIA($\mathcal{B}_2$)-VI},
	ylabel near ticks,
	xlabel=$t$,
    	ylabel={max/min coefficient of $c$}]
    
    \addplot[red,mark=square, mark size=3pt, mark repeat=10,  each nth point=100, filter discard warning=false,  thick] table[y=max_dof, col sep=comma] {CH4N_data.csv};
    \addlegendentry {max dof};
    
    \addplot[blue,mark=o, mark size=3pt, mark repeat=10,  each nth point=100, filter discard warning=false,  thick] table[y=min_dof, col sep=comma] {CH4N_data.csv};
    \addlegendentry {min dof};
\end{axis}
\end{tikzpicture}
\end{subfigure}%
\caption{The free energy (left) and maximum and minimum coefficients (right) of the numerical solution to~\eqref{eq:cahn_hilliard_general} subject to homogeneous Neumann boundary conditions and the initial concentration~\eqref{eq:CH_ex4_init}.}
\label{fig:CH_ex4_neumann_data}
\end{figure}
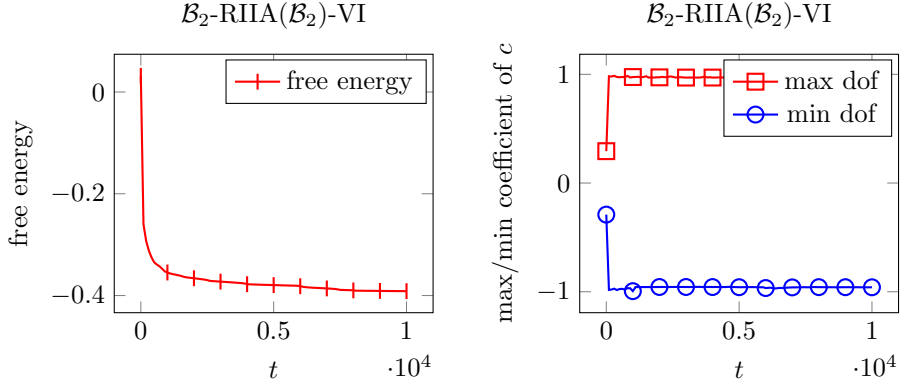

\begin{figure}[ht]
\centering
\begin{subfigure}{0.25\textwidth}
\centering
{\includegraphics[trim={10cm 3cm 19cm 2cm},clip,  width=\textwidth] {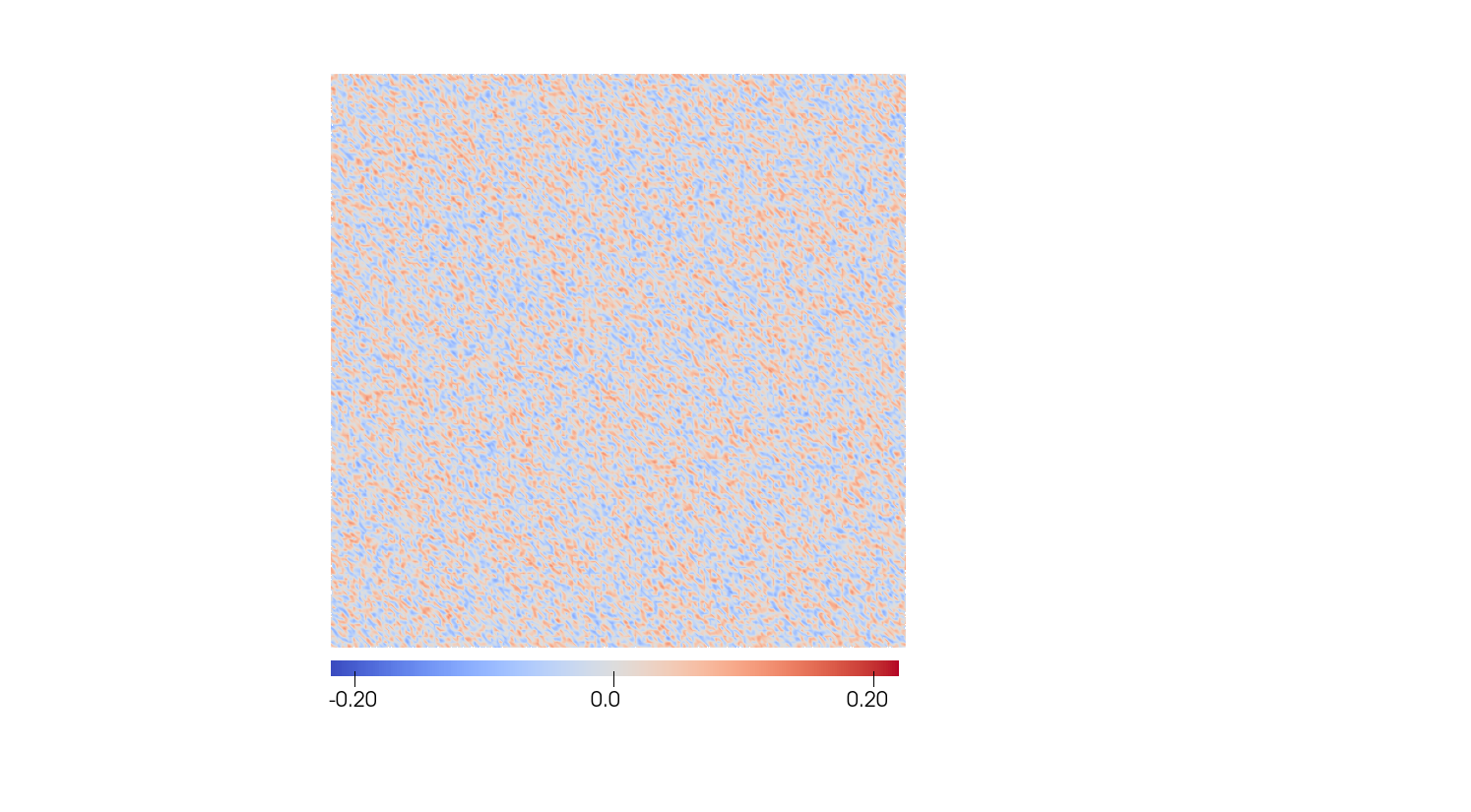}}
\subcaption{$t = 0.0$}
\end{subfigure}%
\begin{subfigure}{0.25\textwidth}
\centering
{\includegraphics[trim={10cm 3cm 19cm 2cm},clip,  width=\textwidth] {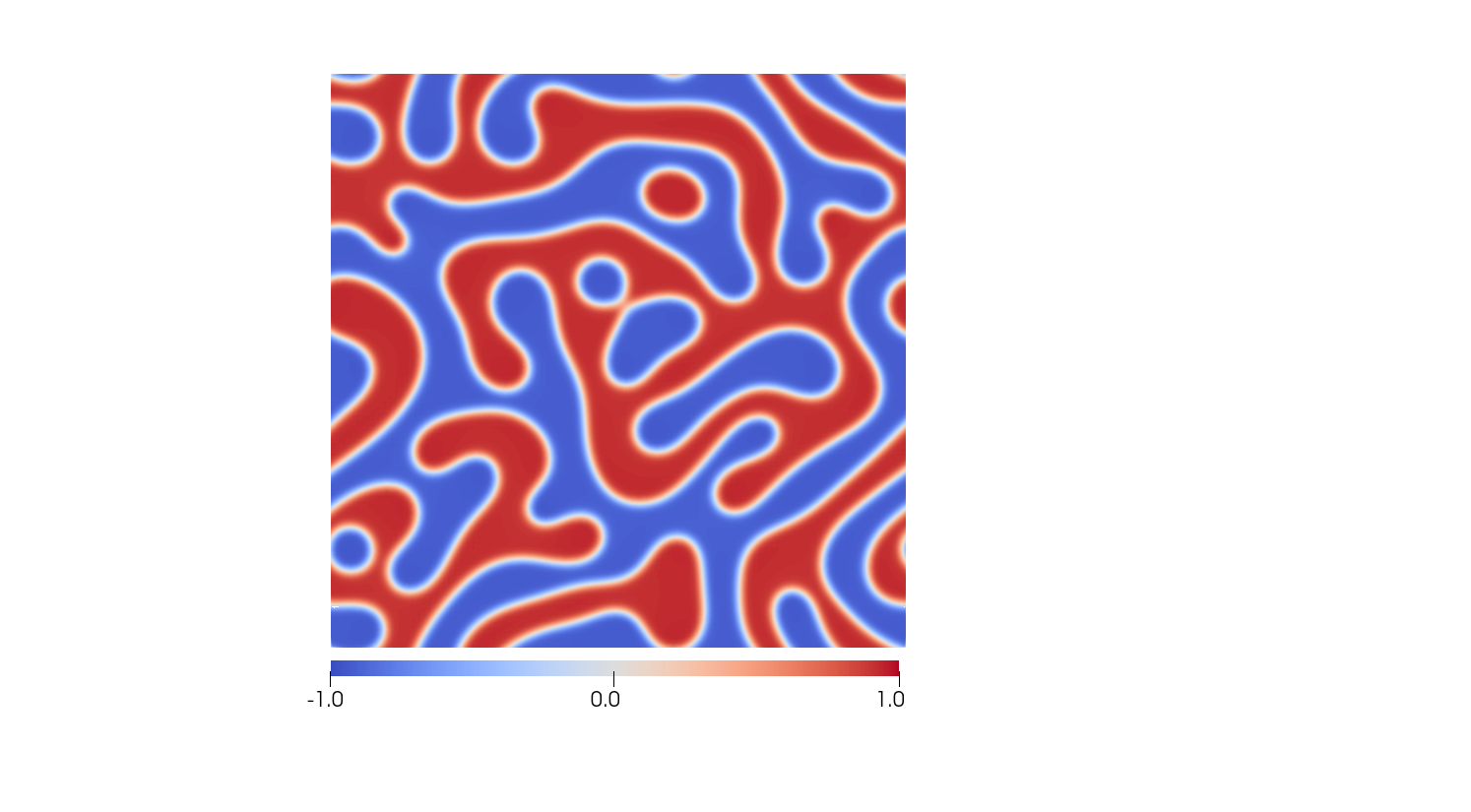}}
\subcaption{$t = 0.01$}
\end{subfigure}%
\begin{subfigure}{0.25\textwidth}
\centering
{\includegraphics[trim={10cm 3cm 19cm 2cm},clip, width=\textwidth] {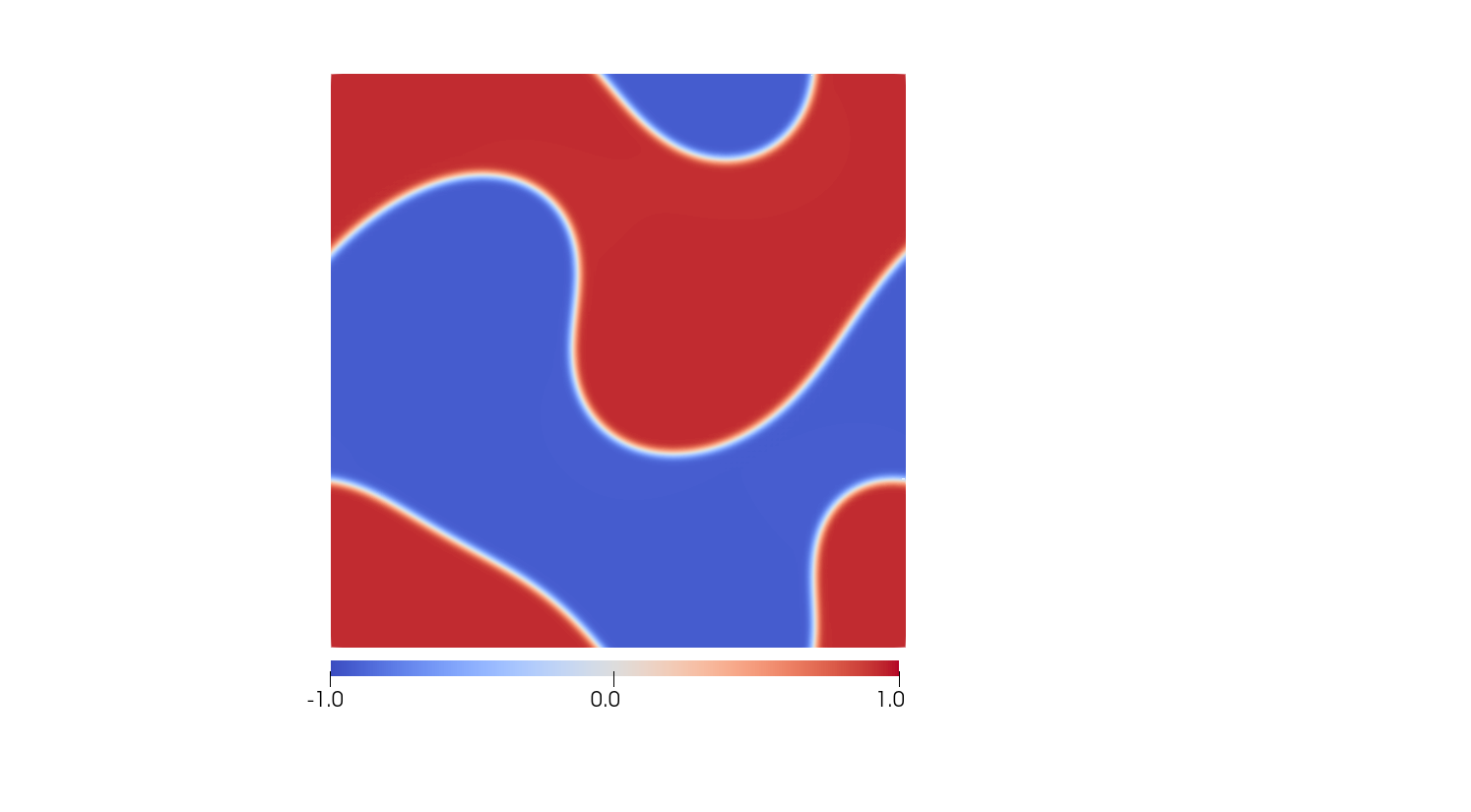}}
\subcaption{$t = 0.5$}
\end{subfigure}%
\begin{subfigure}{0.25\textwidth}
\centering
{\includegraphics[trim={10cm 3cm 19cm 2cm},clip, width=\textwidth] {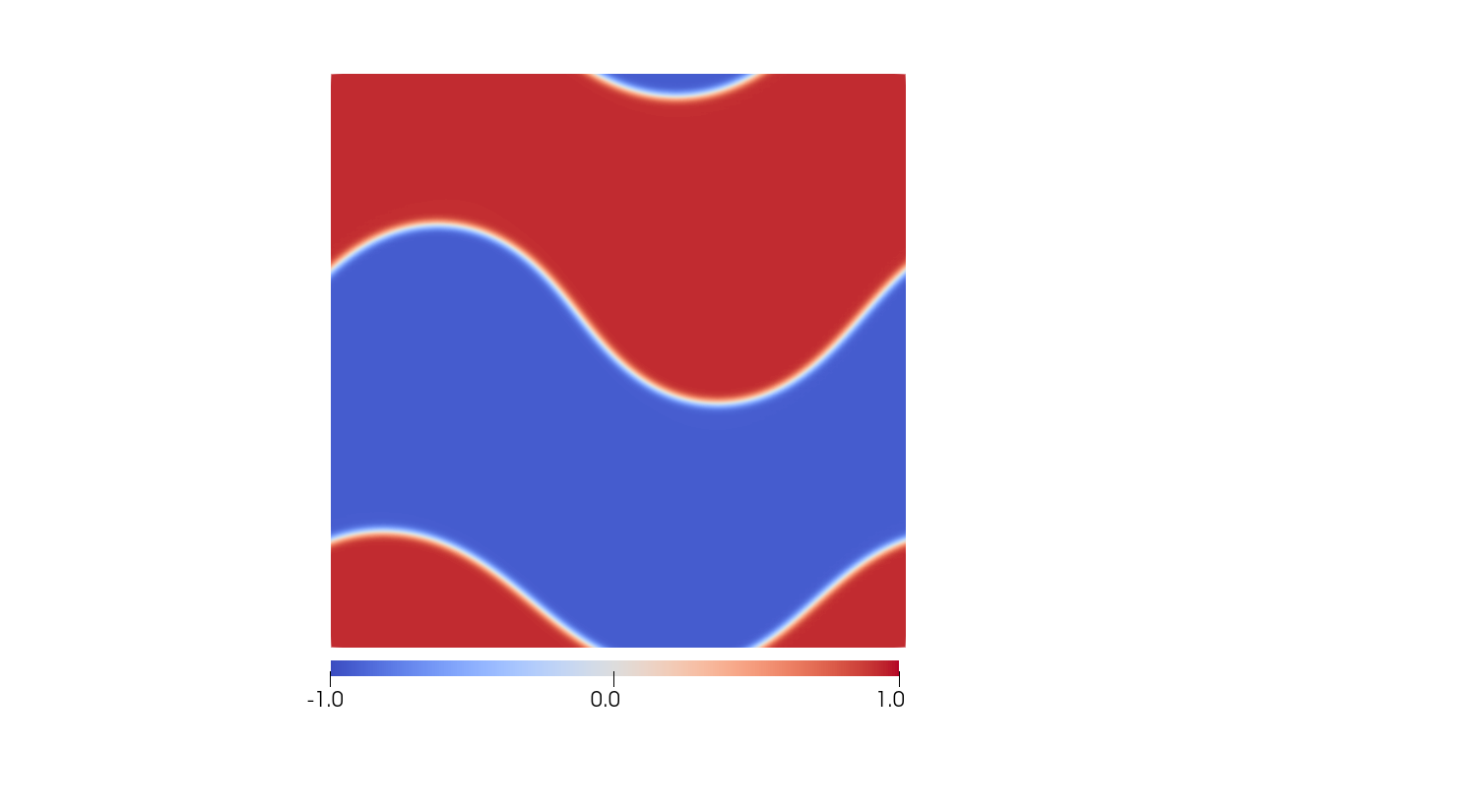}}
\subcaption{$t = 1.0$}
\end{subfigure}%
\caption{Snapshots of the numerical solution to \eqref{eq:cahn_hilliard_general} subject to periodic boundary conditions and the initial concentration~\eqref{eq:CH_ex4_init} at $t =0, 0.01, 0.5, 1.0$.}
\label{fig:CH_ex4_per_snaps}
\end{figure}

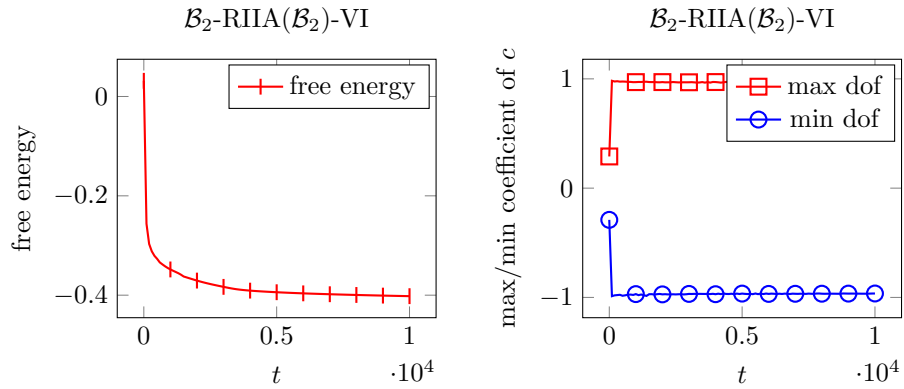
\begin{figure}[ht]

\begin{subfigure}{0.48\textwidth}
\centering
\begin{tikzpicture}
\centering
\begin{axis}[legend style={nodes={scale=1.0, transform shape}},  width=\textwidth, 
        title={$\mathcal{B}_2$-RIIA($\mathcal{B}_2$)-VI},
	ylabel near ticks,
	xlabel=$t$,
    	ylabel={free energy}]
    
    \addplot[red,mark=|, mark size=3pt, mark repeat=10,  each nth point=100, filter discard warning=false, thick] table[y=energy, col sep=comma] {CH4P_data.csv};
    \addlegendentry {free energy};
\end{axis}
\end{tikzpicture}
\end{subfigure}\hspace{0.04\textwidth}%
\begin{subfigure}{0.48\textwidth}
\centering
\begin{tikzpicture}
\centering
\begin{axis}[legend style={nodes={scale=1.0, transform shape}},  width=\textwidth, 
        title={$\mathcal{B}_2$-RIIA($\mathcal{B}_2$)-VI},
	ylabel near ticks,
	xlabel=$t$,
    	ylabel={max/min coefficient of $c$}]
    
    \addplot[red,mark=square, mark size=3pt, mark repeat=10,  each nth point=100, filter discard warning=false,  thick] table[y=max_dof, col sep=comma] {CH4P_data.csv};
    \addlegendentry {max dof};
    
    \addplot[blue,mark=o, mark size=3pt, mark repeat=10,  each nth point=100, filter discard warning=false,  thick] table[y=min_dof, col sep=comma] {CH4P_data.csv};
    \addlegendentry {min dof};
\end{axis}
\end{tikzpicture}
\end{subfigure}%
\caption{The free energy (left) and maximum and minimum coefficients (right) of the numerical solution to~\eqref{eq:cahn_hilliard_general} subject to periodic boundary conditions and the initial concentration~\eqref{eq:CH_ex4_init}.}
\label{fig:CH_ex4_per_data}
\end{figure}

Snapshots of the spinodal decomposition, subject to homogeneous Neumann boundary conditions, are shown in Figure~\ref{fig:CH_ex4_neumann_snaps}, and the free energy and maximum/minimum coefficients are shown in Figure~\ref{fig:CH_ex4_neumann_data}. The results of the same experiment (albeit with different seeding of the random variable) repeated with periodic boundary conditions are shown in show in Figure~\ref{fig:CH_ex4_per_snaps} and Figure~\ref{fig:CH_ex4_per_data}.

\subsubsection{Solver performance}

For each of the above examples, the number of nonlinear iterations is modest in light of the number of steps taken (see Table~\ref{tab:cahn_hilliard_stats}), with only slightly more than one nonlinear iteration required per time step.
\begin{table}[ht]
\begin{center}
\begin{tabular}{ c  c  c  c  c }
 initial condition & boundary conditions & \# steps & time (s) & \# nonlinear its.  \\
 \hline
\eqref{eq:CH_ex4_init} & Homogeneous Neumann & 10,000 & $1.31\cdot 10^{5}$ & 10,076 \\  
\eqref{eq:CH_ex6_init} & Homogeneous Neumann & 10,000 & $1.29\cdot 10^{5}$ & 10,077 \\  
\eqref{eq:CH_ex4_init} & Periodic & 10,000 & $3.84\cdot 10^{5}$ & 10,044 \\  
\eqref{eq:CH_ex6_init} & Periodic & 10,000 & $3.84\cdot 10^{5}$ & 10,045 \\  
\end{tabular}
\end{center}
\caption{Solver statistics for approximating the solution to~\eqref{eq:cahn_hilliard_general} to final time $t = 1.0$ subject to the given initial and boundary conditions.}
\label{tab:cahn_hilliard_stats}
\end{table}

\section{Conclusions and future work}
\label{sec:conc}
We have developed an approach to enforcing bounds constraints in evolution equations by combining variational inequalities with collocation-type implicit Runge--Kutta methods.
These methods support various modes of bounds constraints, and by working in terms of the Bernstein basis, we can obtain methods that uniformly enforce bounds constraints in both space and time.
Our numerical results indicate high-quality solutions that enforce bounds while yielding higher-order approximations.

These results are quite promising, but leave open several important questions.
First of all, fast and scalable solvers will be critical to the long-term success of the method.  Besides extending the Firedrake/PETSc interface~\cite{kirby2018solver} to work correctly with the reduced spaces varying across iterations,  it could be possible to extend the FASCd method~\cite{bueler2024full} to higher order spatial discretizations.  Furthermore, we lack convergence theory for this method.  It could be possible to extend the analysis in~\cite{barrenechea2024nodally, kirby2024high} to obtain estimates for the stage equations at each time step, and then this would have to be combined with a Runge--Kutta error analysis to get fully discrete estimates.